\documentclass{amsart}

\usepackage{amsmath,amsthm,amssymb,latexsym,amscd}
\usepackage{graphicx}

\newcommand{\im}{\mathrm{ Im \;}}

\newcommand{\mmod}{\ \text{mod}\ }

\newcommand{\Z}{\mathbb{Z}}

\newcommand{\cal}{\mathcal}
\newtheorem{thm}{Theorem}
\newtheorem{lem}[thm]{Lemma}
\newtheorem{coro}[thm]{Corollary}
\newtheorem{propo}[thm]{Proposition}
\newtheorem{defi}[thm]{Definition}

\newtheorem{nota}[thm]{Notation}

\newcommand{\tet}[6]{\;Tet
  \left[ 
  \begin{array}{ccc}
 	#1 & #2 & #5 \\ 
 	#3 & #4 & #6  
 	\end{array}
 	\right]}
\newcommand{\admpwise}[2]{\left\{ \begin{array}{ll}
 	#1 & \hbox{, if } 
#2 \hbox{ admissible}
\\ 
  0  & \hbox{, otherwise.}
 	\end{array} \right}

\begin{document}

\title[skein module of the quaternionic manifold]{On the Kauffman bracket skein module\\ of the quaternionic manifold}

\author{Patrick M. Gilmer}
\address{Mathematics Department\\
Louisiana State University\\
Baton Rouge, Louisiana}
\email{gilmer@math.lsu.edu}
\thanks{partially supported by NSF-DMS-0203486}

\author{John M. Harris}
\address{The University of Southern Mississippi\\
Long Beach, Mississippi}
\email{john.m.harris@usm.edu}

\date{December 14, 2004}

\begin{abstract} 
We use recoupling theory to study the Kauffman bracket skein module of the quaternionic manifold over $\Z[A^{\pm 1}]$ localized by inverting all the cyclotomic polynomials.  We prove that the skein module is spanned by five elements.  Using the quantum invariants of these skein elements and the $\Z_2$-homology of the manifold, we determine that they are linearly independent.
\end{abstract}

\maketitle

\section{Introduction}
\label{sec:Introduction}

In \cite{Ka88}, Kauffman presents an elegant construction of the Jones polynomial, an invariant of oriented links in $S^3$, by constructing a new invariant, the Kauffman bracket polynomial.  The Kauffman bracket is an invariant of unoriented framed links in $S^3$, defined by the following skein relations:

\begin{enumerate}
	\item $\left<\begin{minipage}{0.5in}\includegraphics[width=0.5in]{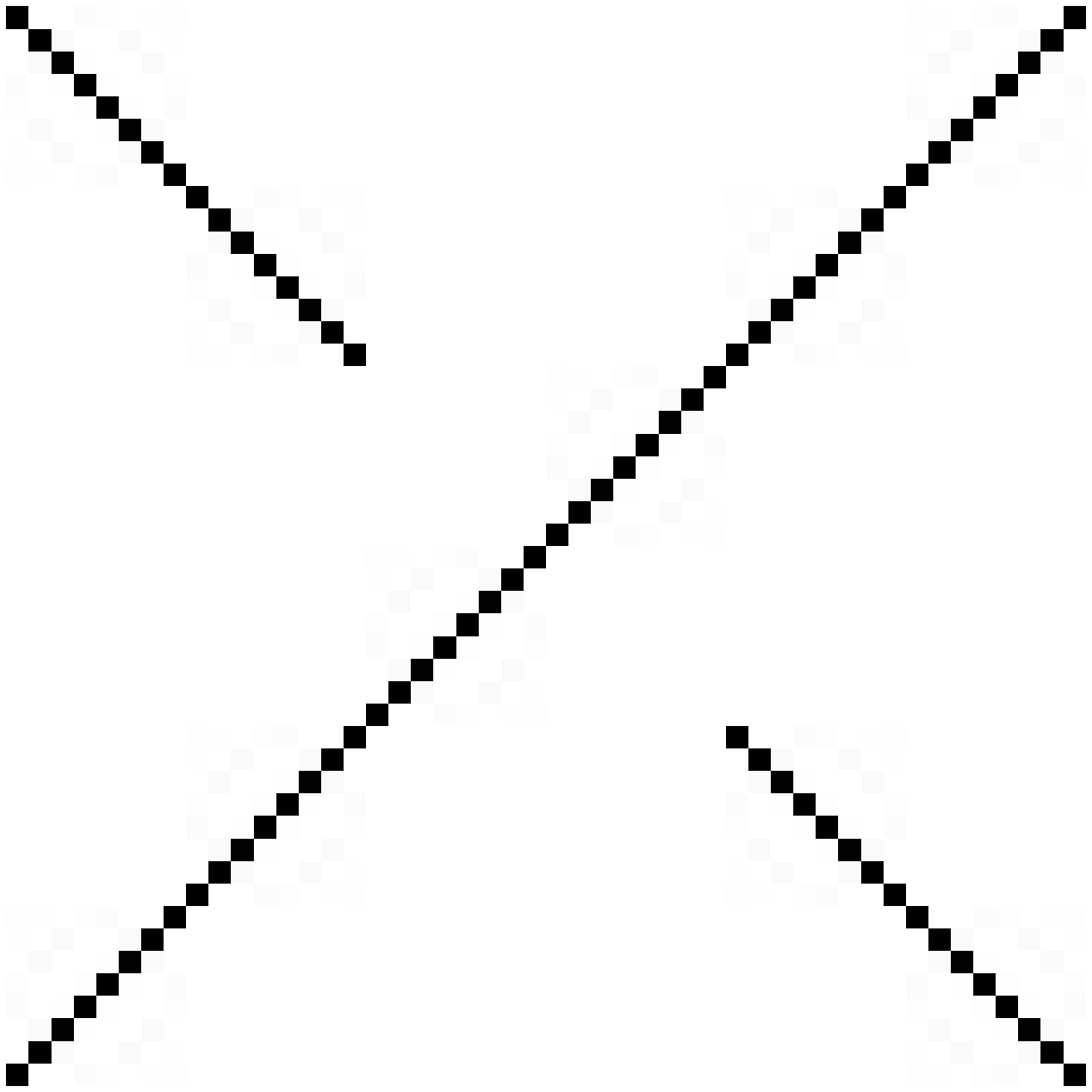}\end{minipage}\right> = A \left<\begin{minipage}{0.5in}\includegraphics[width=0.5in]{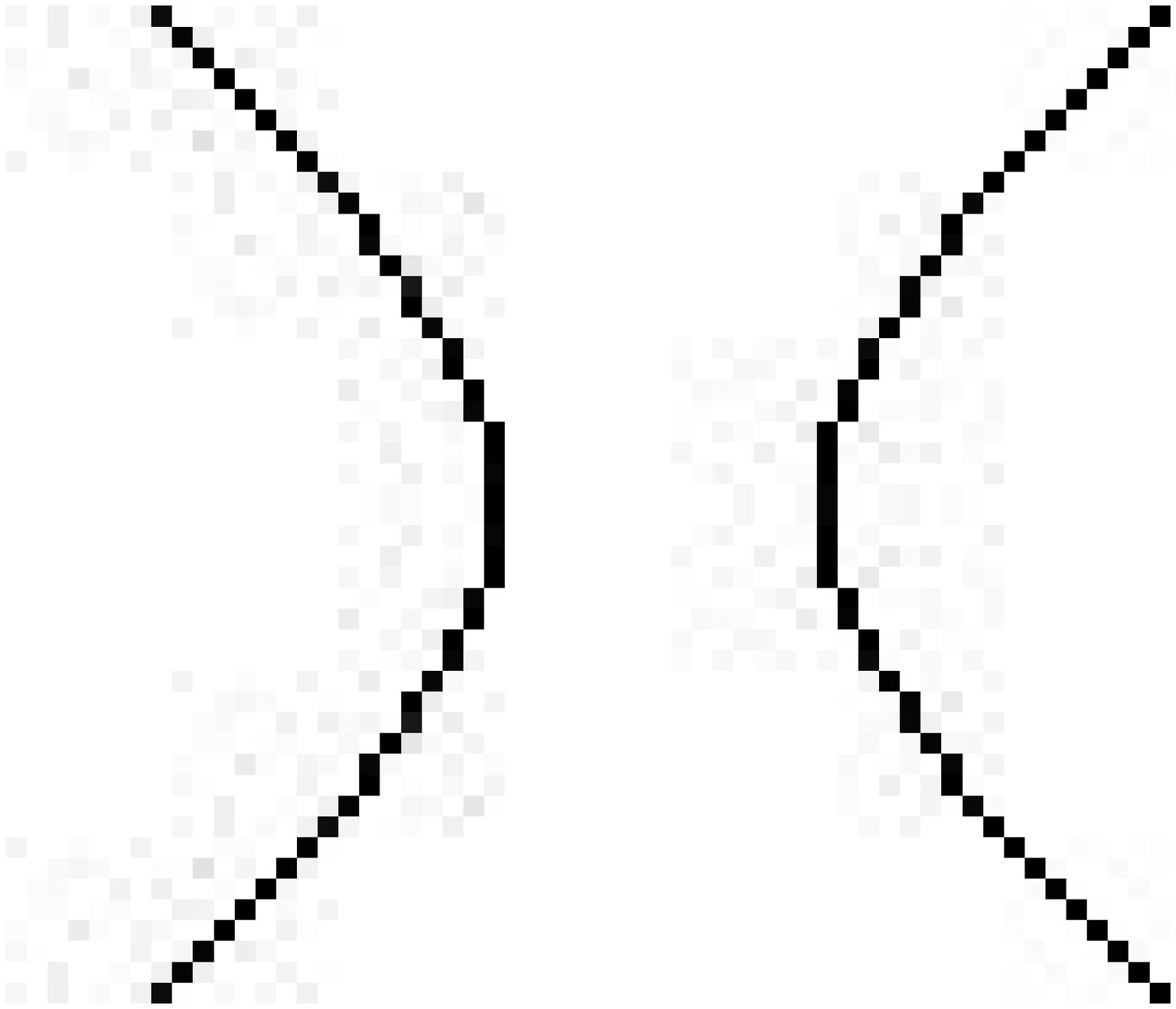}\end{minipage}\right> + A^{-1} \left<\begin{minipage}{0.5in}\includegraphics[width=0.5in]{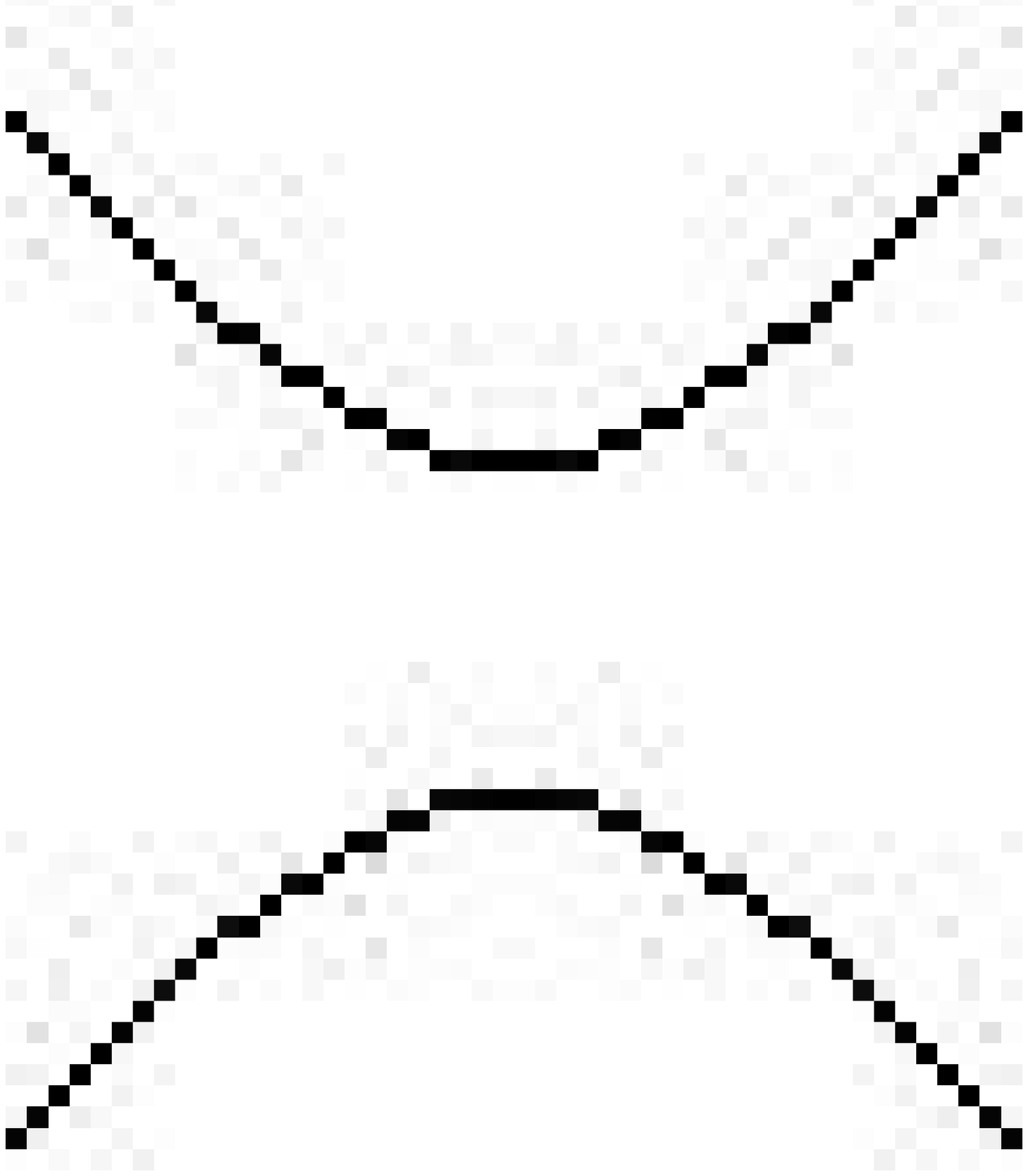}\end{minipage}\right>$
	\item $<L \cup \hbox{unknot}> = (-A^{-2} - A^2) <L>$
\end{enumerate}
 
For the invariant to be well-defined, one also must normalize it by choosing a value for the empty link.  $<\hbox{empty link}> = 1$, for instance.

Alternatively, we can use the skein relations to construct a module of equivalence classes of links in $S^3$, or, for that matter, in any oriented 3-manifold.  See Przytycki (\cite{Pr91}) and Turaev (\cite{Tu88}).

\begin{defi}
Let $N$ be an oriented 3-manifold, and let $R$ be a commutative ring with identity, with a specified unit $A$.  The Kauffman bracket skein module of $N$, denoted $S(N;R,A)$, or simply $S(N)$, is the free $R$-module generated by the framed isotopy classes of unoriented links in $N$, including the empty link, quotiented by the skein relations which define the Kauffman bracket.
\end{defi}

Since every crossing and unknot can be eliminated from a link in $S^3$ by the skein relations, $S(S^3)$ is generated by the empty link. Kauffman's argument that his bracket polynomial is well-defined shows that $S(S^3)$ is free on the empty link.

Hoste and Przytycki have, in fact, computed the skein modules of all of the closed, oriented manifolds of genus 1:  $S(L(p,q))$ , which is free on $\left\lfloor\frac{p}{2}\right\rfloor + 1$ generators (\cite{HP93}), and $S(S^1 \times S^2) \cong \Z[A^{\pm 1}]\oplus
(\bigoplus^\infty_{i=1}\Z[A^{\pm 1}]/(1-A^{2i+4}))$, (\cite{HP95}). 

They have also computed the skein modules of $I$-bundles over surfaces (\cite{HP93}).

Additionally, Bullock has found a presentation for the complement of a $(2, q)$ torus knot in \cite{Bu95}, and has determined whether or not the 
 skein module of the result of integral surgery on a trefoil is finitely generated in \cite{Bu97}. Bullock and LoFaro have computed the skein module of the exteriors of the twist knots in \cite{BuLo}.

All of the computations mentioned above have been carried out with $R=\Z[A^{\pm1}]$.

Let $M$ be the quaternionic manifold, $S^3$ quotiented by the action of the quaternionic group $Q_8$:  $M$ can be obtained by identifying opposite faces of a cube with one-quarter twists.  

\begin{figure}
$$
\begin{minipage}{1.5in}\includegraphics[width=1.5in]{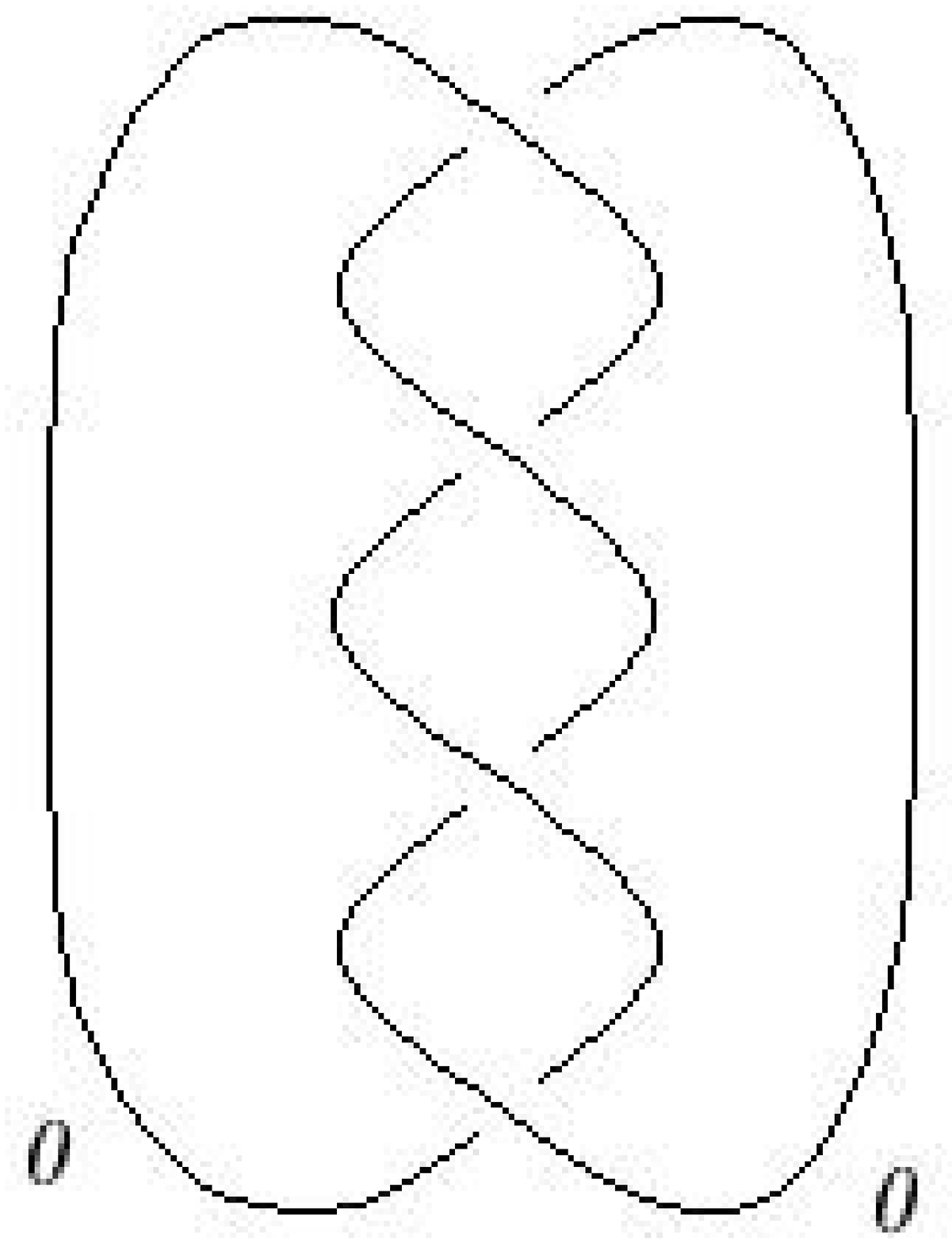}\end{minipage}
\hspace{1.5in}
\begin{minipage}{1.5in}\includegraphics[width=1.5in]{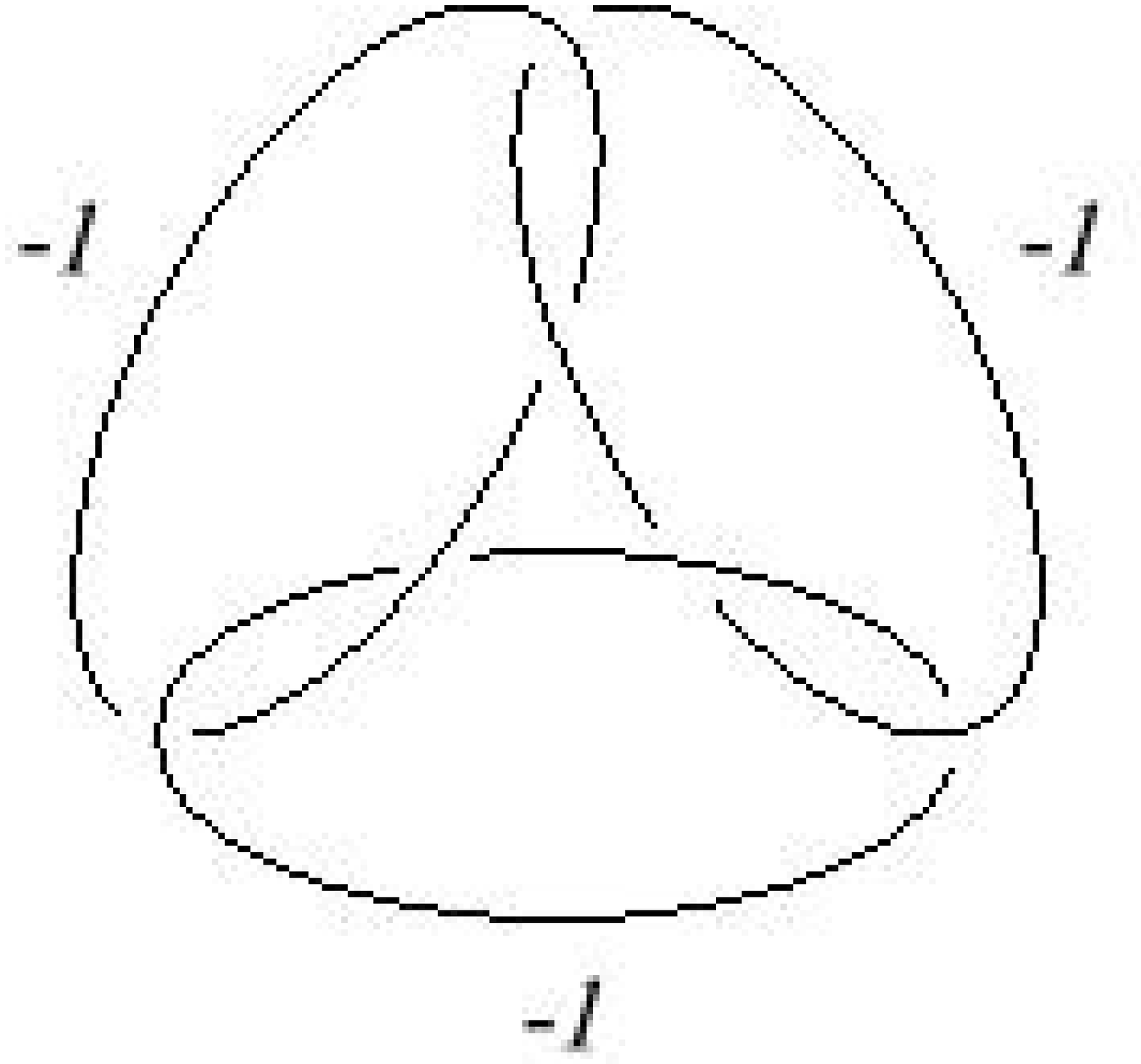}\end{minipage}
$$
\caption{Surgery descriptions of the quaternionic manifold}
\label{fig:quatsurgery}
\end{figure}

$M$ is the three-fold cover of $S^3$ branched over the trefoil, and is an irreducible manifold of genus 2.  Rolfsen gives surgery descriptions of this manifold in \cite{Ro76}.  See Figure \ref{fig:quatsurgery}.

Let ${\cal R}$ denote $\Z[A^{\pm1}]$, localized by inverting the multiplicative set generated by the elements of the set $\{ A^n-1 | n \in \Z^+ \} .$  Equivalently, we could say that we are localizing by inverting all cyclotomic polynomials in $A$ or all quantum integers (see below for
the
 definition of a quantum integer).  
\begin{thm}
$S(M;{\cal R},A)) \approx {\cal R}^5$.  
\end{thm}

\begin{coro}
$S(M;\Z[A^{\pm1}],A))$ has at most cyclotomic torsion, i.e. every torsion element is annihilated by a product of cyclotomic polynomials.  
\end{coro}

In fact, we conjecture that $S(M;\Z[A^{\pm1}],A))$ is torsion-free.
The quaternionic manifold is the first closed, irreducible, genus-two 3-manifold whose Kauffman bracket skein module (over a localization of  $\Z[A^{\pm1}]$) has been computed.

From now on, we consider all skeins to be over ${\cal R}$ and denote $S(N;{\cal R},A)$ by $S(N)$.

Before we proceed, we note that the defining relations of the Kauffman bracket skein module respect $\Z_2$-homology.  Since $H_1(M;\Z_2) = \Z_2 \oplus \Z_2$, $S(M)$ is a direct sum of four submodules $S_1(M)$, $S_i(M)$, $S_j(M)$, and $S_k(M)$.  

Furthermore, the permutation group on three elements acts on the manifold permuting the three non-trivial homology classes.  Thus $S_i(M)$, $S_j(M)$, and $S_k(M)$ are isomorphic.  We will see that $S_1(M)$ is isomorphic to $ {\cal R}^2,$ and that $S_i(M)$, $S_j(M)$, and $S_k(M)$ are isomorphic to
${\cal R}.$

\section{Preliminary Remarks}
\label{sec:PreliminaryRemarks}
We draw $H_2$ as a ball before the handles are attached from top to bottom along the dashed lines.  Equivalently, one can view the diagrams as pictures of handlebodies whose handles have been cut.

\begin{figure}
$$
\begin{minipage}{2.4in}\includegraphics[width=2.4in]{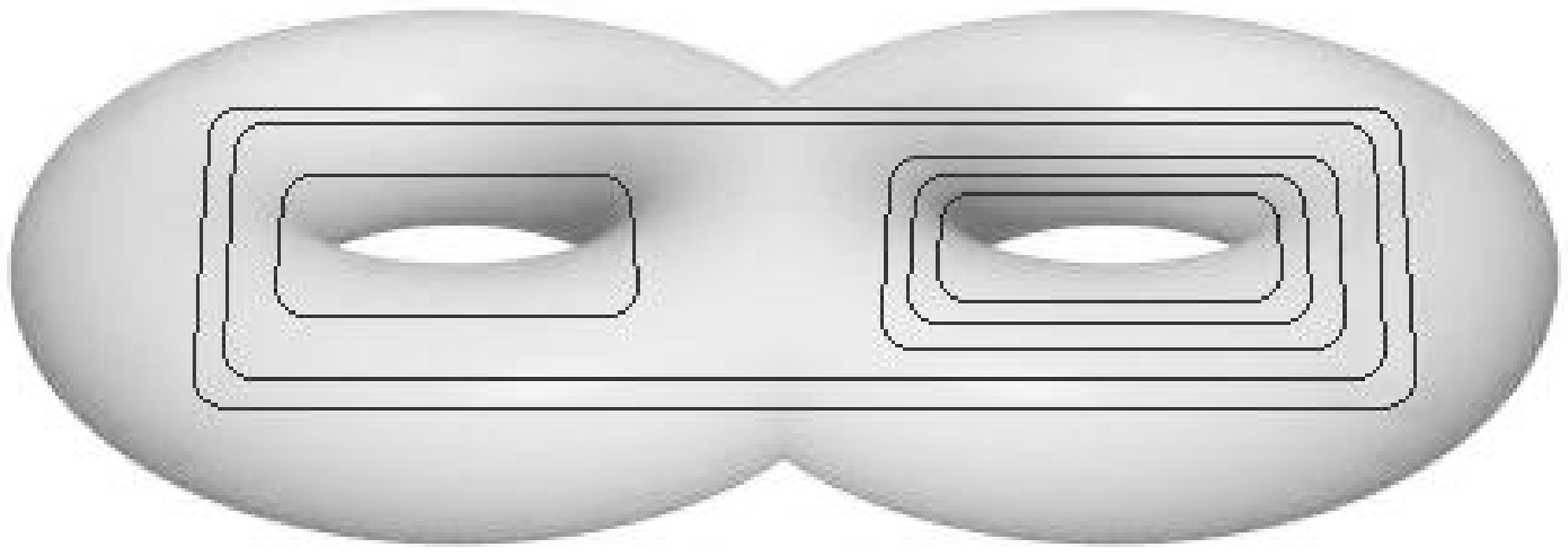}\end{minipage}
\hspace{0.5in} \begin{minipage}{2in}\includegraphics[width=2in]{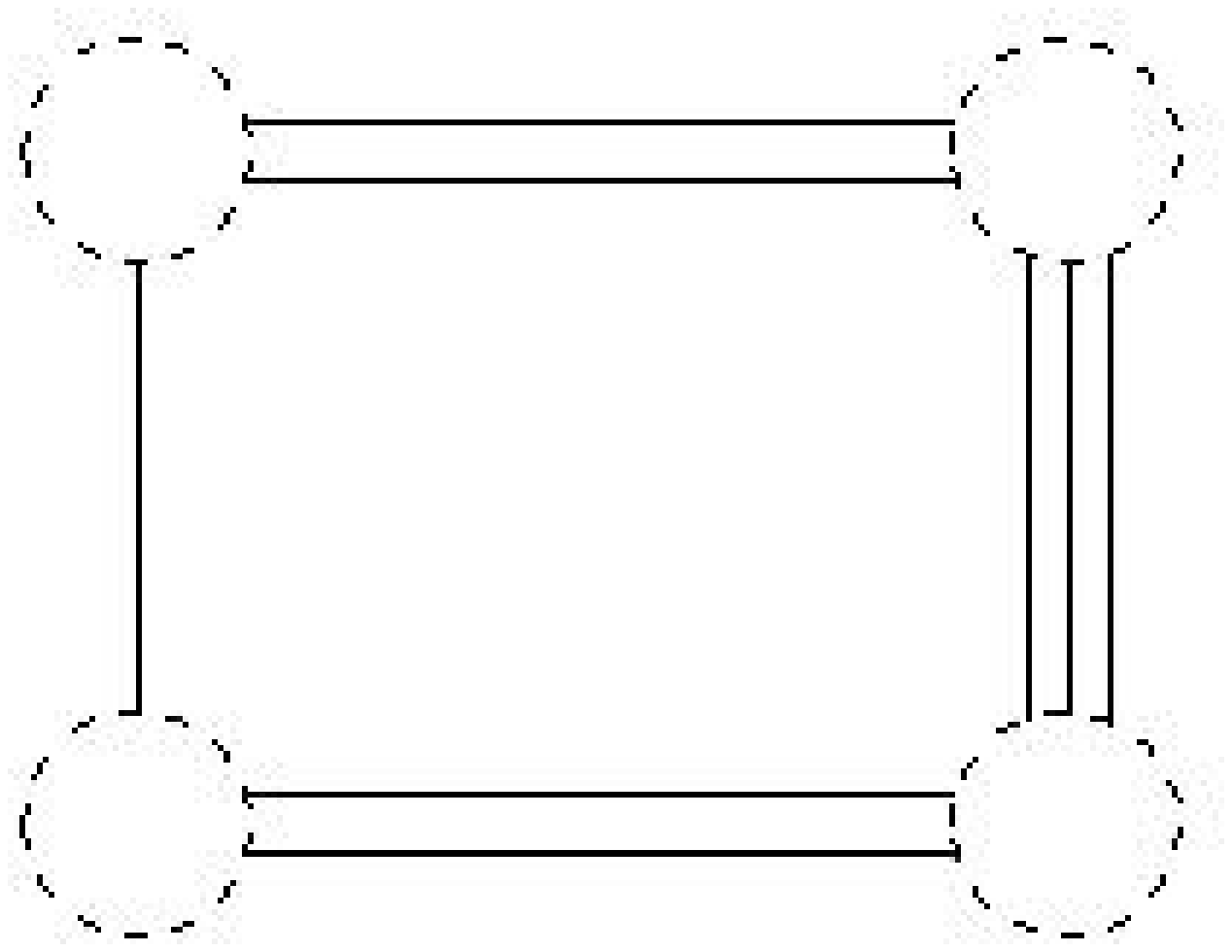}\end{minipage}$$
\caption{The link $x y^2 z^3$ in the standard diagram and the cut-open diagram of $H_2$}
\label{fig:basis}
\end{figure}

The standard basis $B$ for the Kauffman bracket skein module of the solid handlebody $H_g$ of genus g is the collection of all links without crossings in the standard diagram of $H_g$.  For $H_2$, we use Bullock's algebraic notation ${x^i y^j z^k}$ (\cite{Bu97}).  See Figure \ref{fig:basis}.  

\begin{defi}
If $|a-b| \leq c \leq a+b$ and $a+b+c \equiv 0 \mmod2$, then the triple $\{a,b,c\}$ is said to be admissible.
\end{defi}

\begin{figure}
$$\begin{minipage}{0.6in}\includegraphics[width=0.6in]{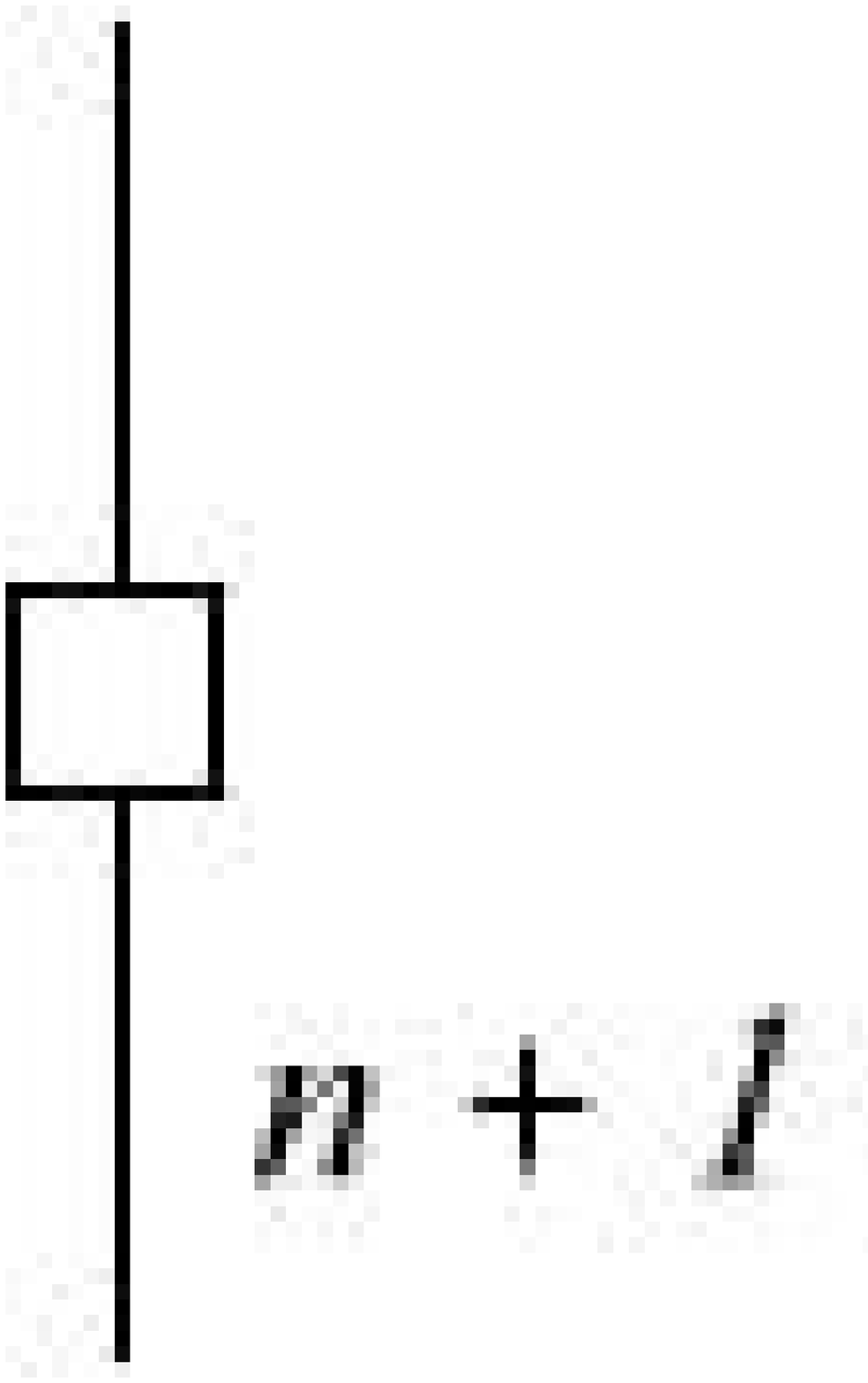}\end{minipage} = \begin{minipage}{0.6in}\includegraphics[width=0.6in]{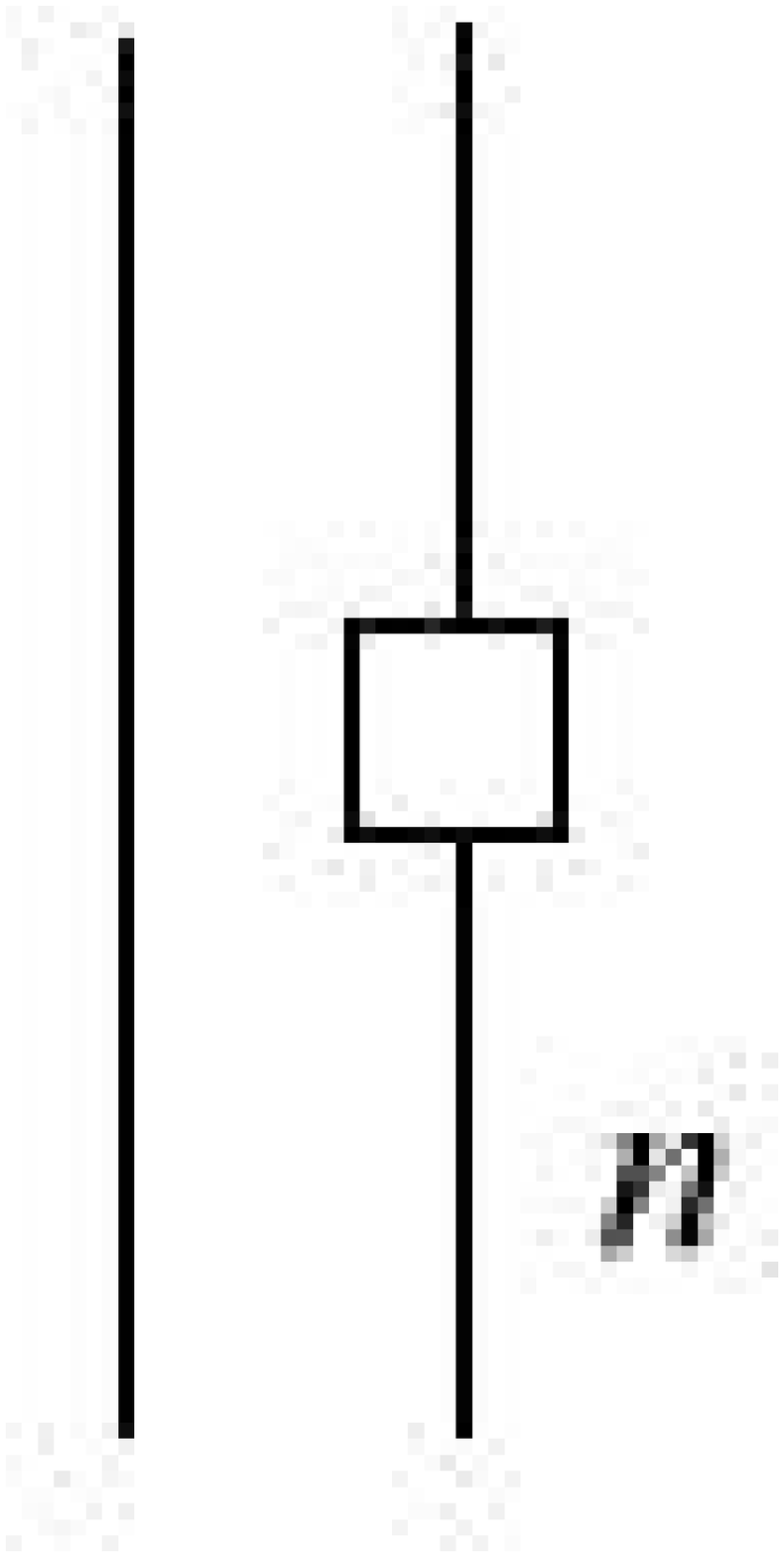}\end{minipage} - \frac{\Delta_{n-1}}{\Delta_n}  \begin{minipage}{0.6in}\includegraphics[width=0.6in]{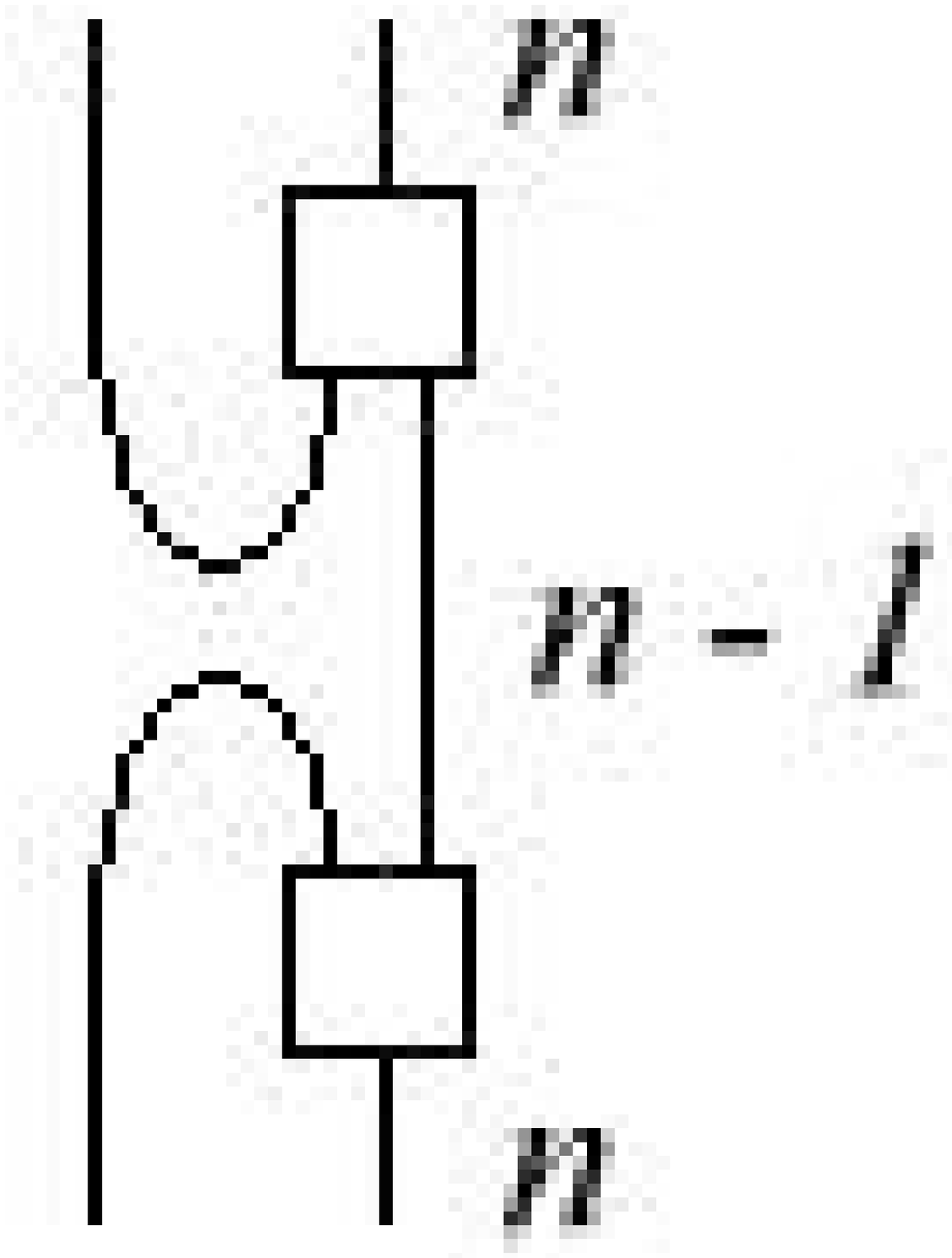}\end{minipage},$$ 
$$\hbox{where } \Delta_n = \begin{minipage}{0.8in}\includegraphics[width=0.8in]{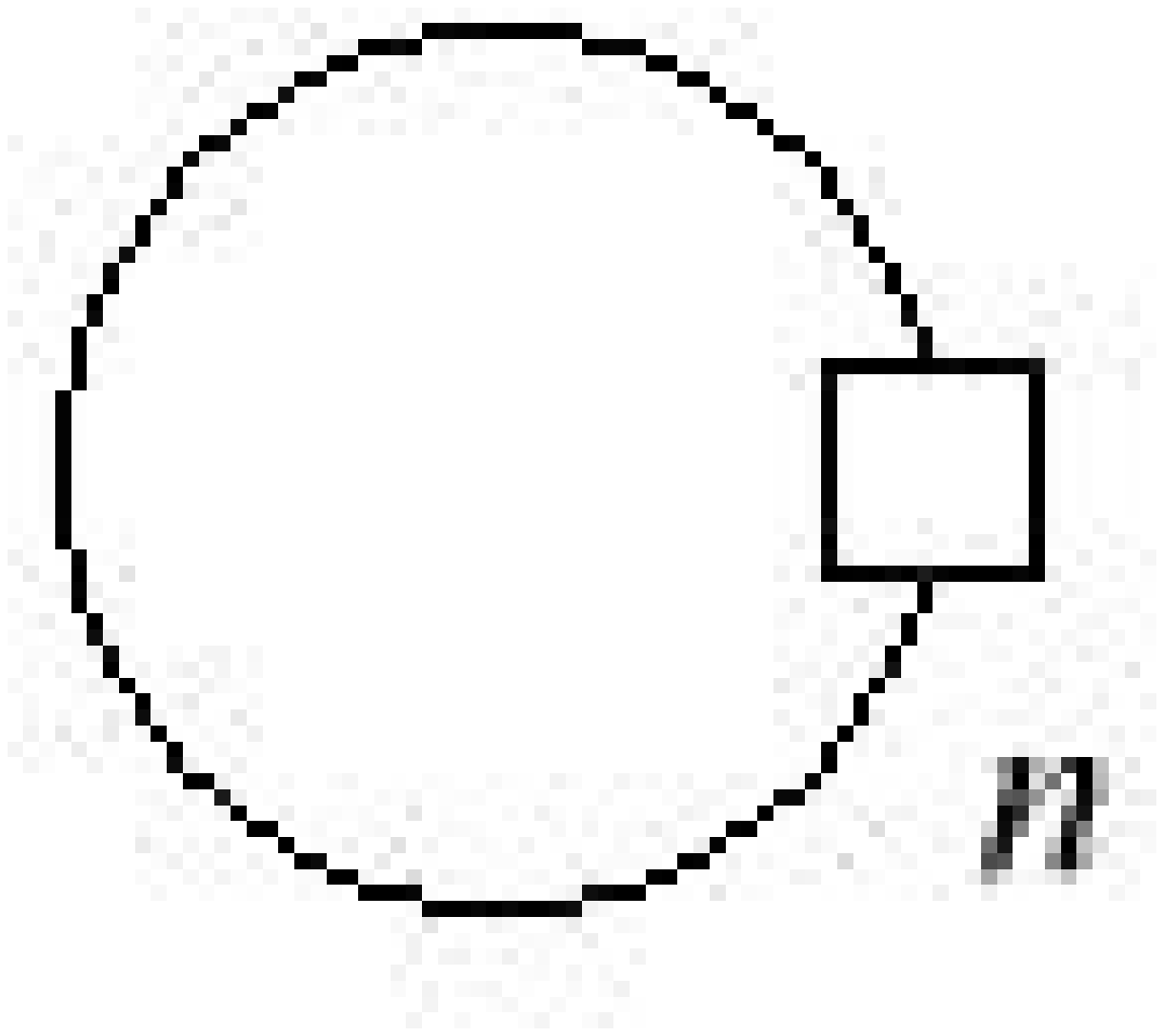}\end{minipage}$$
\caption{Definition of Jones-Wenzl idempotents}
\label{fig:idempotent}
\end{figure}

An arc labelled $n$ represents the $n$th Jones-Wenzl idempotent $f_n$:  $f_0$ is an empty tangle, $f_1$ is a single arc, and $f_n$ is a linear combination of $(n,n)$-tangles defined by the recursive relation in Figure \ref{fig:idempotent}.

Note that, for $f_n$ to be defined, $\Delta_{k-1}$ must be invertible in ${\cal R}$ for all $k \leq n$.  This holds, as $\Delta_{k-1} = (-1)^{k-1} [k]$ where  the quantum integer $[k] = A^{2-2k} + A^{6-2k} + \cdots + A^{2n-k} + A^{2k-2}$.  In Figure \ref{fig:idempotent}, the presence of a idempotent is indicated by a small rectangle.  Hereafter, the rectangles will be dropped.  Moreover, edges representing $f_1$ will be left unlabelled.

\begin{figure}
$$\begin{minipage}{0.8in}\includegraphics[width=0.8in]{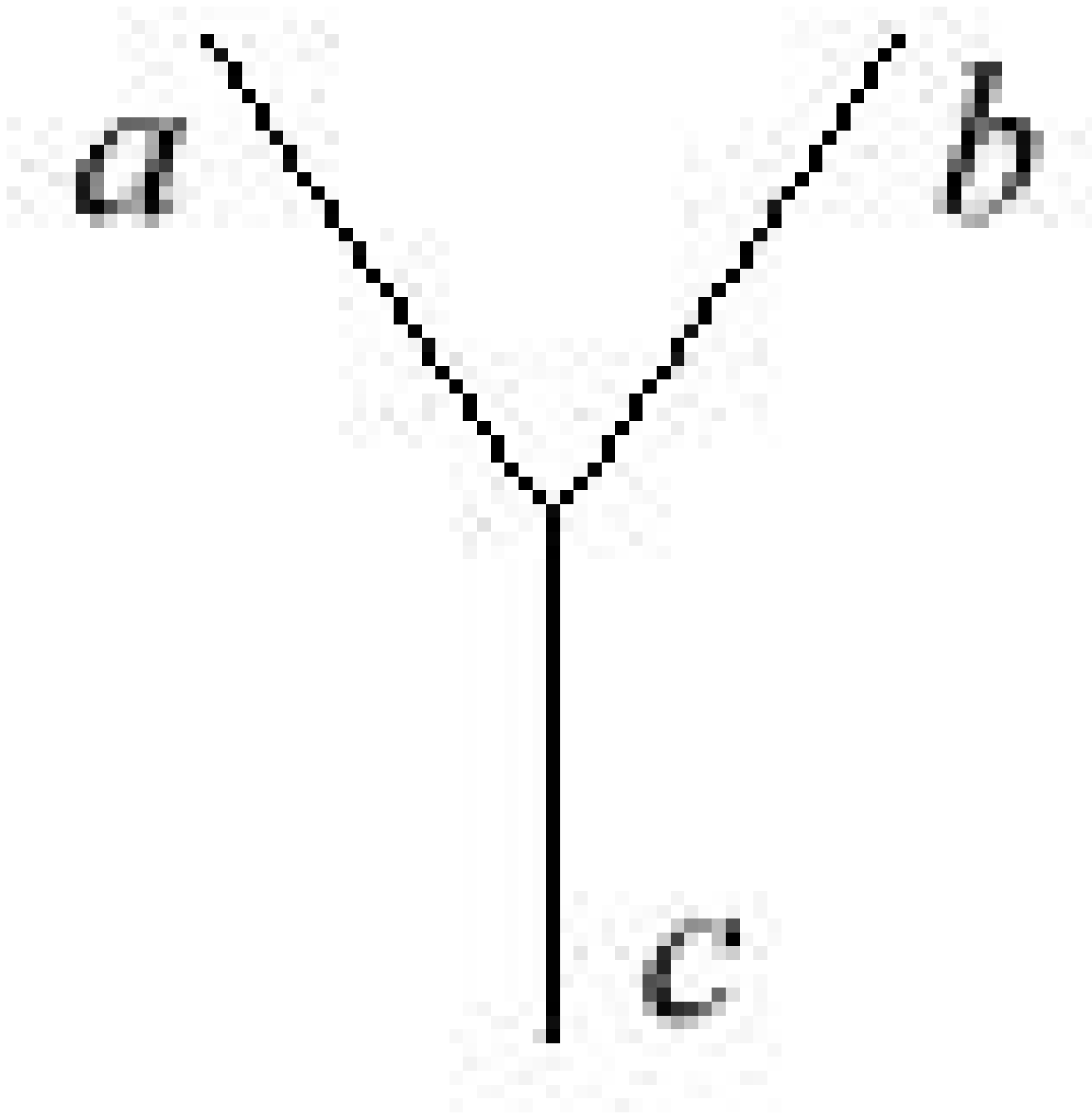}\end{minipage} = \begin{minipage}{0.8in}\includegraphics[width=0.8in]{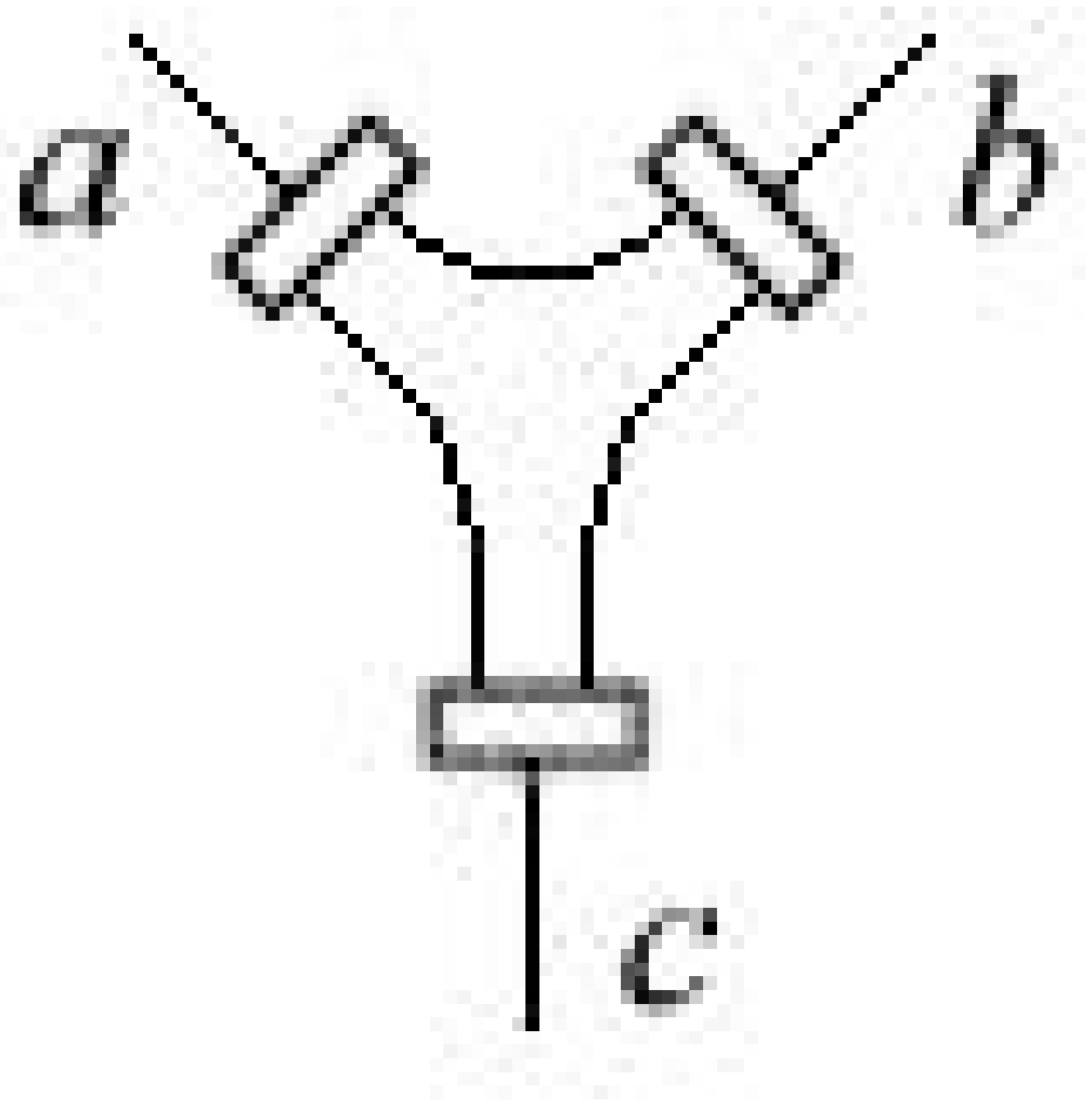}\end{minipage}$$
\caption{Trivalent graphs represent linear combinations of links} 
\label{fig:trivalent}
\end{figure}

It is convenient to extend our view of skein modules to include banded trivalent graphs, whose edges are labelled so that at each vertex, the labels form an admissible triple.  See \cite[4.5]{BHMV2} for the precise definition of a banded trivalent graph:  all trivalent graphs we discuss will be assumed to have such a banding, coming from a regular neighborhood in the surface in which they are drawn.  These graphs represent linear combinations of links, as expressed in Figure \ref{fig:trivalent}.

In \cite{MV94}, Masbaum and Vogel describe an algorithm for reducing every trivalent graph in a 3-ball to a scalar multiple of the empty link, and trivalent graphs are often used to denote this scalar.  In addition to the labelled unknot in Figure \ref{fig:idempotent}, two of these scalars appear frequently enough to merit symbols:

$$\begin{minipage}{0.8in}\includegraphics[width=0.8in]{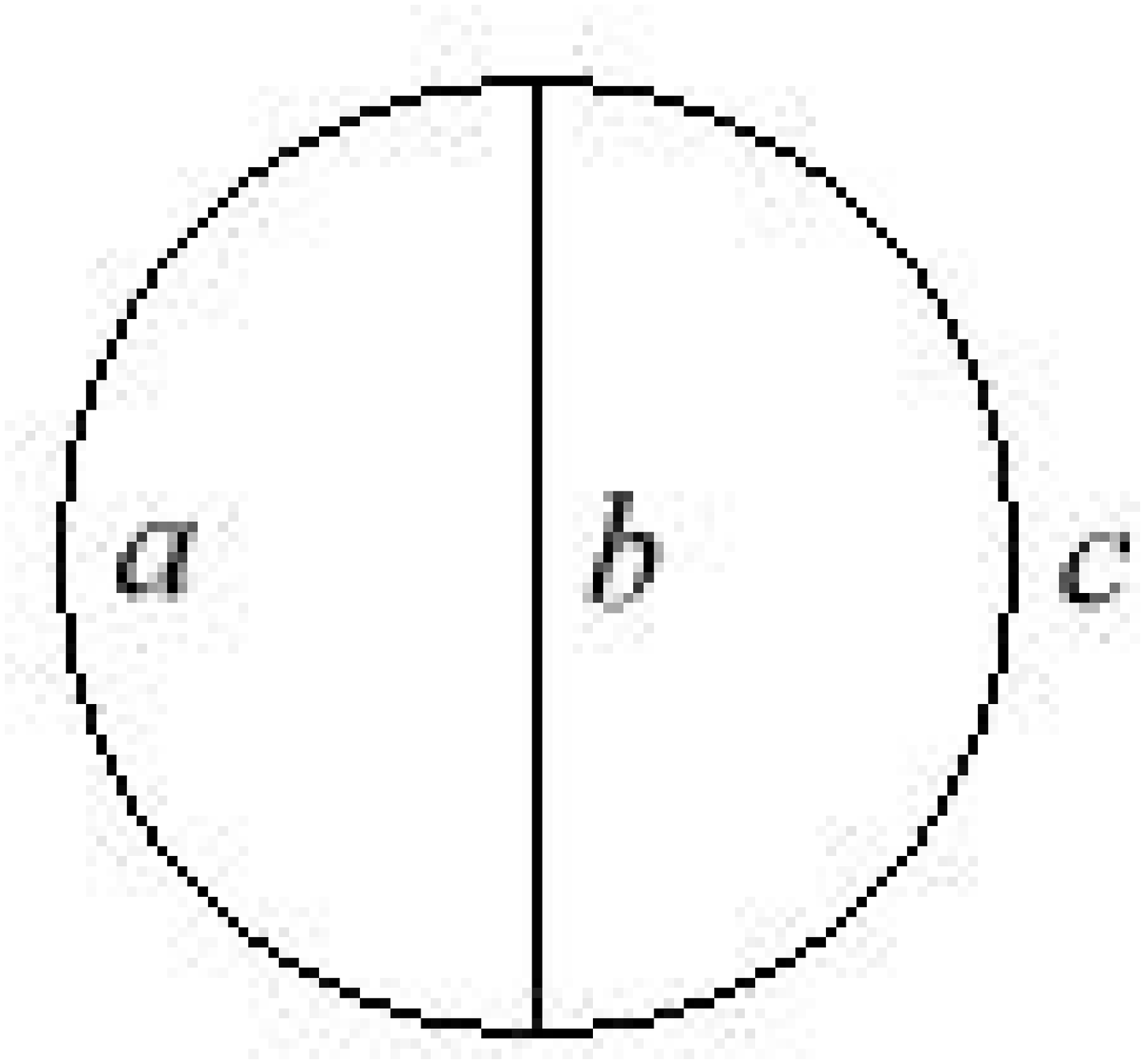}\end{minipage} = \theta(a,b,c) \hbox{ and}$$
$$\begin{minipage}{1.0in}\includegraphics[width=1.0in]{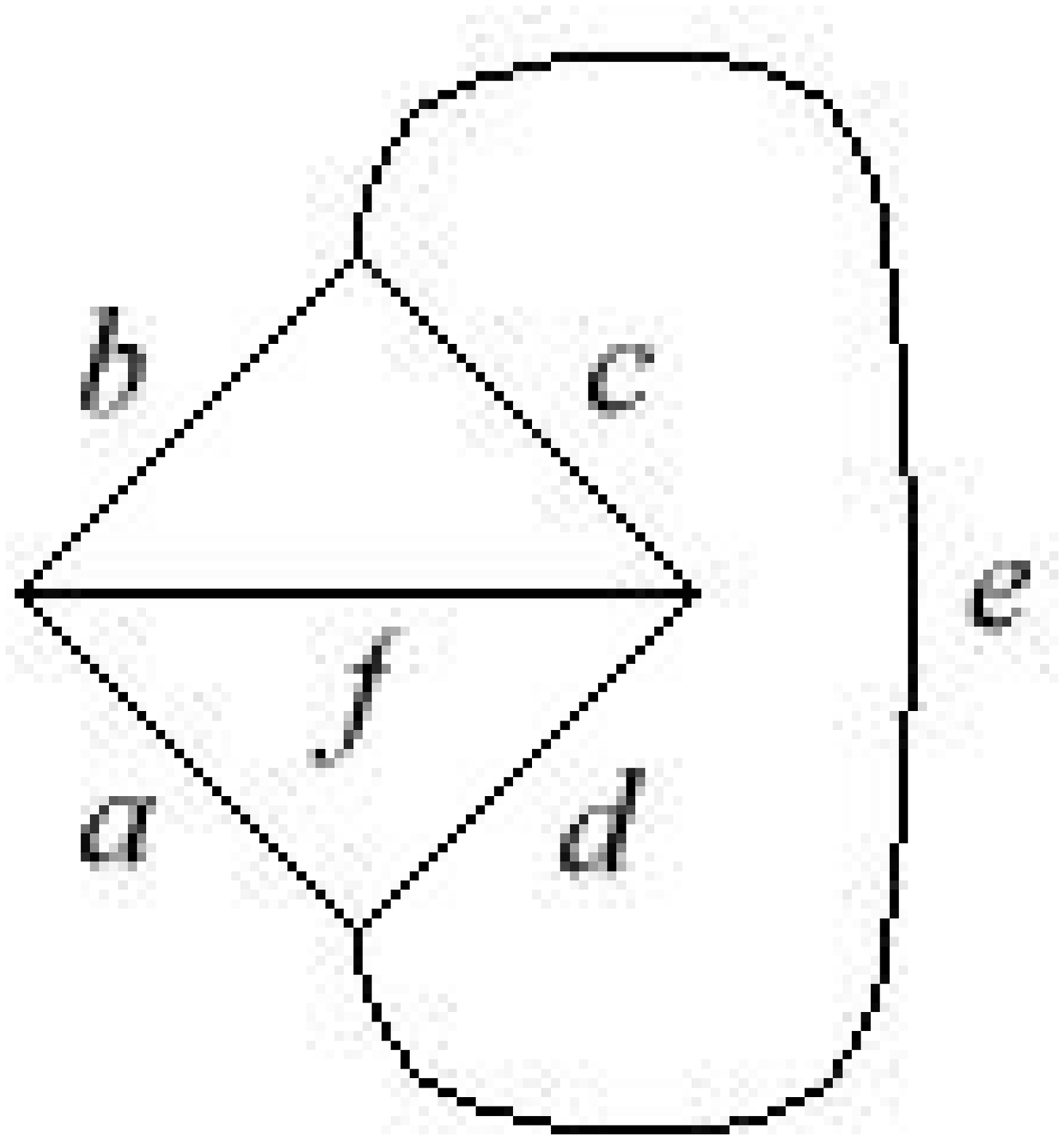}\end{minipage} = \tet a b c d e f .$$
\medskip

Additionally, we have the following useful local moves:
$$\begin{minipage}{0.8in}\includegraphics[width=0.8in]{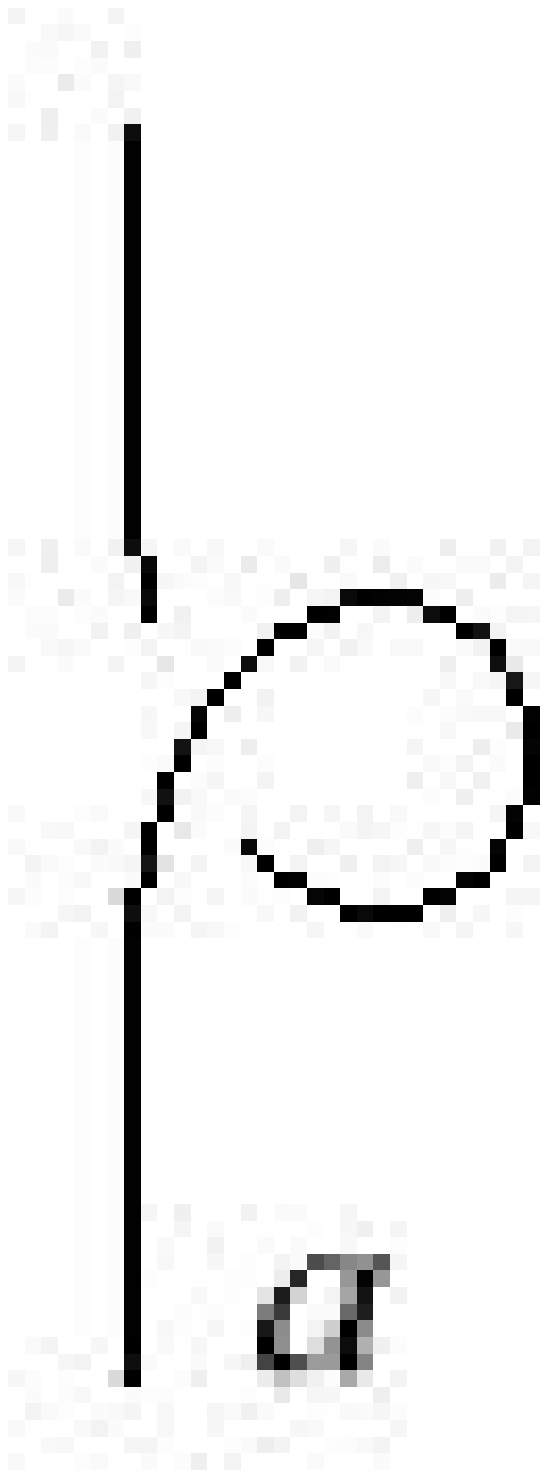}\end{minipage} = (-A)^{a(a+2)}  \begin{minipage}{0.8in}\includegraphics[width=0.8in]{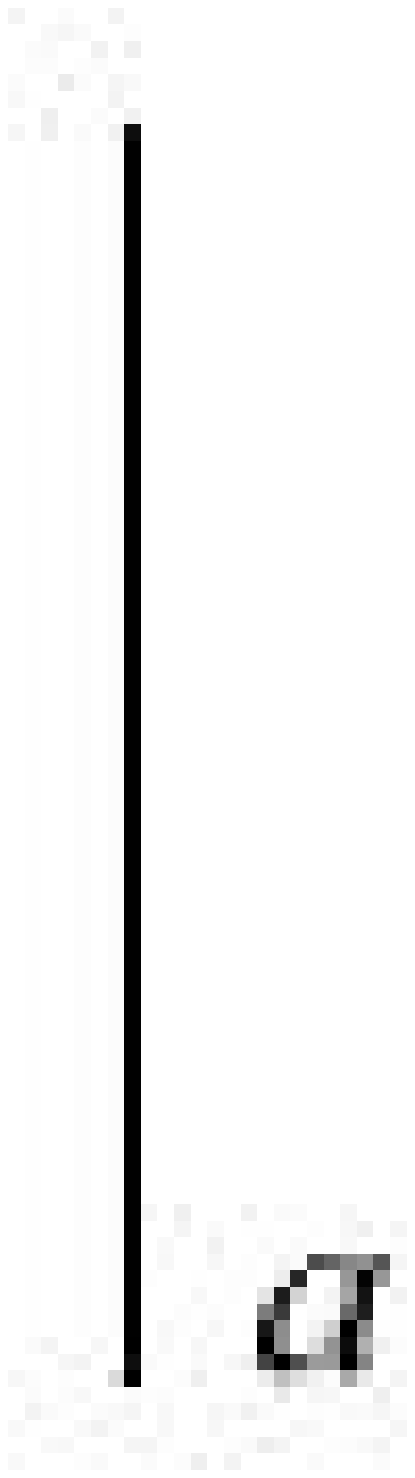}\end{minipage},$$
$$\begin{minipage}{0.8in}\includegraphics[width=0.8in]{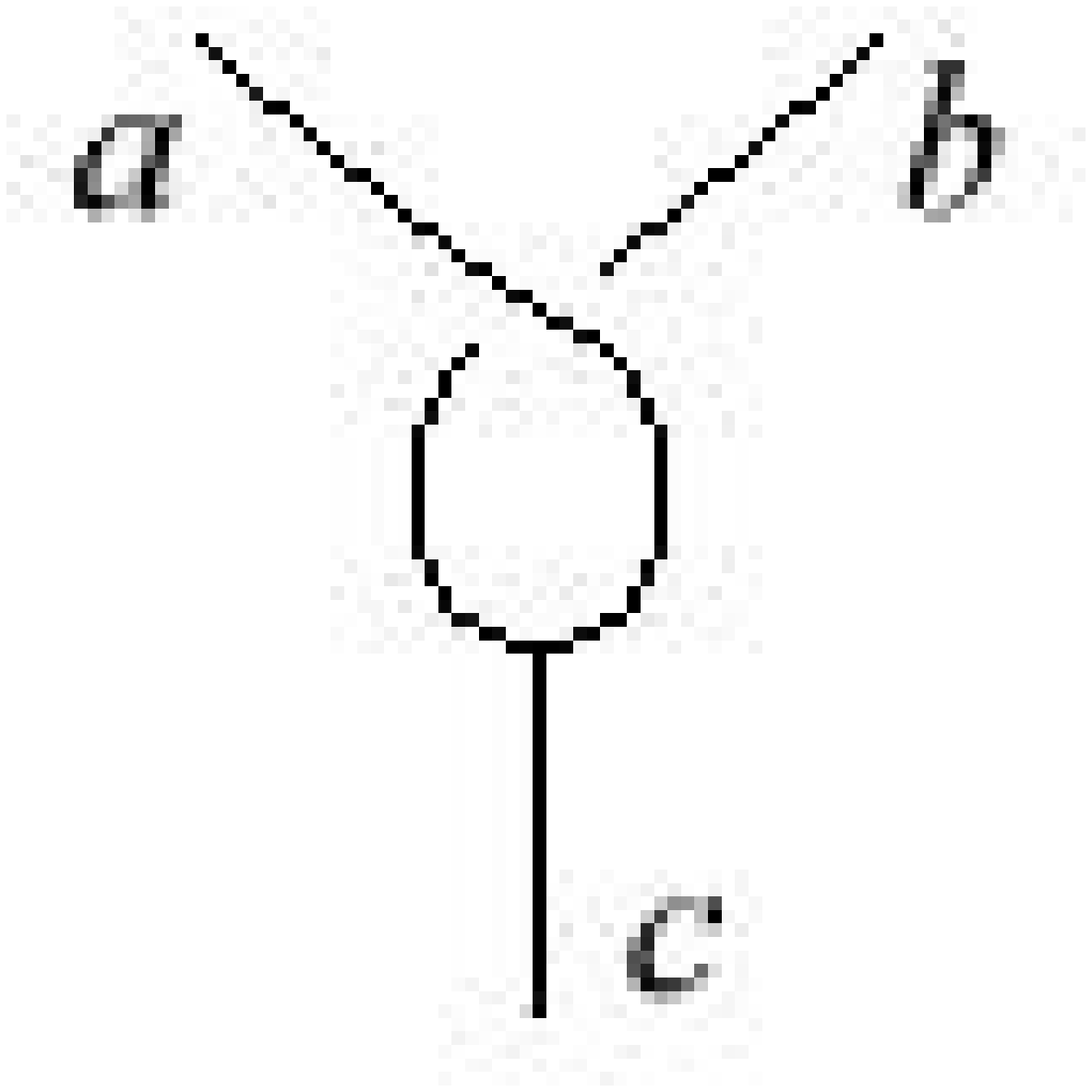}\end{minipage} = \lambda_c^{a \;b} \begin{minipage}{0.8in}\includegraphics[width=0.8in]{graphics/trivalent}\end{minipage},$$
$$\begin{minipage}{0.8in}\includegraphics[width=0.8in]{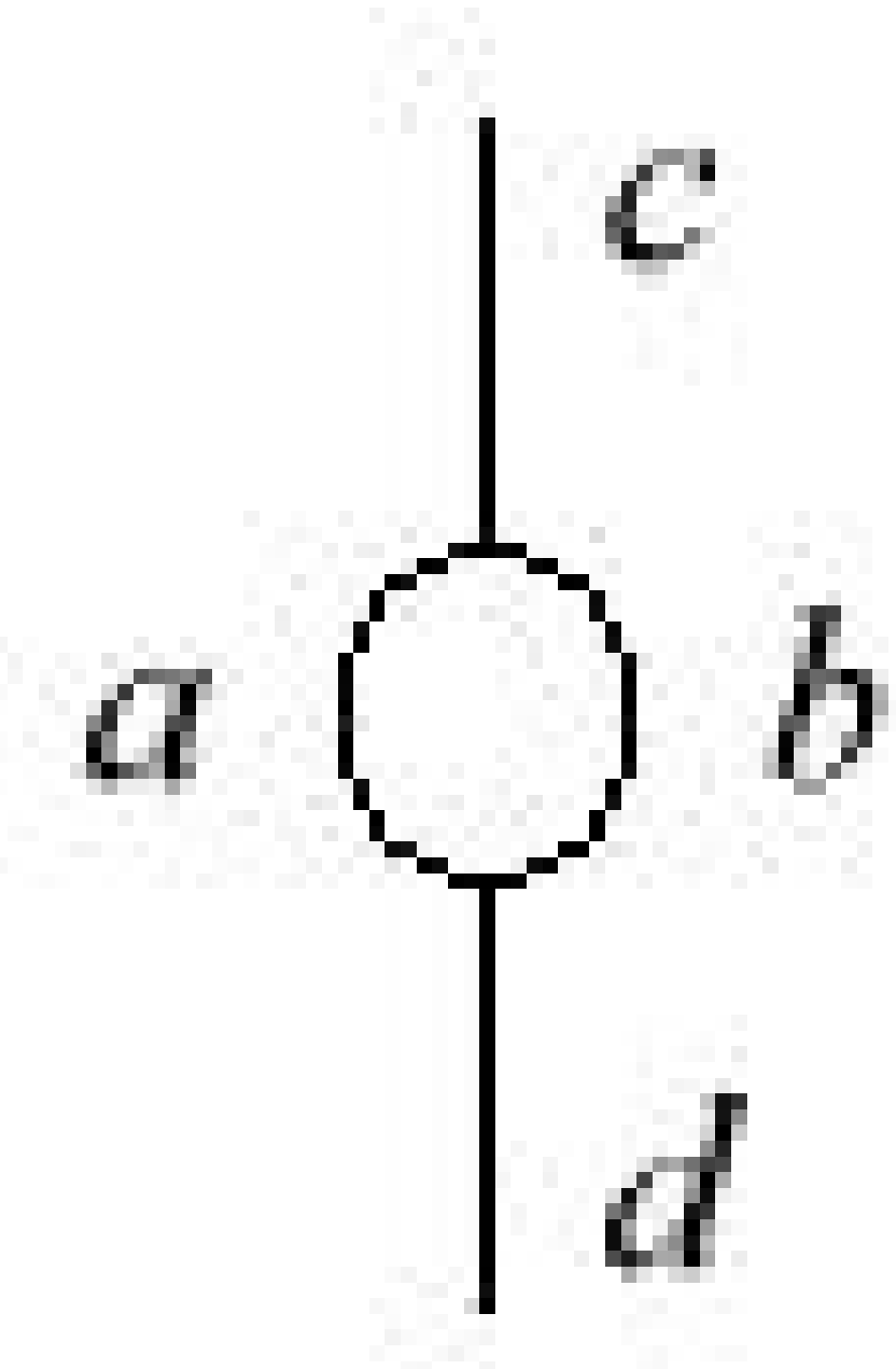}\end{minipage} = \delta_d^c \;\frac{\theta(a,b,c)}{\Delta_c}   \begin{minipage}{0.8in}\includegraphics[width=0.8in]{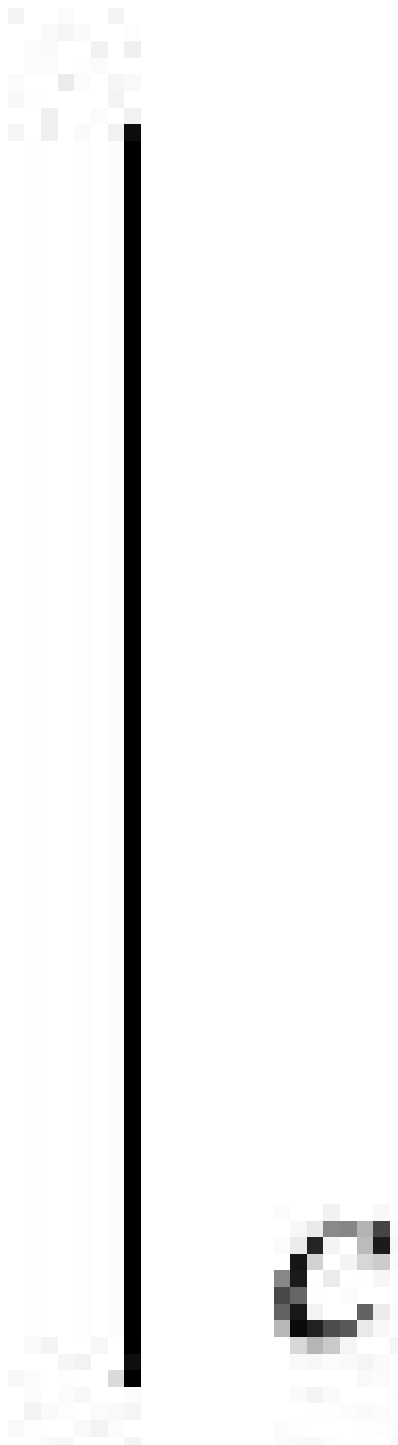}\end{minipage}, \hbox{ and}$$
$$\begin{minipage}{1.0in}\includegraphics[width=1.0in]{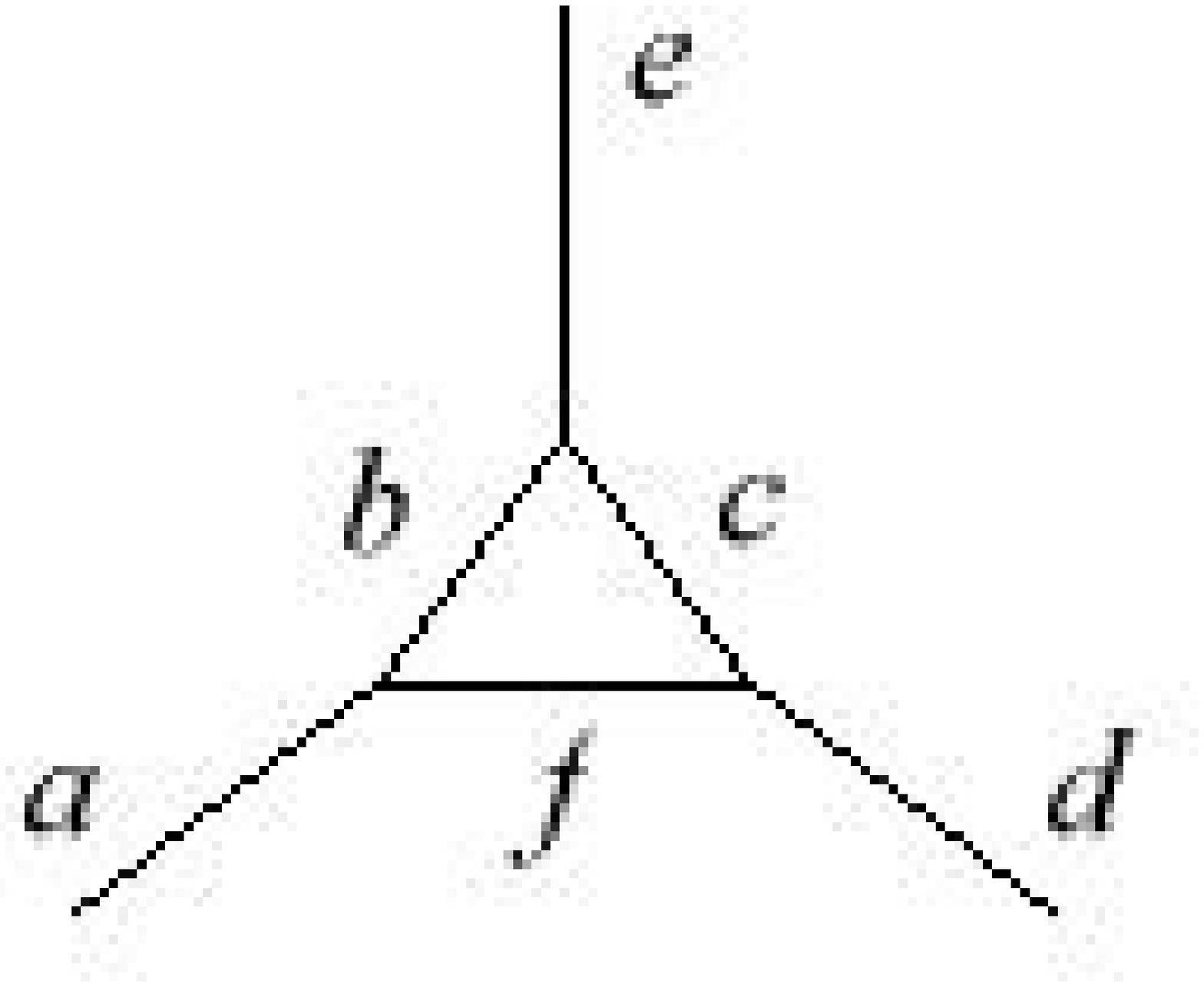}\end{minipage} = \admpwise {\frac{\tet a b c d e f}{\theta(a,d,e)}  \begin{minipage}{0.8in}\includegraphics[width=0.8in]{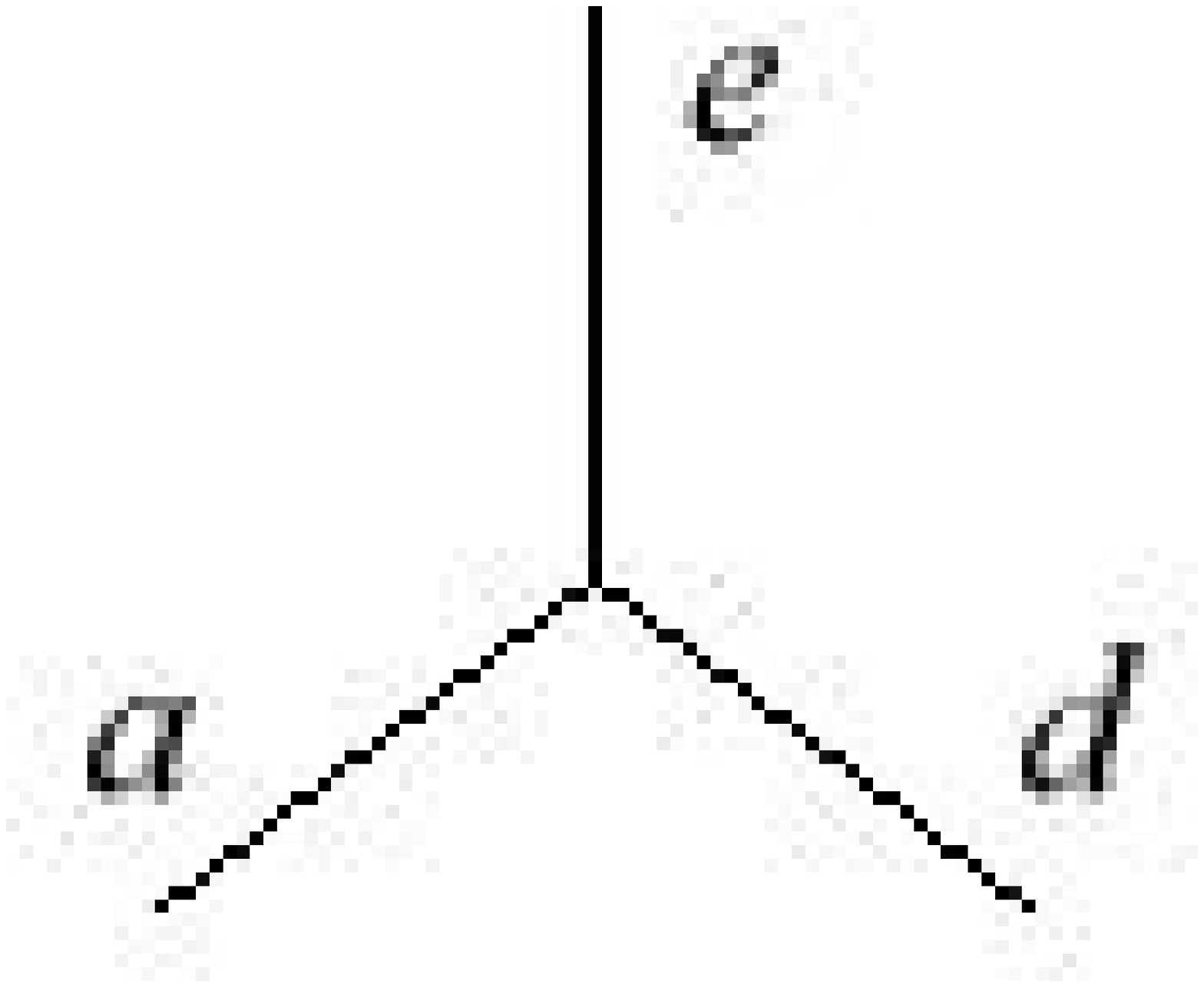}\end{minipage}} {(a,d,e)} .$$

Explicit formulas for these scalars in terms of $A$ are given in \cite{MV94} and \cite{KL94}.  The only skein-theoretic derivation of the formula for the tet is given in \cite{MV94}.  Our notation is the same as in \cite{KL94}.

Working with trivalent graphs is often much easier than working directly with links, thanks mainly to the following identity:

\begin{thm}{(Fusion Formula)}
$$\begin{minipage}{0.8in}\includegraphics[width=0.8in]{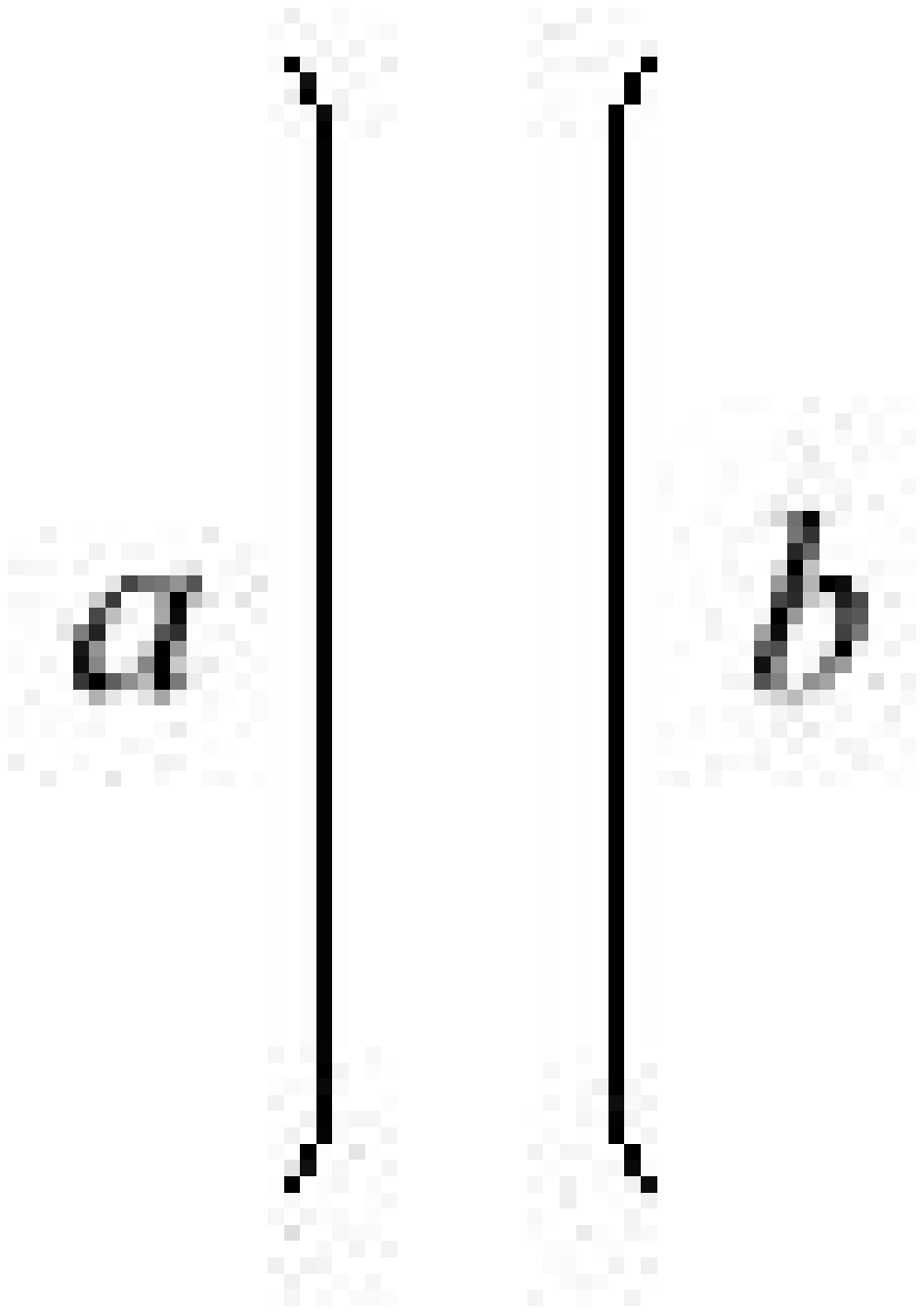}\end{minipage}=\sum_i \frac{\Delta_i}{\theta(a,b,i)} \begin{minipage}{0.8in}\includegraphics[width=0.8in]{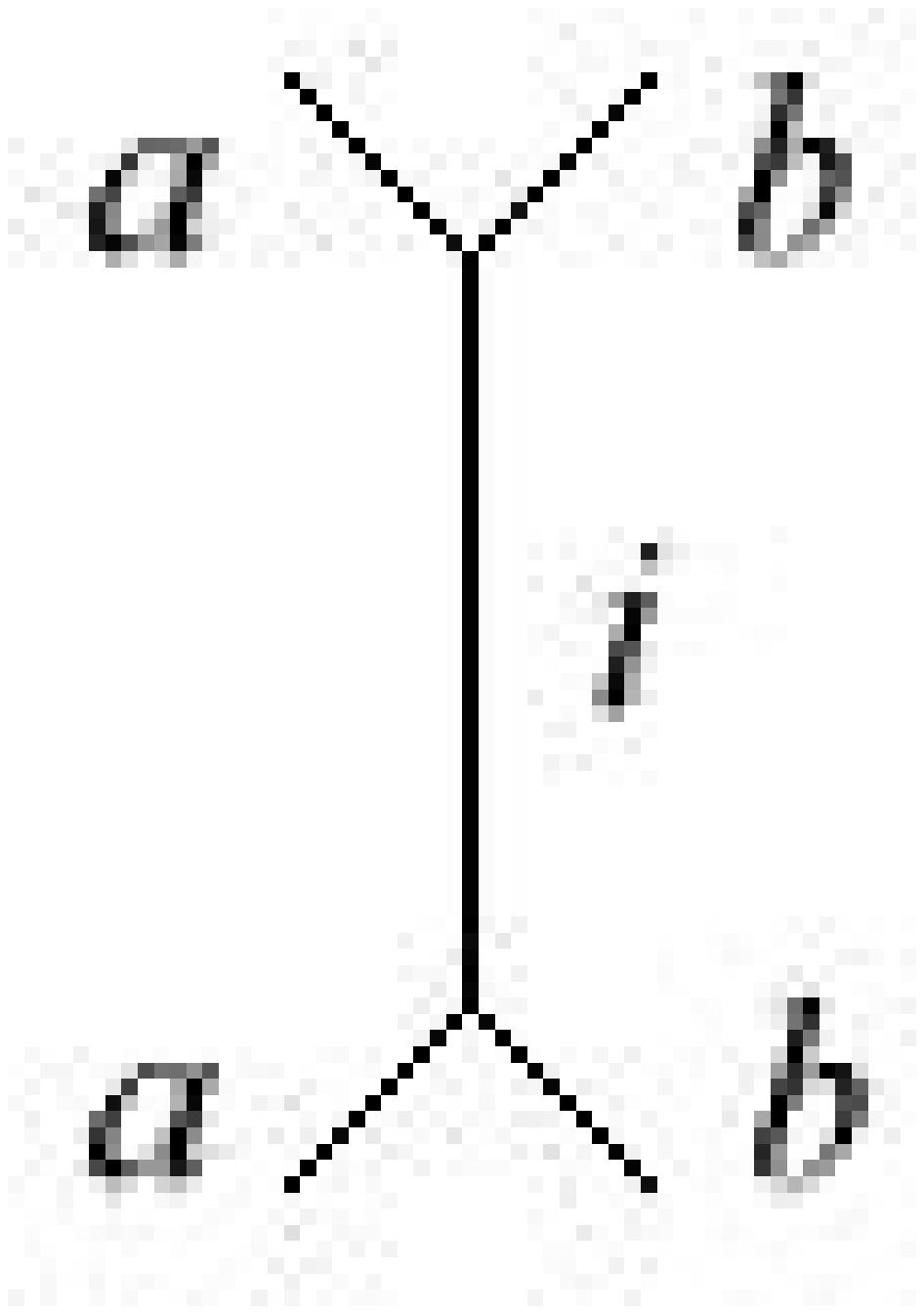}\end{minipage},$$ where the sum is over all admissible labellings.
\end{thm}

We also have the following well-known theorem.  See, for example, \cite{BFK00}. 

\begin{thm}{(Sphere Lemma)}
\label{thm:sphere}

If a sphere intersects a skein element in exactly 1 labelled arc
 , then  
$$\begin{minipage}{0.5in}\includegraphics[width=0.5in]{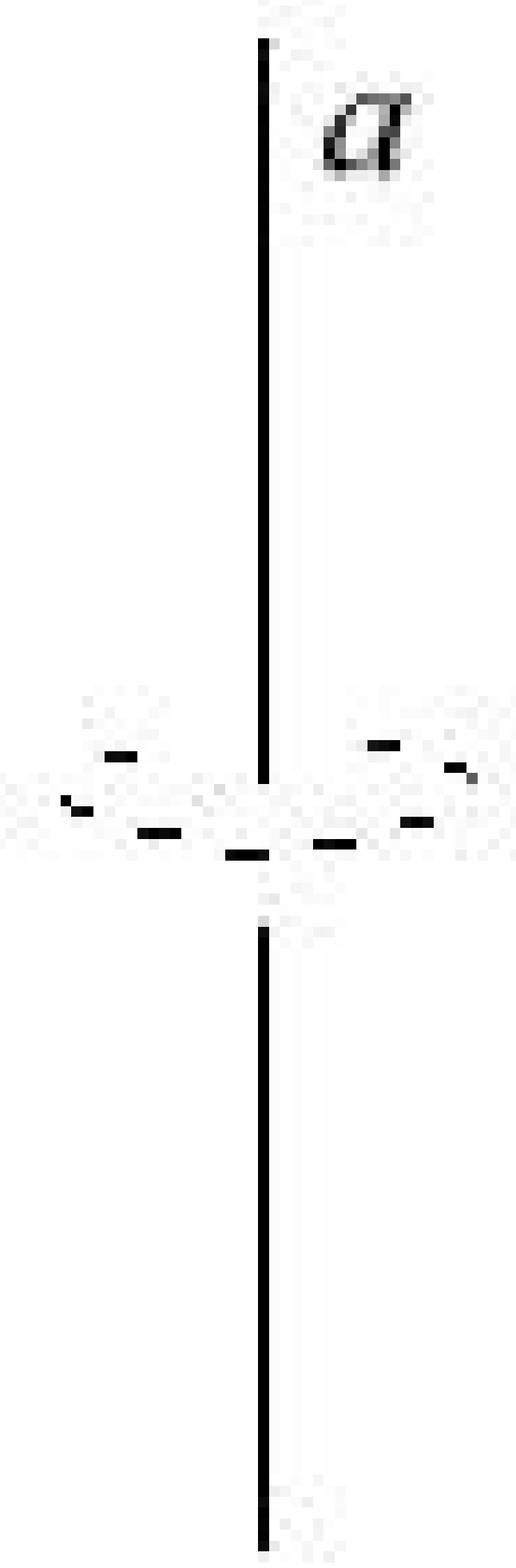}\end{minipage}= 0.$$
\end{thm}

The idea of using the fusion formula in conjunction with the sphere lemma is also well-known.  See, for example, \cite{Ro95}.  

\begin{thm}
If a sphere intersects a skein element in exactly 2 labelled arcs, then 
$$\begin{minipage}{0.75in}\includegraphics[width=0.75in]{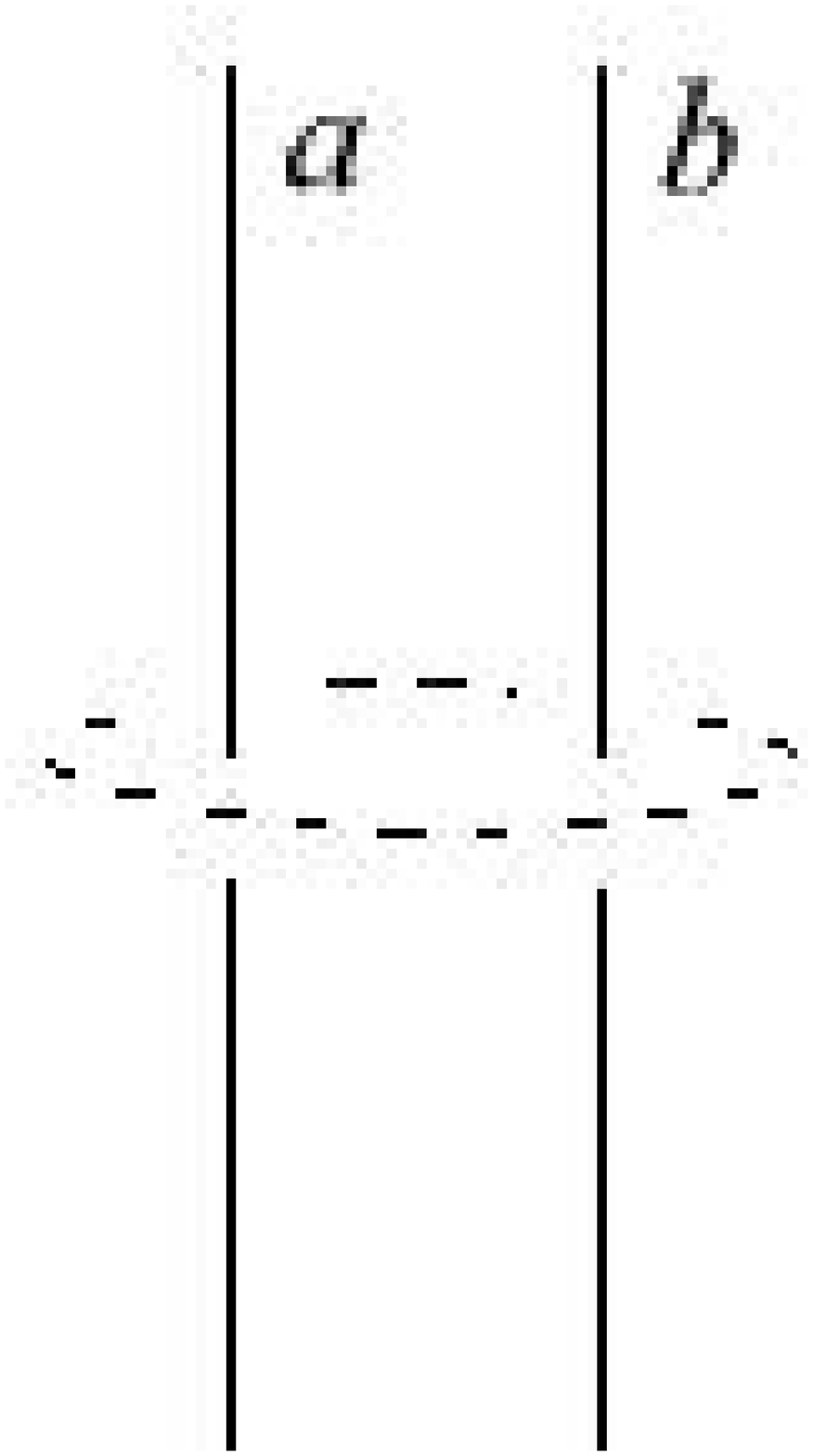}\end{minipage}=\frac{\delta_b^a}{\Delta_a} \begin{minipage}{0.75in}\includegraphics[width=0.75in]{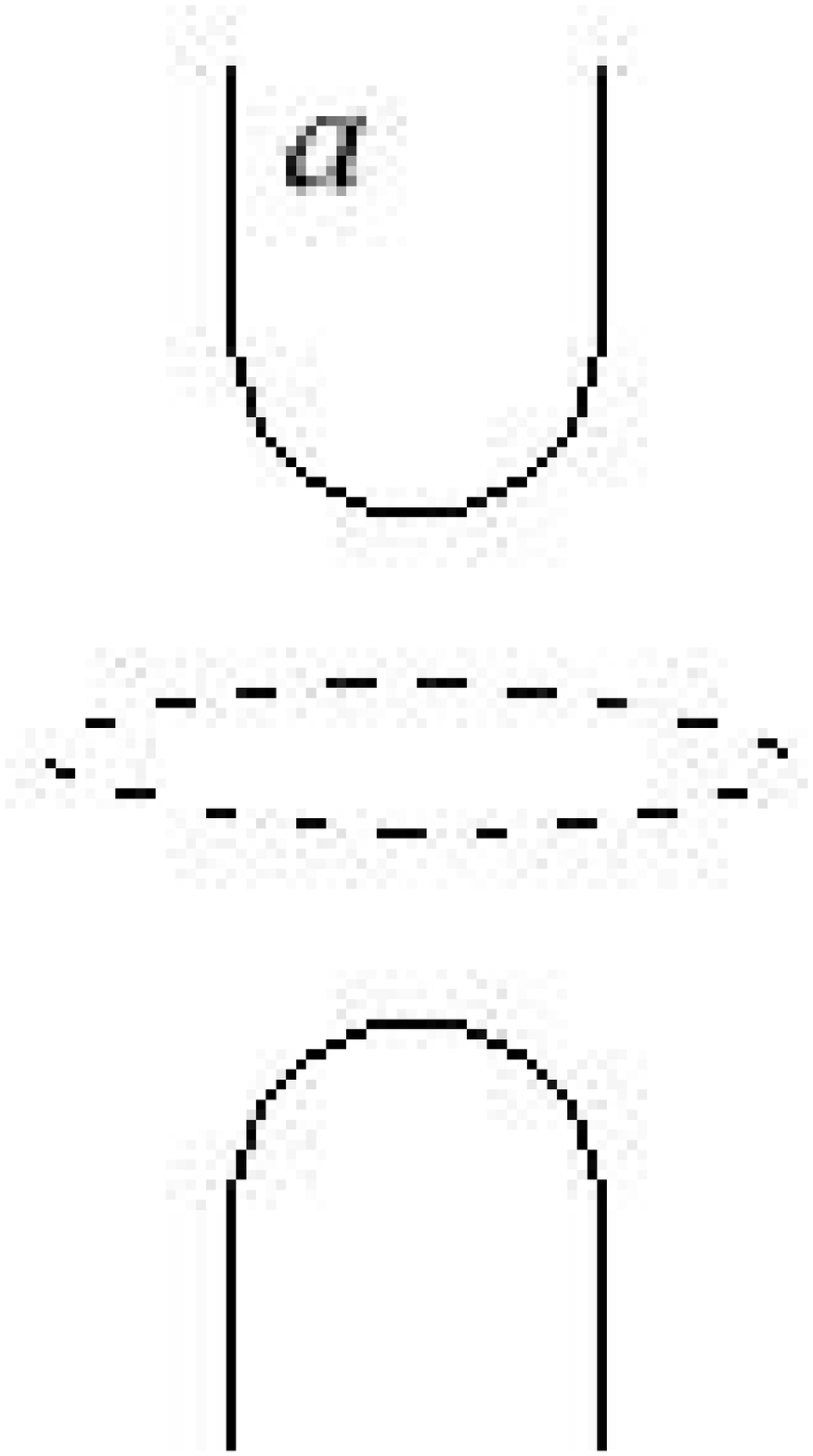}
\end{minipage}.$$

If a sphere intersects a skein element in exactly 3 labelled arcs, then 
$$\begin{minipage}{1.0125in}\includegraphics[width=1.0125in]{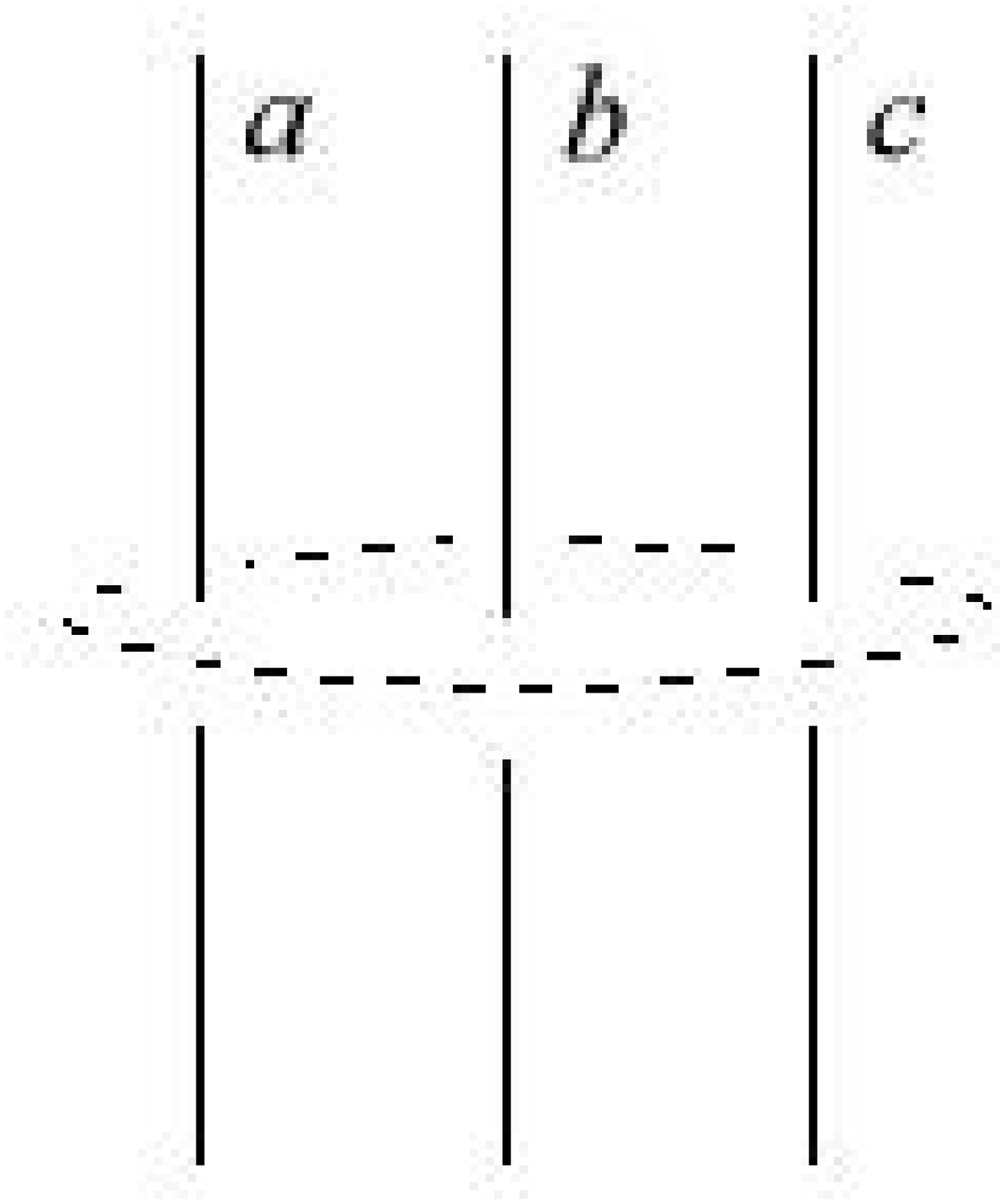}\end{minipage}=\admpwise {\frac{1}{\theta(a,b,c)} \begin{minipage}{1.0125in}\includegraphics[width=1.0125in]{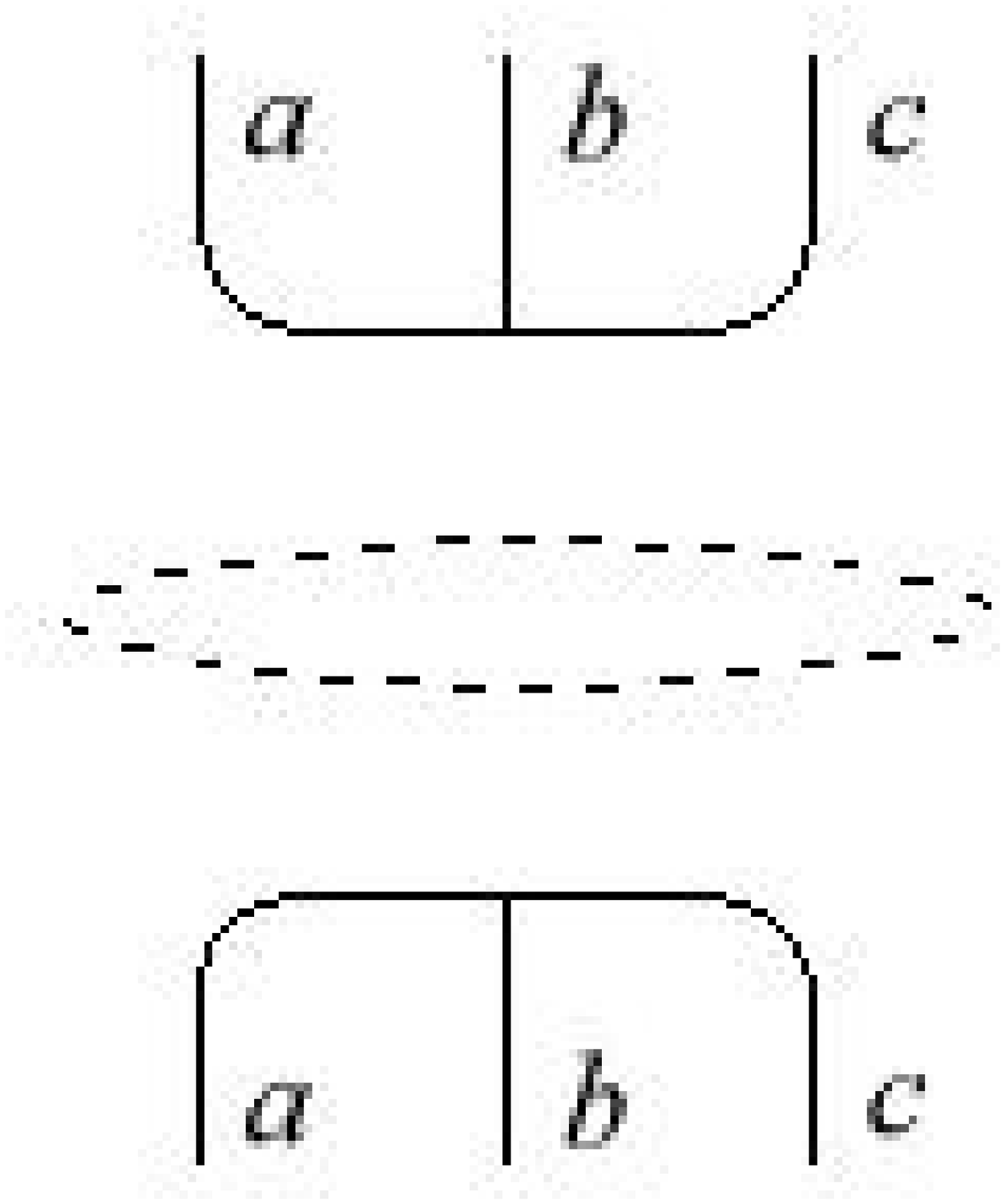}
\end{minipage}} {(a,b,c)} .$$

If a sphere intersects a skein element in exactly $n>3$ arcs, then
$$\begin{minipage}{1.4in}\includegraphics[width=1.4in]{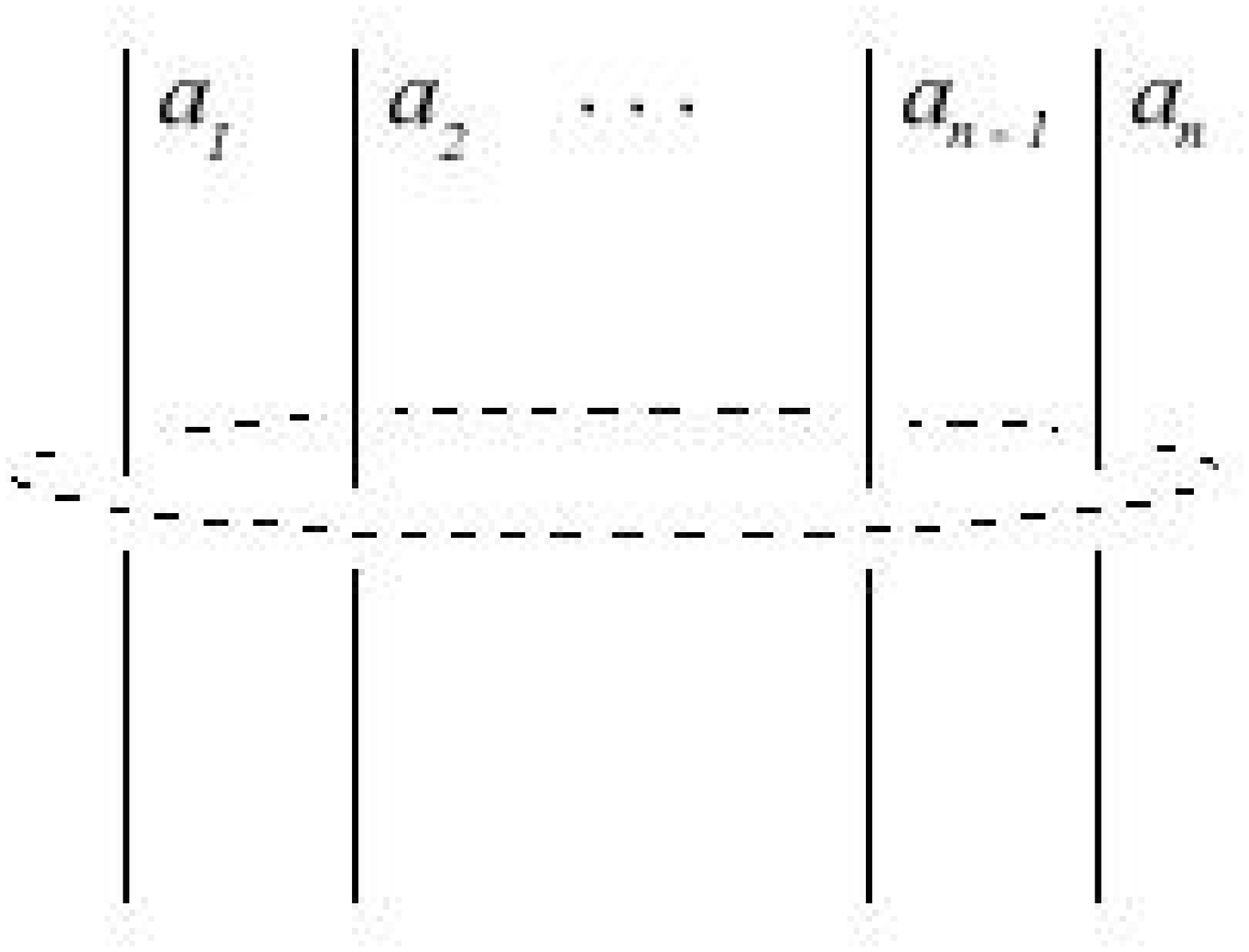}\end{minipage}=\sum \frac{1}{ \begin{minipage}{1.4in}\includegraphics[width=1.4in]{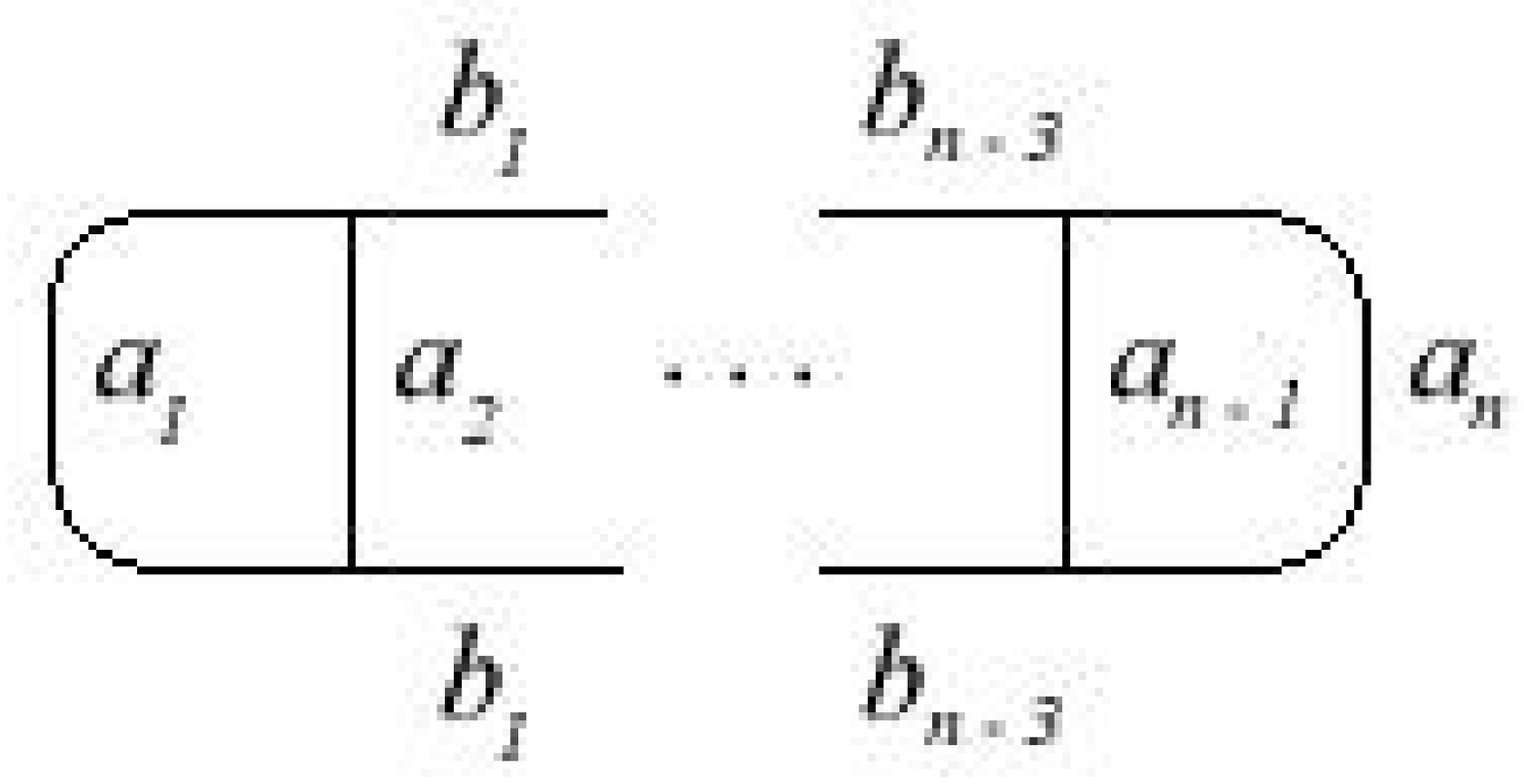}
\end{minipage}} \begin{minipage}{1.6in}\includegraphics[width=1.6in]{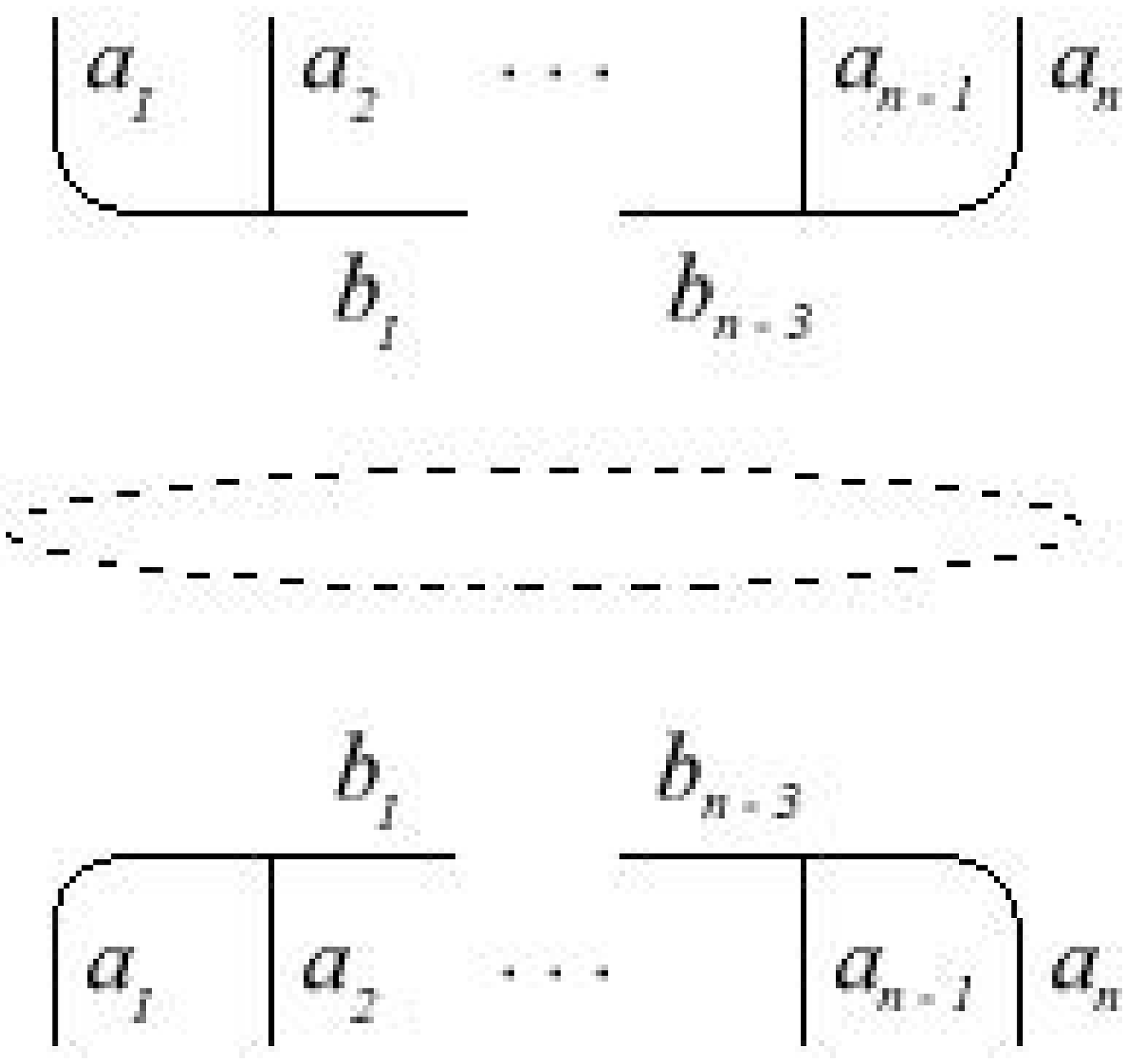}
\end{minipage},$$ where the sum is over all admissible labellings.
\end{thm}

\noindent {\em Proof.} The proof proceeds by induction on $n$.
The case $n=1$ is the sphere lemma.  Using the fusion formula, we can reduce the number of strands by 1 and then apply the inductve hypothesis.  The factor resulting from fusion in each term can be absorbed into its coefficient, and for $n>2$, the coefficient can be rewritten in the desired form with a bigon move.  The result follows.  \qed

Hereafter, we will omit the cumbersome piecewise notation:  when $\theta(a,b,c)$ appears in the denominator of an fraction, the fraction is to be taken as zero if $(a,b,c)$ is not admissible.  We also extend the Kronecker delta notation as follows:  $\delta(a,b,c)$ is 1 when $(a,b,c)$ is admissible and 0 otherwise.

\begin{figure}
$$ \begin{minipage}{2.5in}\includegraphics[width=2.5in]{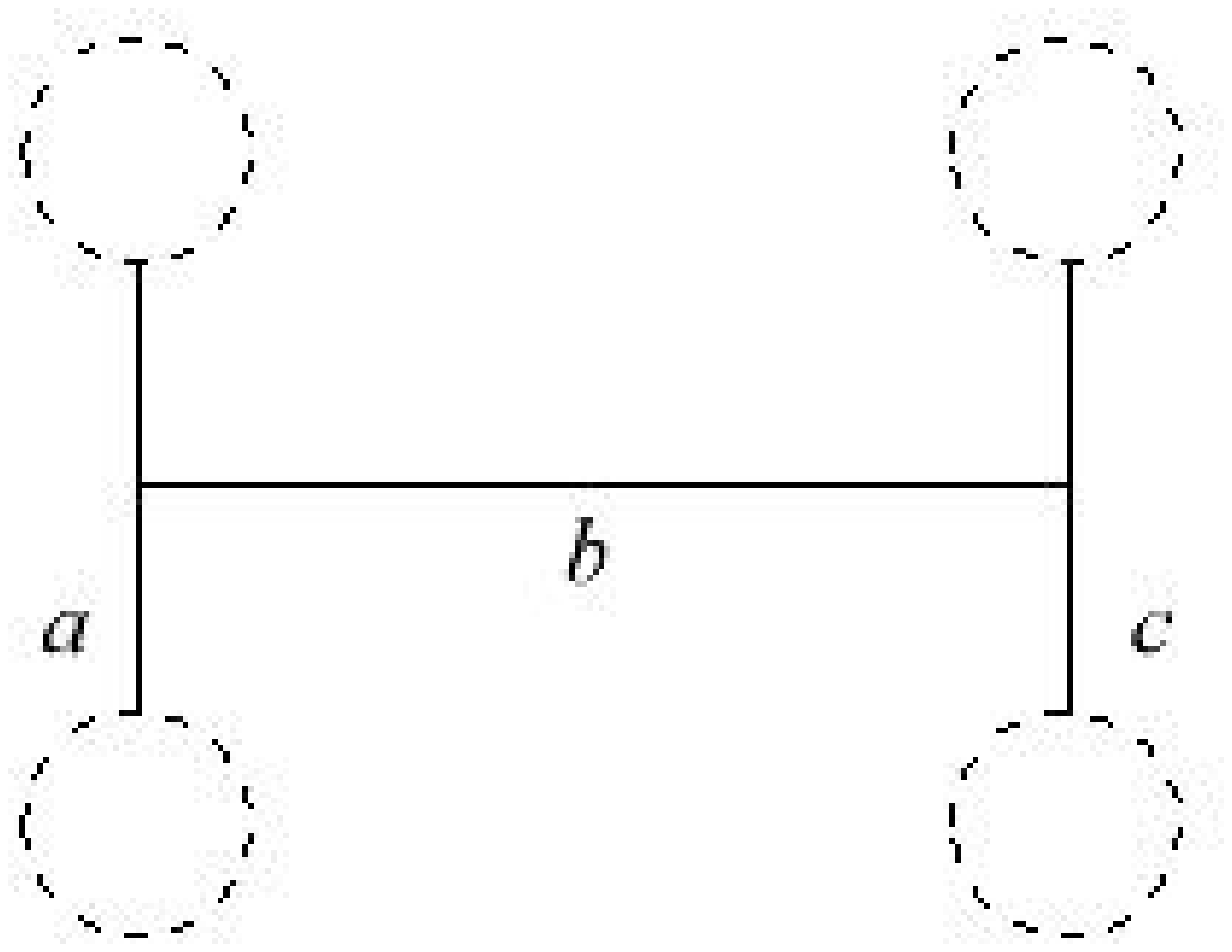}\end{minipage}$$
\caption{The trivalent graph $(a,b,c)$ in $H_2$}
\label{fig:spine}
\end{figure}

For our work, we will use an alternate basis for $S(H_2)$:  the admissible labellings of a trivalent graph which is itself a deformation retract of $H_2$.  We use the notation $(a,b,c)$:  see Figure \ref{fig:spine}.  Bullock, Frohman, and Kania-Bartoszynska describe this basis, for $R = \mathbb{C}$ and $|A| \neq 1$, in \cite{BFK00}.

Note that, for $(a,b,c)$ to be admissibly labelled, $b \leq 2a$, $b \leq 2c$, and $b$ must be even.

Note that $S(S^1 \times S^2) \cong {\cal R}$.  Przytycki has shown in \cite{Pr00} that the skein module of the connected sum of two manifolds is the tensor product of their skein modules.  Hence, $S((S^1 \times S^2) \# 
(S^1 \times S^2)) \cong {\cal R}$.  In \cite{Gi98}, the first author defined an ${\cal R}$-valued bracket evaluation of admissibly labelled trivalent graphs in a connected sum of copies of $S^1 \times S^2$ which essentially coincides with this isomorphism, using fusion as above for the evaluation.  

\begin{figure}
\begin{eqnarray*}
\begin{minipage}{2in}\includegraphics[width=2in]{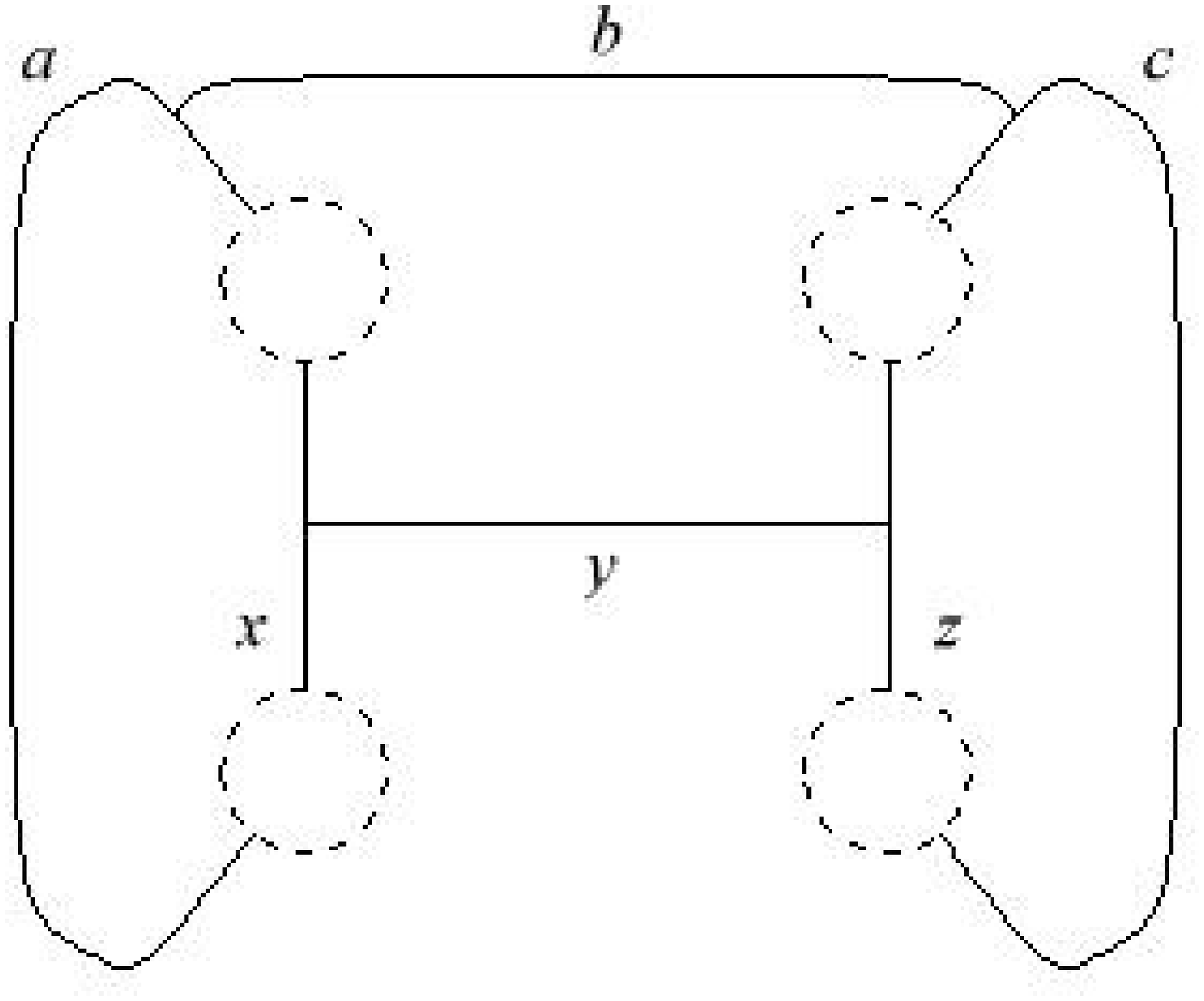}\end{minipage}
 & = & \frac{\delta^a_x \;\delta^b_y  \;\delta^c_z}{\Delta_x \;\Delta_y \;\Delta_z}
\begin{minipage}{1.8in}\includegraphics[width=1.8in]{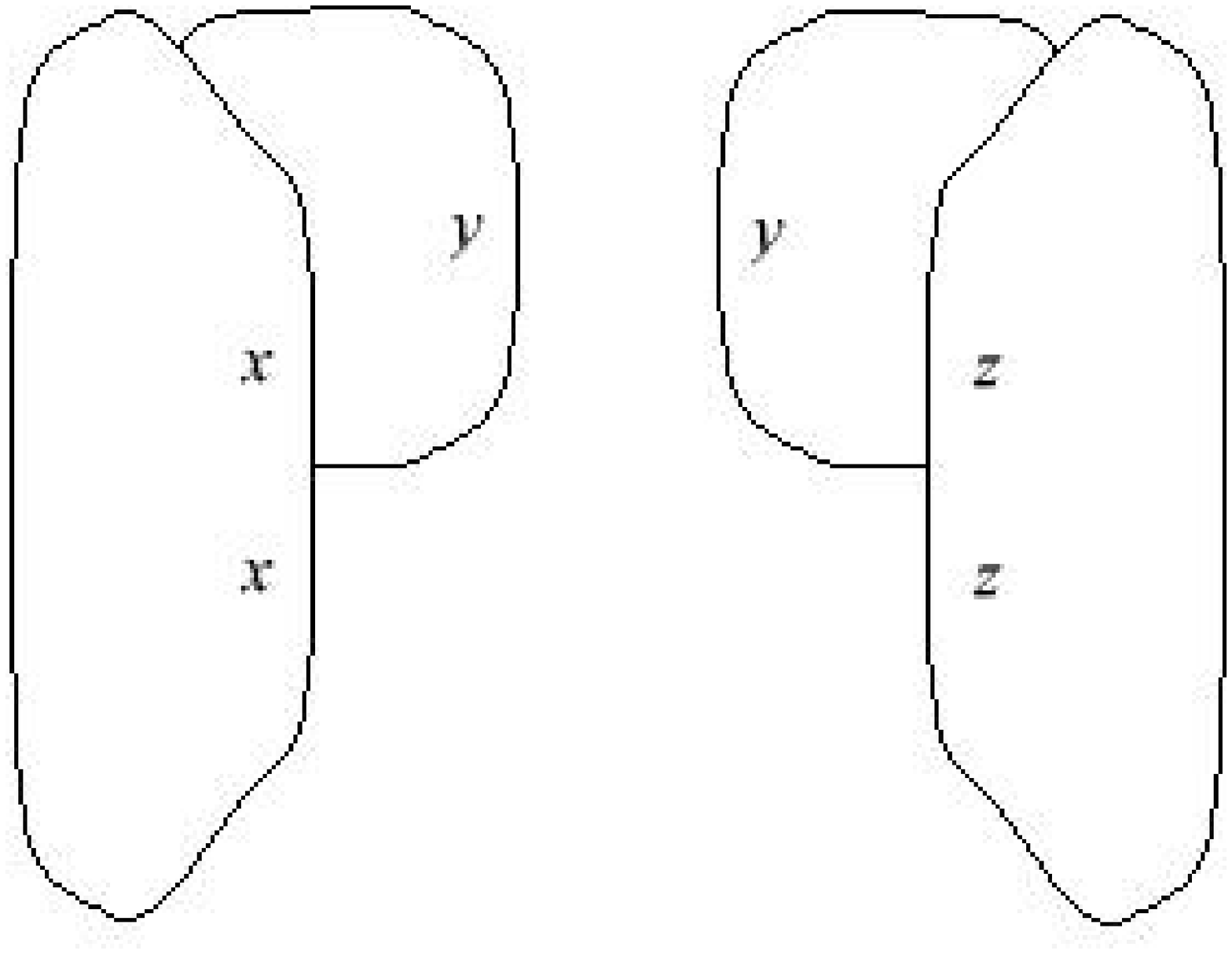}\end{minipage}\\
 & = & \frac{\delta^a_x \;\delta^b_y  \;\delta^c_z \;\theta(x,x,y) \;\theta(y,z,z)}{\Delta_x \;\Delta_y \;\Delta_z}
\end{eqnarray*}
\caption{Computation of $<(x,y,z),(a,b,c)>$}
\label{fig:norm}
\end{figure}

By gluing a copy of $H_2$ to itself with orientation reversed along the identity map on the boundaries, we obtain a Hermitian form $<,>:S(H_2) \times S(H_2)\longrightarrow {\cal R}$, 
and $\{(a,b,c)\}$ is orthogonal with respect to this product.  (Here ${\cal R}$ is equipped with the involution which sends $A$ to $A^{-1}$.)  See Figure \ref{fig:norm}.  In the figure, the graph in the outer, undrawn handlebody has been pushed into the inner handlebody.  This form is closely related to the Yang-Mills measure on the skein algebra of a surface discussed by Bullock, Frohman, and Kania-Bartoszynska in \cite{BFK00}. The first author discussed the form $<,>$ and the orthogonality of our basis of trivalent graphs in a course on quantum topology in the fall of 2001.

We now prove that the alternate basis described above is indeed a basis, as there does not seem to be a proof in the literature.  In the proof, we use the principle of well-founded induction:

A well-founded order on a set $X$ is a partial order such that every nonempty subset of $X$ has a minimal element.  The principle of well-founded induction states that given a property $p$ defined on a well-founded ordered set $X$, if $p$ holds for every minimal element of $X$, and if, for every $y \in X$, $p$ holds for $y$ if $p$ holds for every $x<y$ in $X$, then $p$ holds for every element of $X$.  See, for example, \cite{Mo89}.

\begin{propo}
$\{(a,b,c)\}$ is a basis of $S(H_2)$.
\end{propo}
 
\noindent {\em Proof.} Let $\leq_B$ denote the partial ordering on the standard basis $B$ defined by 
$$x^{i_1} y^{j_1} z^{k_1} \leq_B x^{i_2} y^{j_2} z^{k_2} \Leftrightarrow 
i_1 \leq i_2, j_1 \leq j_2, \hbox{ and } k_1 \leq k_2.$$

Since $\leq_B$ is well-founded, we can show that $\{(a,b,c)\}$ spans $S(H_2;{\cal R},A)$ by well-founded induction:  

The empty link belongs to $\{(a,b,c)\}$, and so it can certainly be written as a linear combination of elements of $\{(a,b,c)\}$.

Suppose that every $x^{i_1} y^{j_1} z^{k_1} <_B x^i y^j z^k$ can be written  
as a linear combination of elements of $\{(a,b,c)\}$.  $(i + j, 2 j, k + j) = x^i y^j z^k + $ a linear combination of lesser terms, so 
$x^i y^j z^k$ can be expressed as a linear combination of elements in $\{(a,b,c)\}$.  Hence, $\{(a,b,c)\}$ spans $S(H_2)$.

Now suppose that a linear combination $\sum C_{x,y,z} (a,b,c) = 0$.    
Then, for each $(x,y,z)$, $$C_{x,y,z} <(x,y,z),(x,y,z)> = <\sum C_{x,y,z} (a,b,c), (x,y,z)> = 0.$$
Hence, $C_{x,y,z} = 0$ for each $(x,y,z)$, and so, $\{(a,b,c)\}$ is linearly independent. \qed

By viewing a closed 3-manifold in terms of its handle decomposition, we can hope to analyze the structure of its skein module.  A 
closed, oriented, genus-two 3-manifold is built by adding 2 solid cylinders, or 1-handles, to a ball, or 0-handle, and then by attaching 2 thickened disks, or 2-handles, and closed up by adding a 3-handle.  Before the 2-handles are attached, the manifold is a solid handlebody.  As each 2-handle is added, a set of relations is introduced among the skein elements, namely those obtained from sliding arcs over the newly attached thickened disk.  The final 3-handle has no effect on the skein module.

In this way, one obtains a presentation of the skein module.  The 
generators are the basis elements of the solid handlebody, and the 
relations are given by the ways in which arcs may slide across the 
2-handles.  When applied to links in the manifolds, this is the most 
common method for generating a presentation of the module.  See, for 
example, Hoste and Przytycki (\cite{HP93},\cite{HP95}), and Bullock 
(\cite{Bu95}, \cite{Bu97}.

 As Masbaum has pointed out in \cite{Ma96}, a good way to insure that one has found a complete set of relations is to use generators of the relative skein modules of the handlebody relative two points on the boundary along each of the attaching curves to index the relations. However, we do not take this route.  Instead, we pick relatively simple relations which are not guaranteed to be complete.  Then, we use a separate argument to show that the resulting generators are linearly independent.  

We take this alternate route because the relations that are obtained following the first route are too difficult to analyze completely.  However, the second author \cite{Ha03} did study these relations, and using Mathematica \cite{Wo}, found convincing evidence that the skein module of $M$ over the field of rational functions in $A$ should be five-dimensional.

\begin{figure}
$$
\begin{minipage}{2.5in}\includegraphics[width=2.5in]{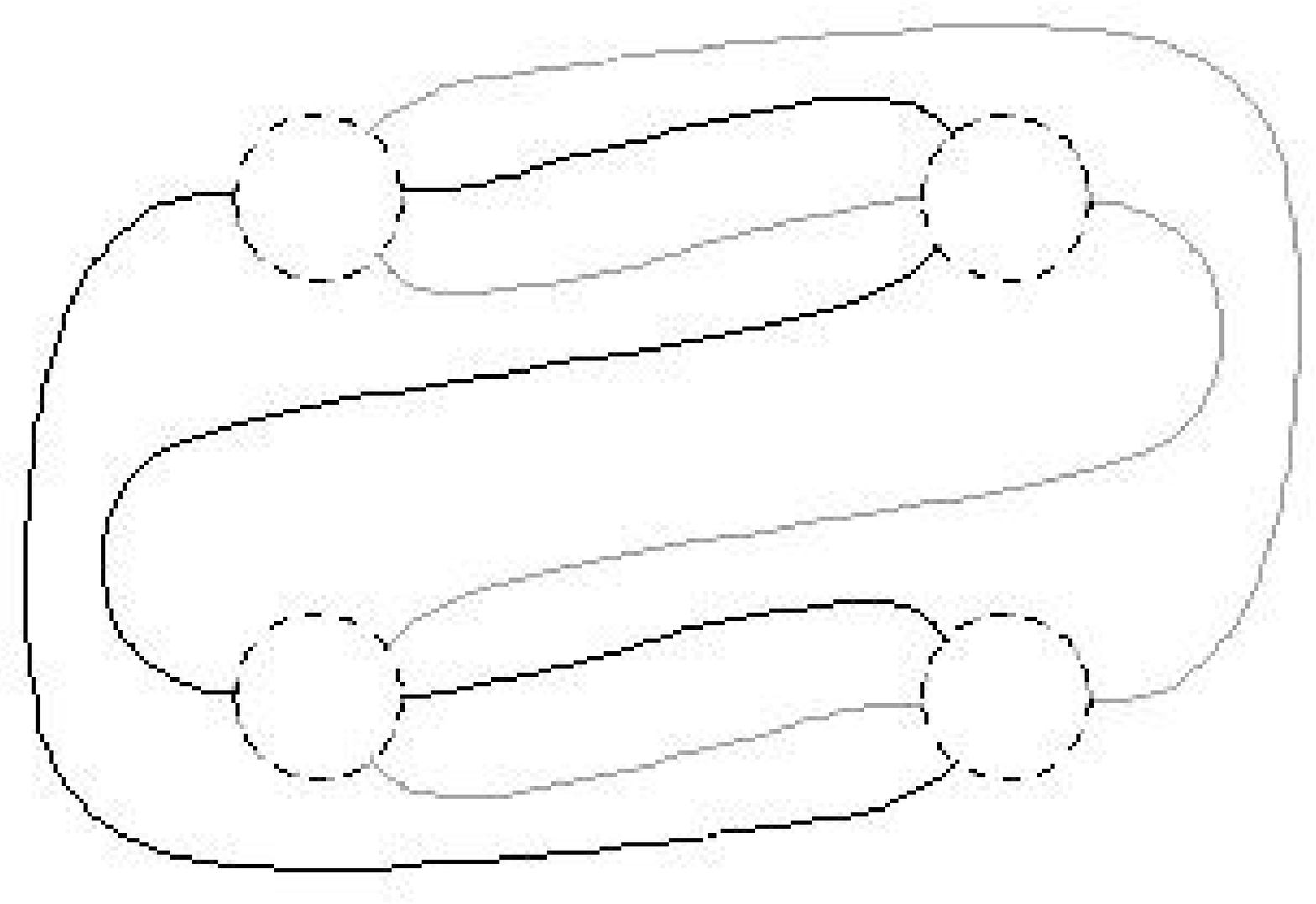}\end{minipage}
$$
\caption{A Heegaard splitting of the quaternionic manifold}
\label{fig:quatsplit}
\end{figure}

Following the method of Rolfsen applied to the Poincare homology sphere in \cite{Ro76}, we can construct a Heegaard splitting of the quaternionic manifold from its surgery description.  See Figure \ref{fig:quatsplit}.  The curves mark the boundaries of the attached 2-handles.

\section{Spanning}
\label{sec:Spanning}

\begin{figure}
\begin{eqnarray*}
\begin{minipage}{1.9in}\includegraphics[width=1.9in]{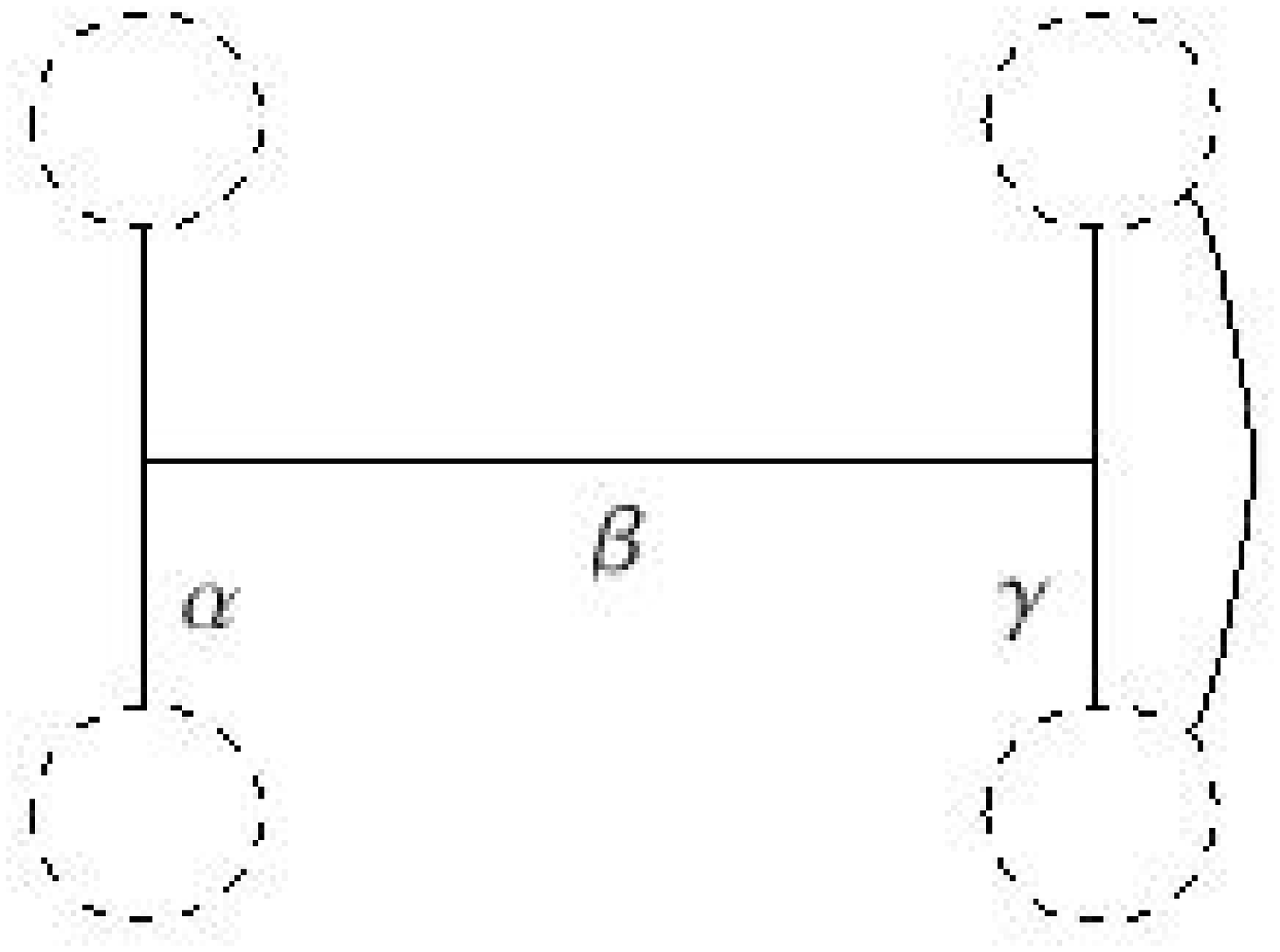}\end{minipage}
 & =  \begin{minipage}{1.9in}\includegraphics[width= 1.9in]{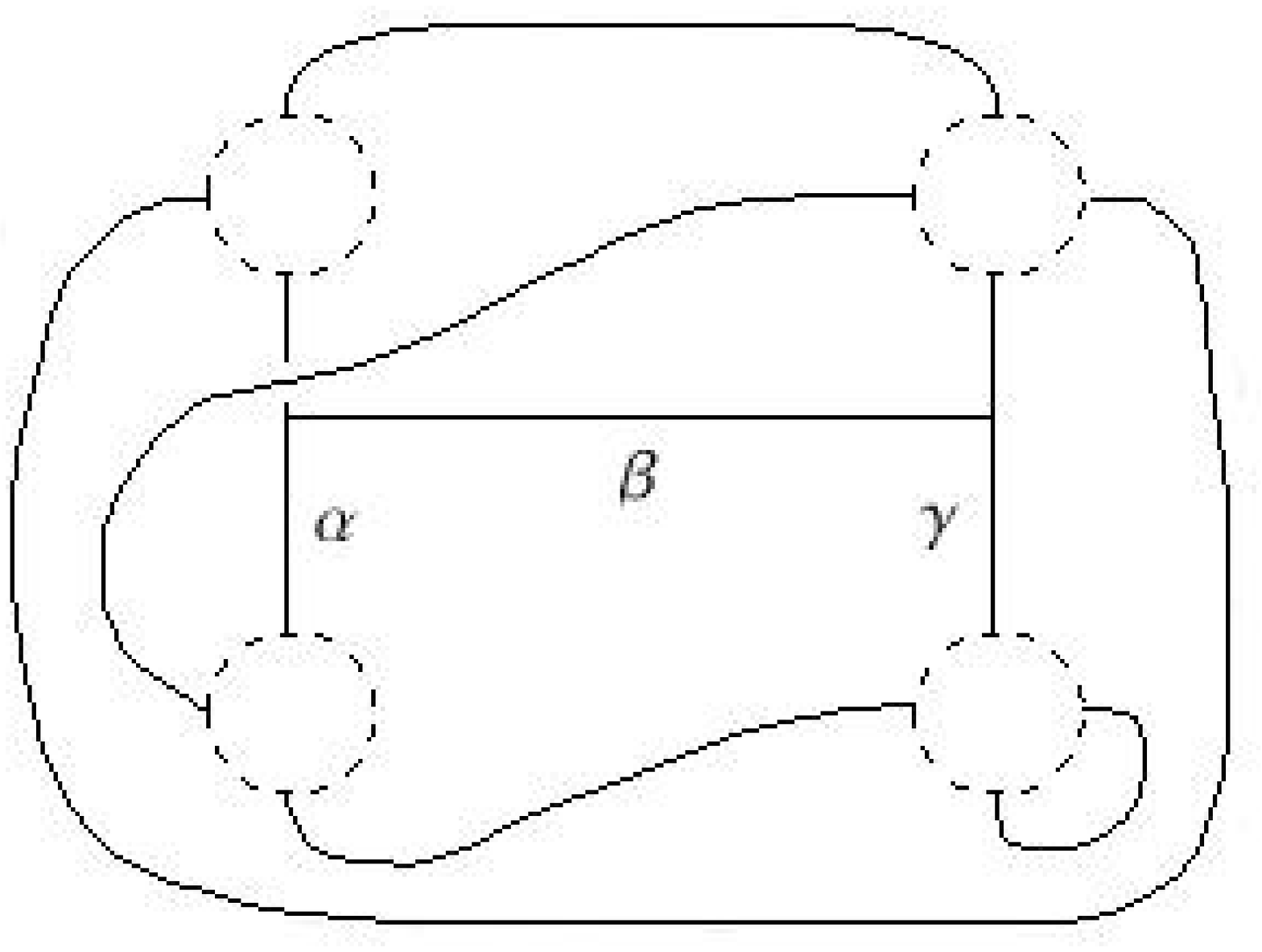}\end{minipage}\\ 
  =  \begin{minipage}{1.9in}\includegraphics[width=1.9in]{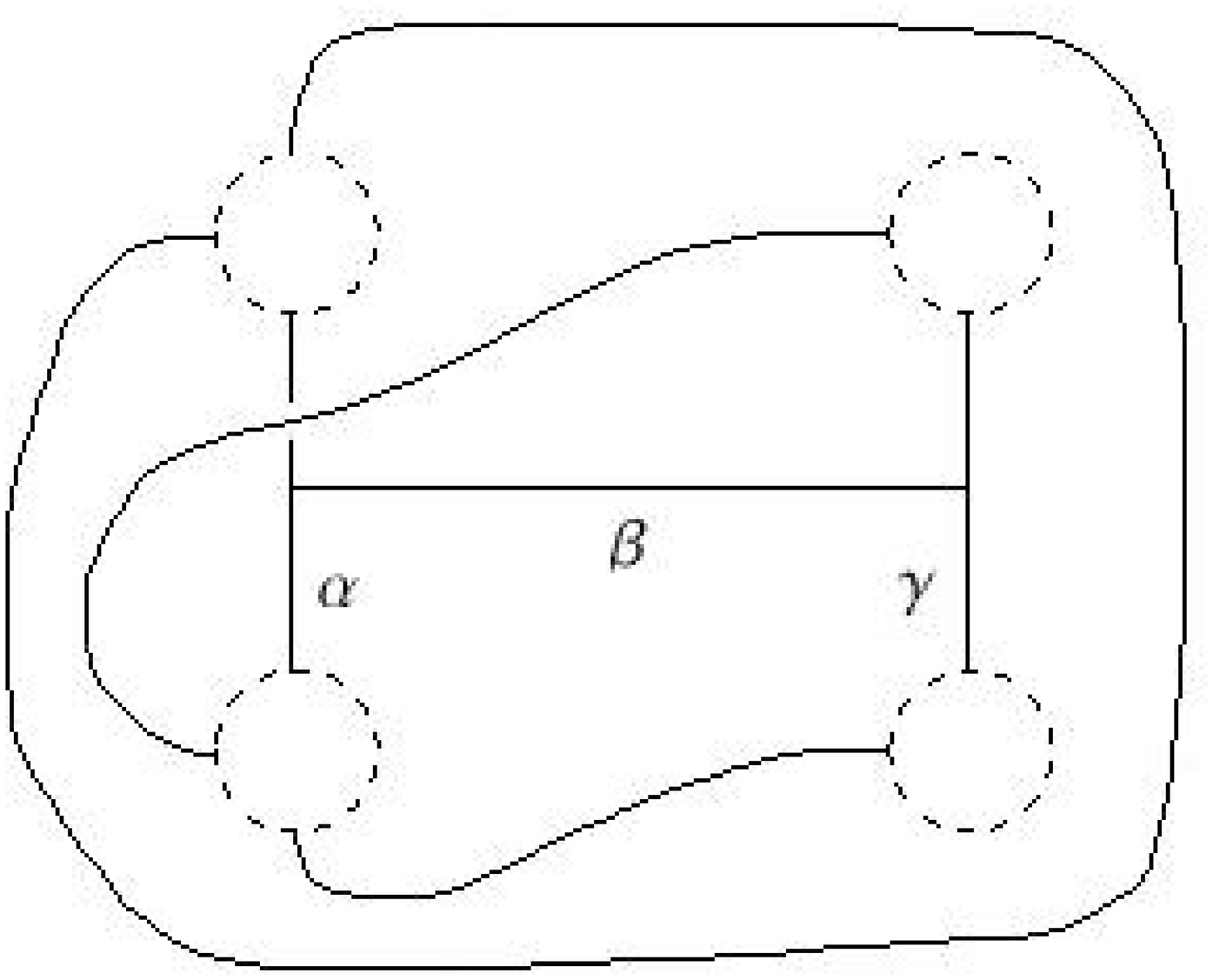}\end{minipage} 
  = &  \begin{minipage}{1.9in}\includegraphics[width=1.9in]{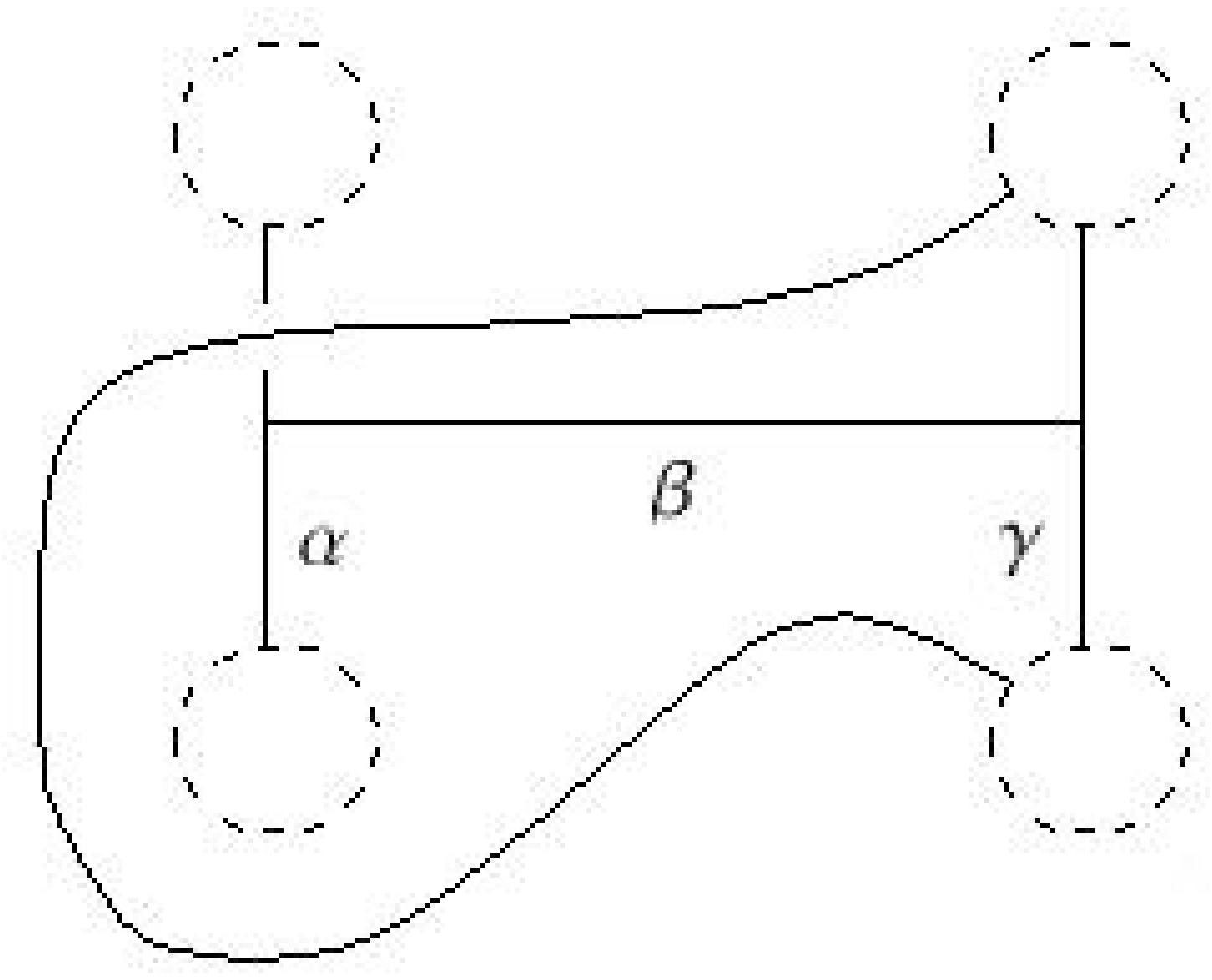}\end{minipage}
\end{eqnarray*}
\caption{
Slide 1;  Slide 2 
is obtained by rotating figure $180^\mathrm{o}$, leaving labels in place}
\label{fig:slide1}
\end{figure}

\begin{figure}
\begin{eqnarray*}
\begin{minipage}{1.9in}\includegraphics[width=1.9in]{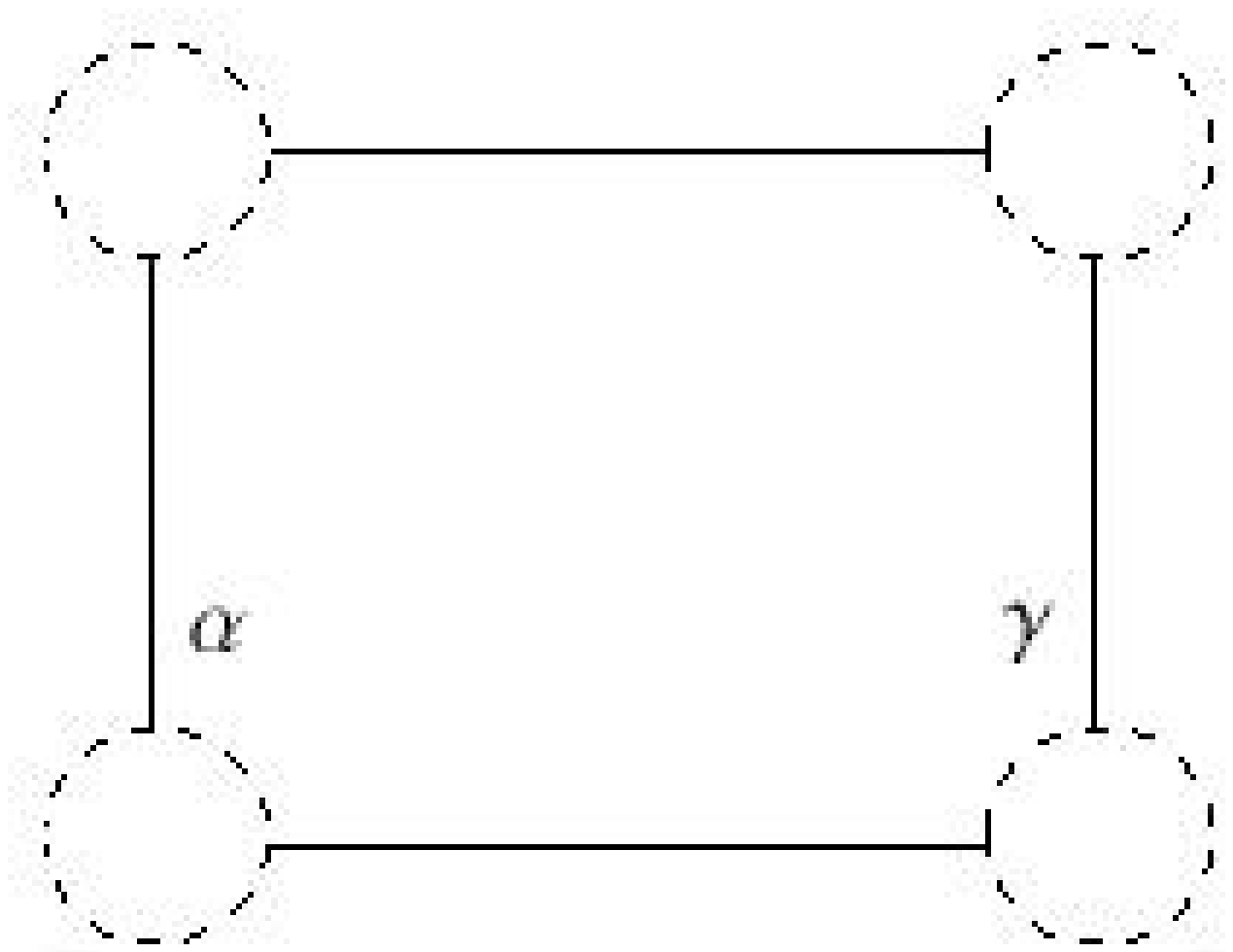}\end{minipage}
 & = & \begin{minipage}{1.9in}\includegraphics[width=1.9in]{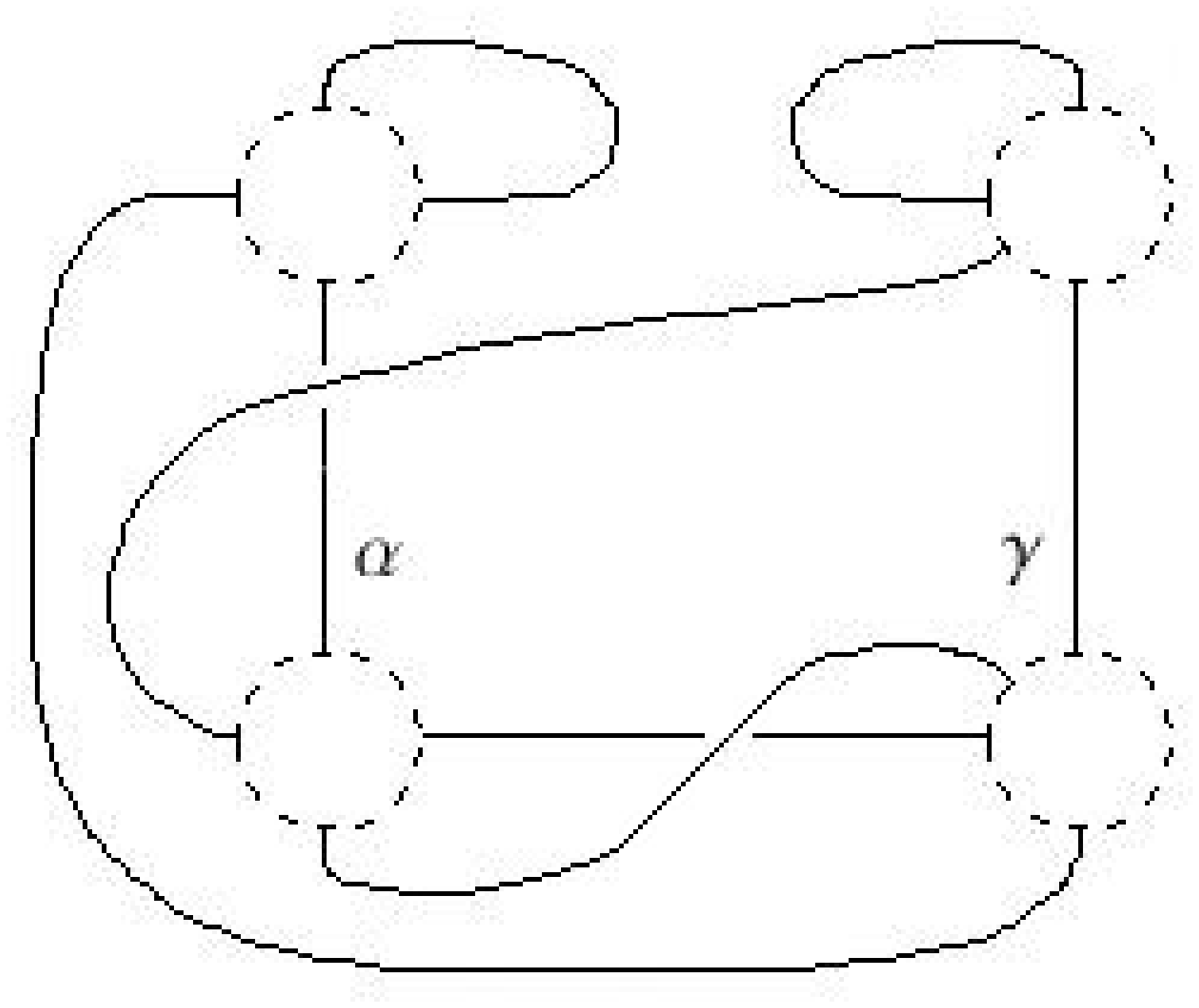}\end{minipage}\\       
 & = & \begin{minipage}{1.9in}\includegraphics[width=1.9in]{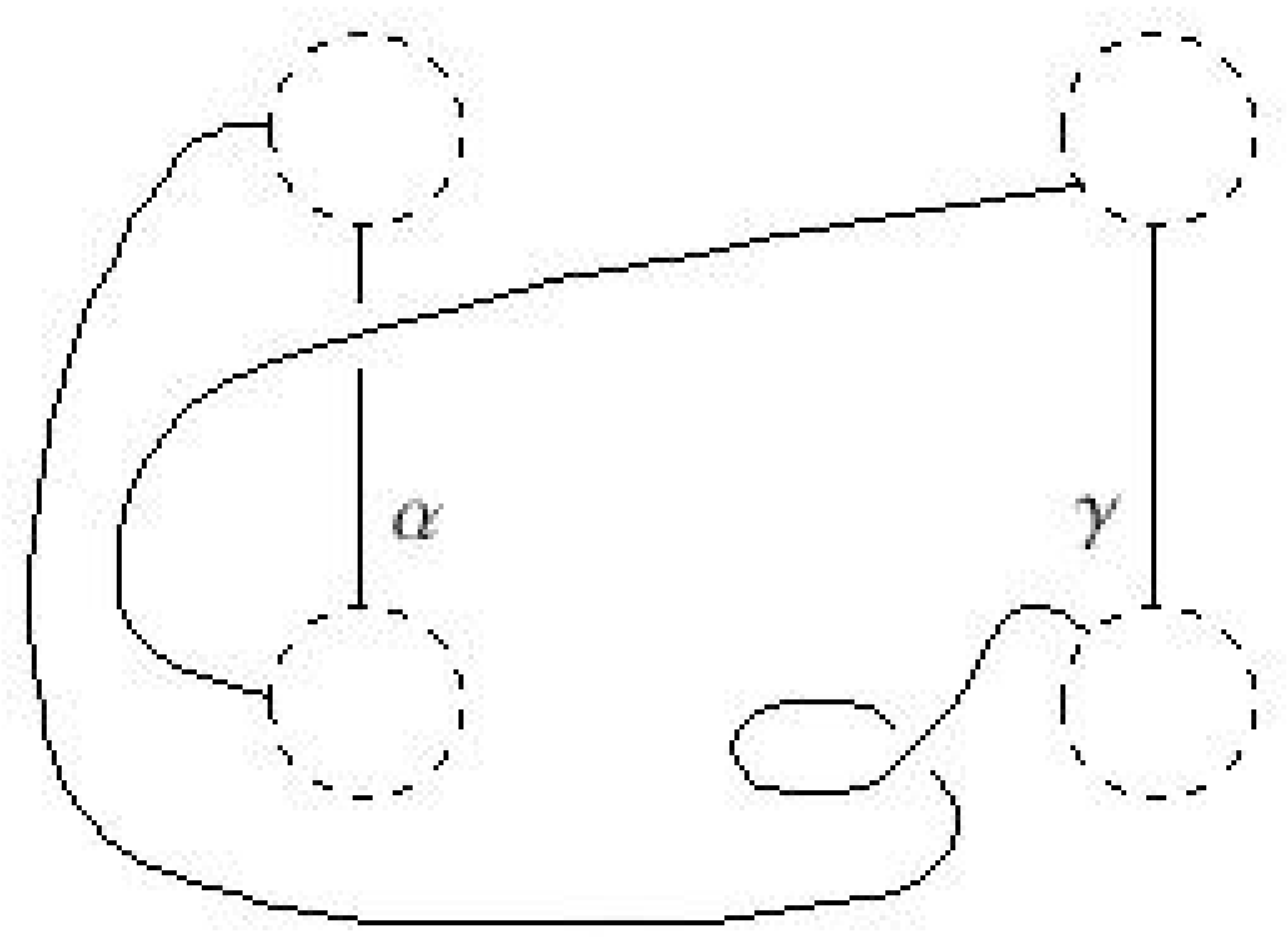}\end{minipage}
\end{eqnarray*}  
\caption{ Slide 3}
\label{fig:slide3}
\end{figure}

\begin{figure}
\begin{eqnarray*}
\begin{minipage}{2.2in}\includegraphics[width=2.2in]{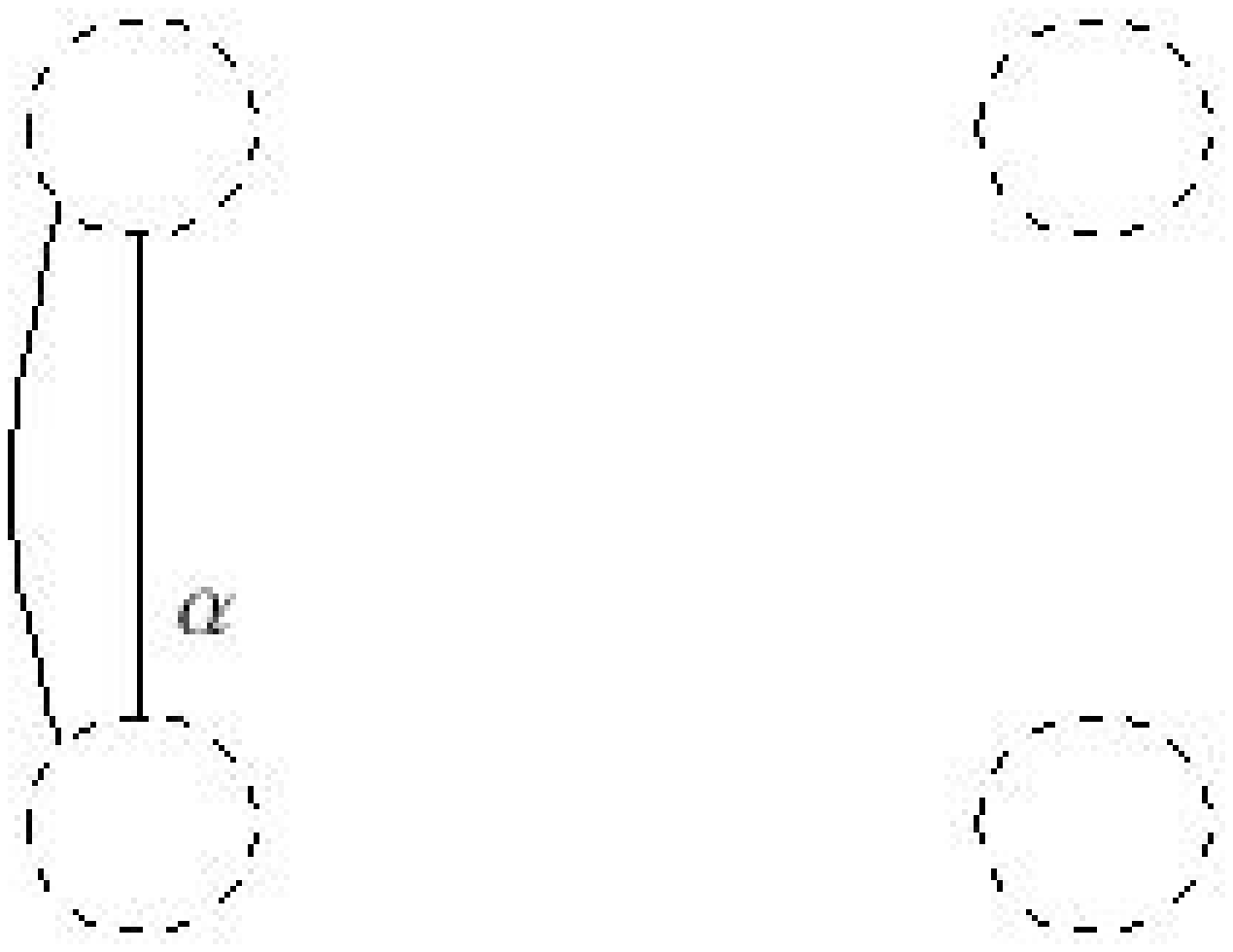}\end{minipage}
 & = & \begin{minipage}{2.2in}\includegraphics[width=2.2in]{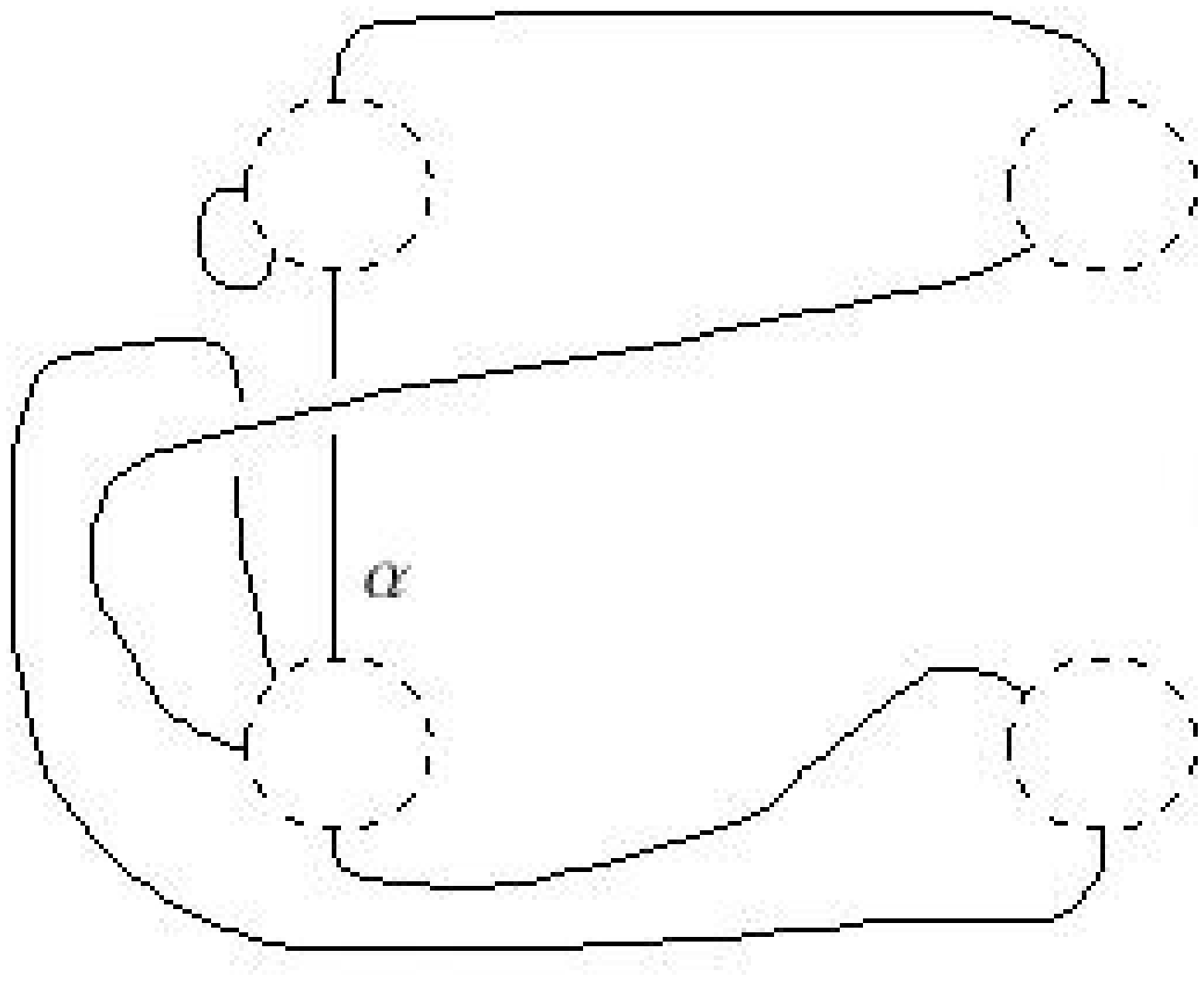}\end{minipage}\\        
 & = & \begin{minipage}{2.2in}\includegraphics[width=2.2in]{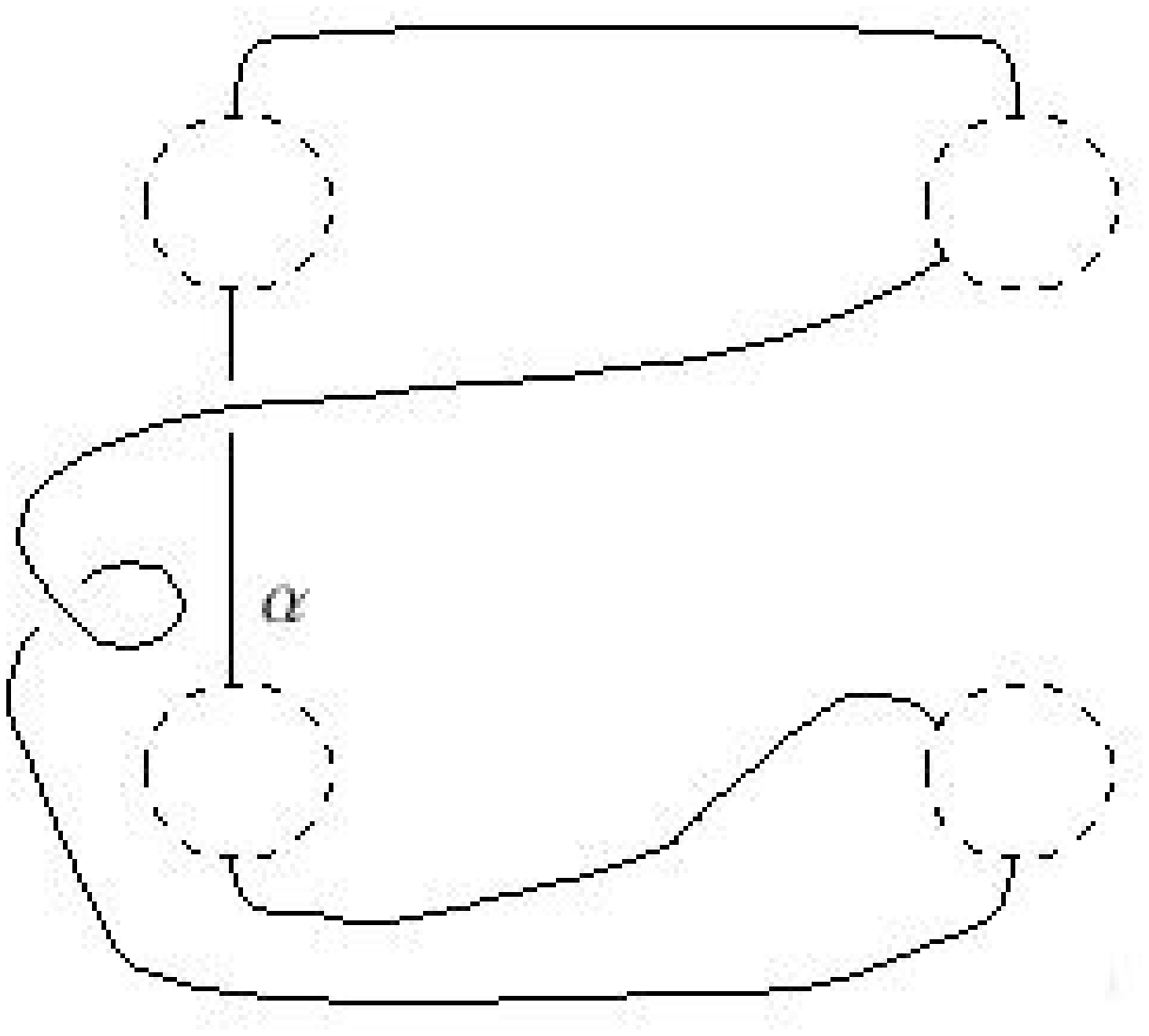}\end{minipage}
\end{eqnarray*}  
\caption{
Slide 4;  Slide 5 is obtained by rotating figure $180^\mathrm{o}$, leaving labels in place}
\label{fig:slide4}
\end{figure}

\begin{figure}
\begin{eqnarray*}
\begin{minipage}{2.2in}\includegraphics[width=2.2in]{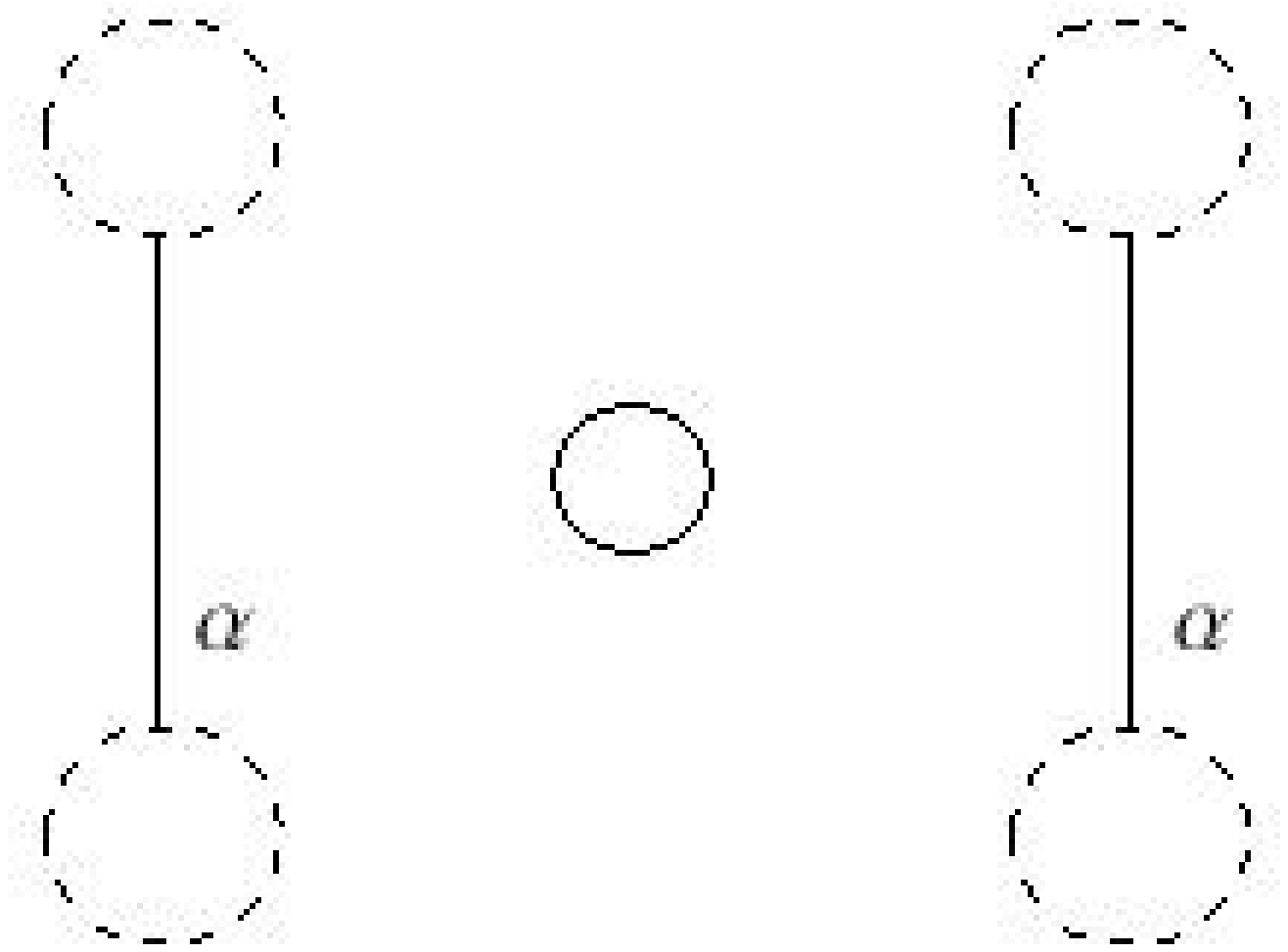}\end{minipage}
 & = &   \begin{minipage}{2.2in}\includegraphics[width=2.2in]{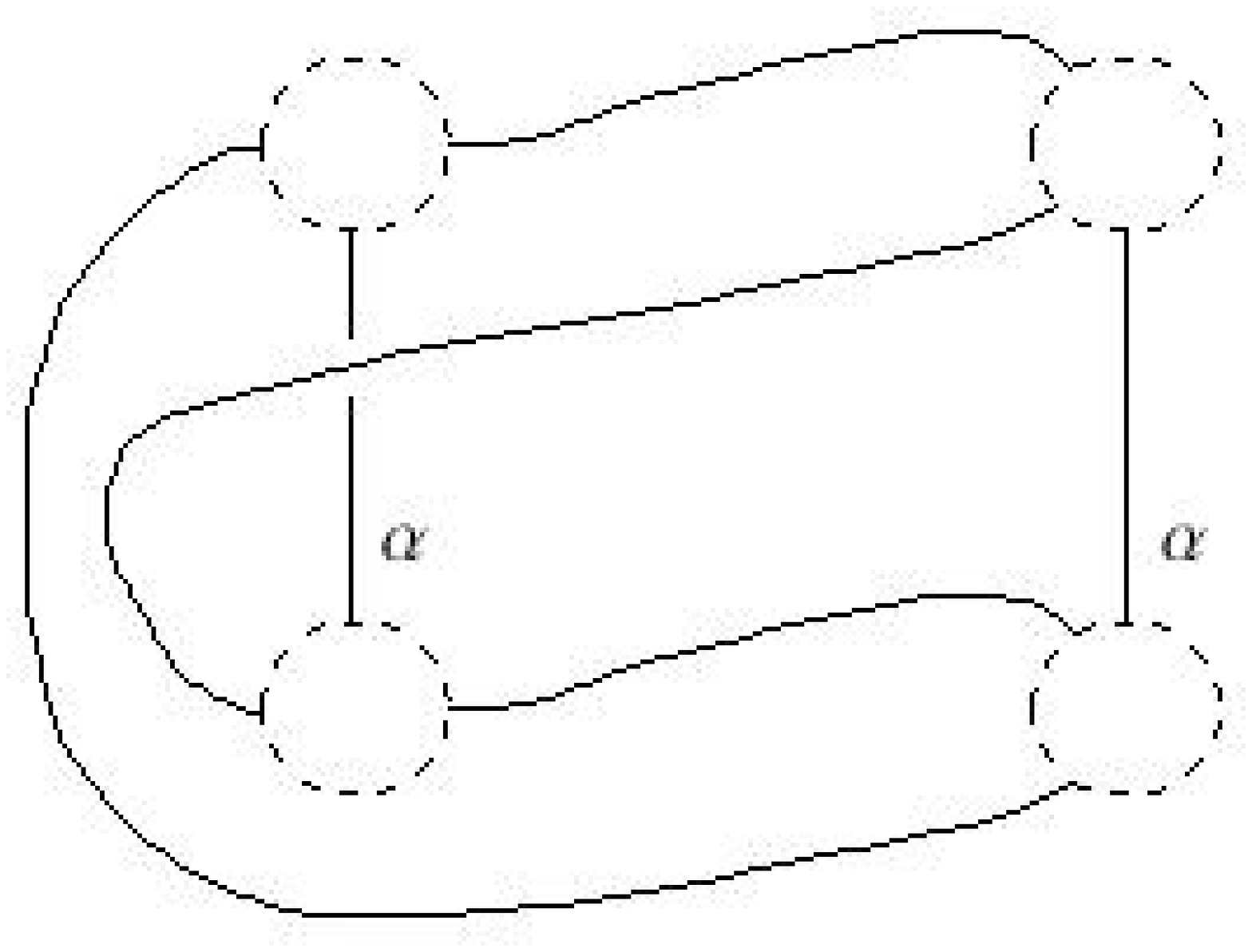}\end{minipage}
\end{eqnarray*}
\caption{Slide 6}
\label{fig:slide6}
\end{figure}

In this section, we show the skein module is generated by five elements.  To this end, we present in Figures \ref{fig:slide1},  \ref{fig:slide3}, \ref{fig:slide4}, and \ref{fig:slide6} six sets of 
slides over the attached 2-handles, yielding relations written in terms of our orthogonal basis.

Note that the diagrams of the second and fifth slides are just the diagrams of the first and fourth, respectively, rotated 180 degrees.  Hence, the computations of the coefficients for the first and fourth relations yield the coefficients for the second and fifth as well.

\begin{figure}
\begin{eqnarray*}
\begin{minipage}{2in}\includegraphics[width=2in]{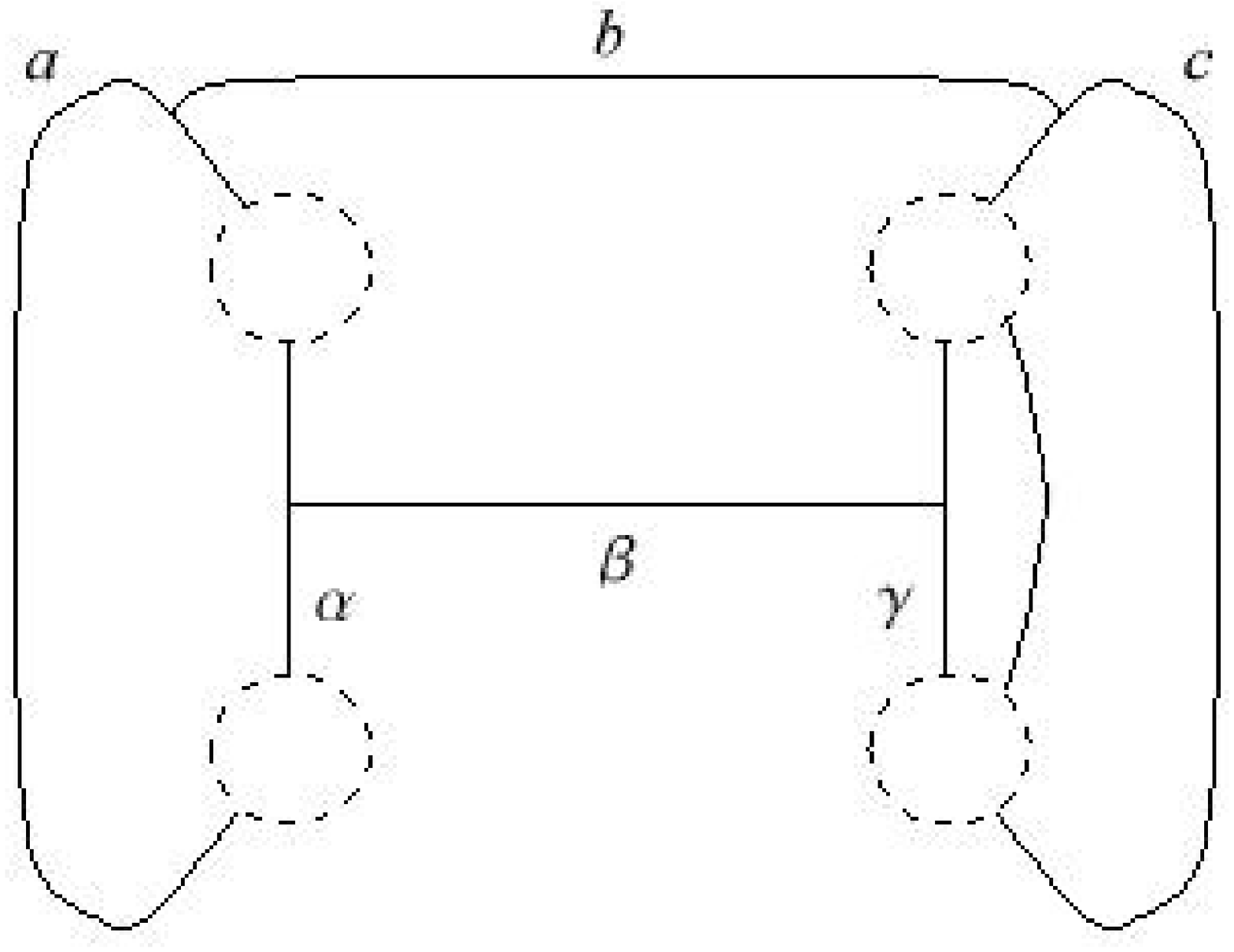}\end{minipage}
 & = & \frac{\delta^a_\alpha \;\delta^b_\beta \; \begin{minipage}{2in}\includegraphics[width=2in]{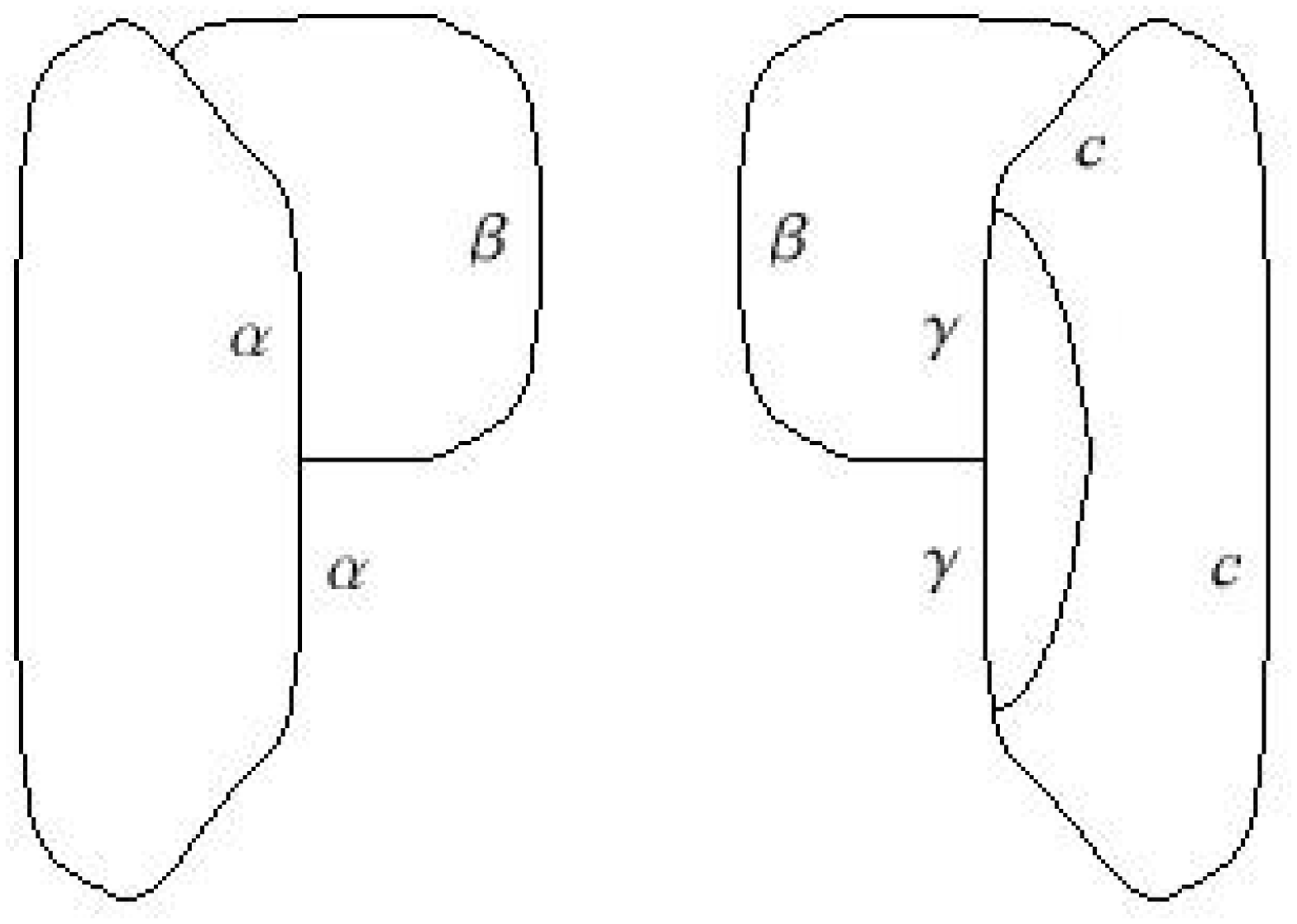}\end{minipage} }{\Delta_\alpha \;\Delta_\beta \;\theta(c,\gamma,1)}\\ \\
 & = &  \frac{\delta^a_\alpha \;\delta^b_\beta \;\theta(\alpha,\alpha,\beta) \tet c \gamma \gamma c \beta 1 }
 	{\Delta_\alpha \;\Delta_\beta \;\theta(c,\gamma,1)}
\end{eqnarray*}
\caption{
Computing the coefficients for the left side of Relation 1}
\label{fig:lhsrel1}
\end{figure}

\begin{figure}
\begin{eqnarray*}
\begin{minipage}{1.9in}\includegraphics[width=1.9in]{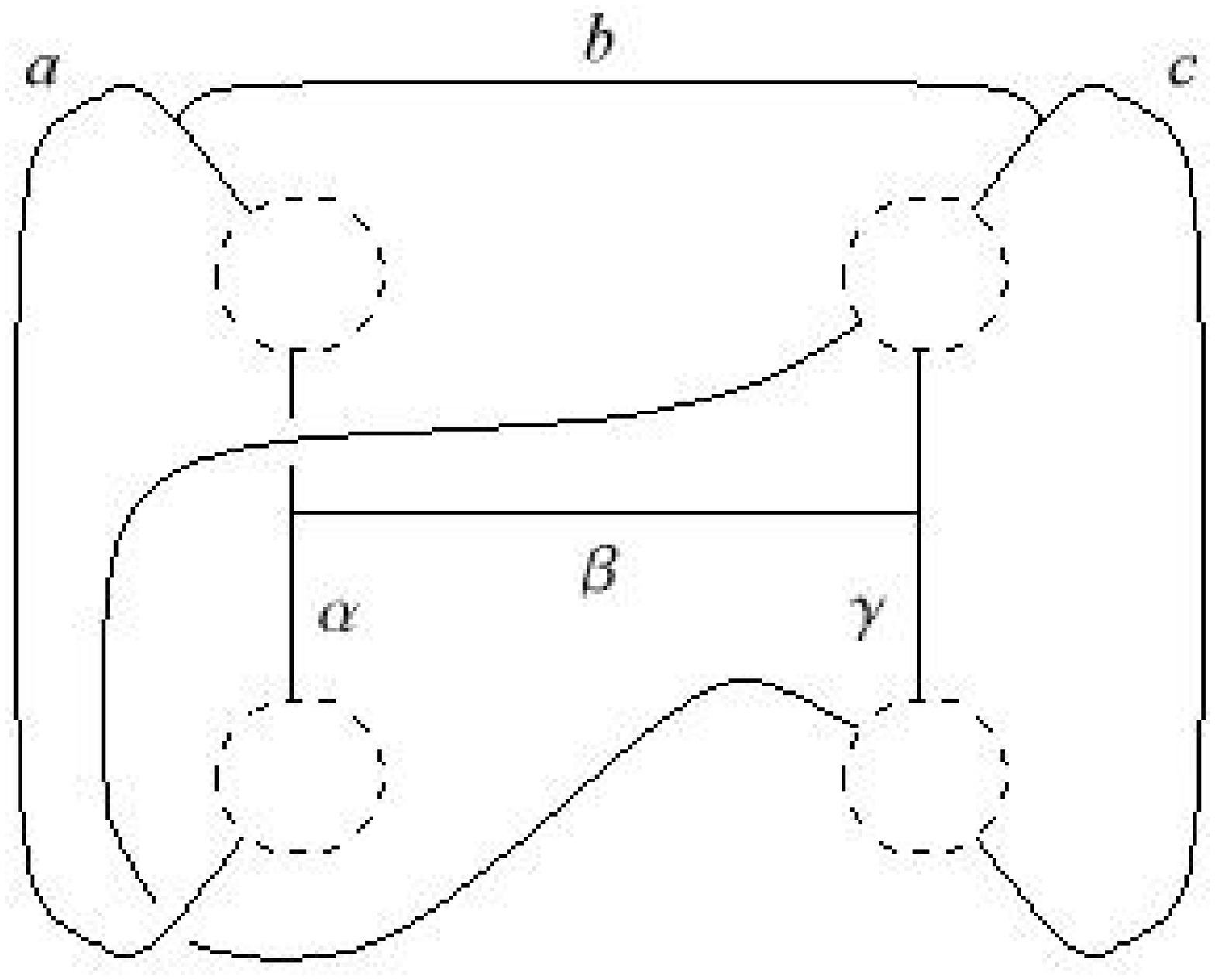}\end{minipage}
 & = & \sum_{r \;adm} \frac{\delta_\alpha^a \Delta_r  \begin{minipage}{1.9in}\includegraphics[width=1.9in]{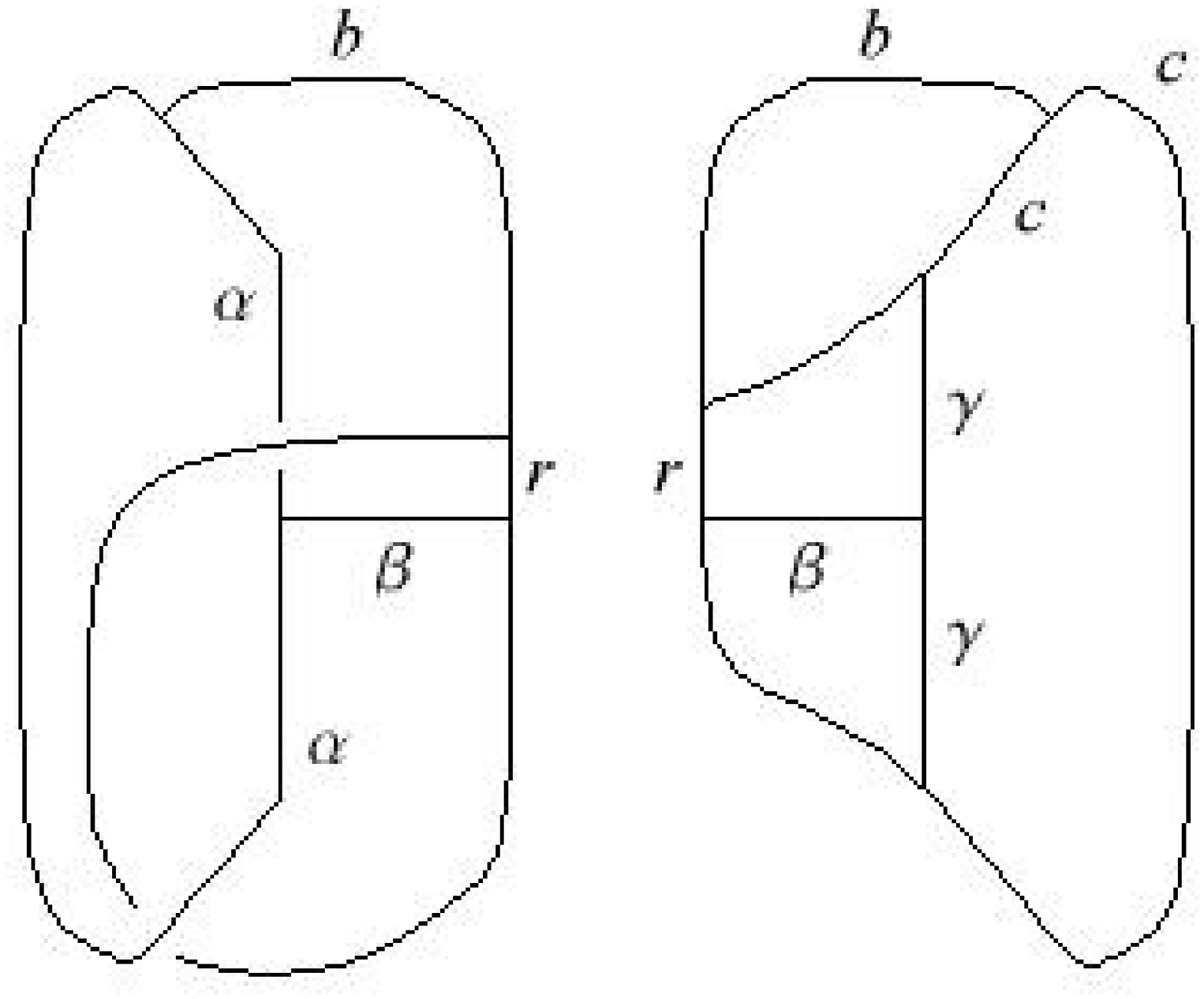}\end{minipage} }{\Delta_\alpha \;\theta(\beta,r,1) \;\theta(b,r,1) \;\theta(c,\gamma,1)}
\end{eqnarray*}
where
$$\begin{minipage}{1.175in}\includegraphics[width=1.175in]{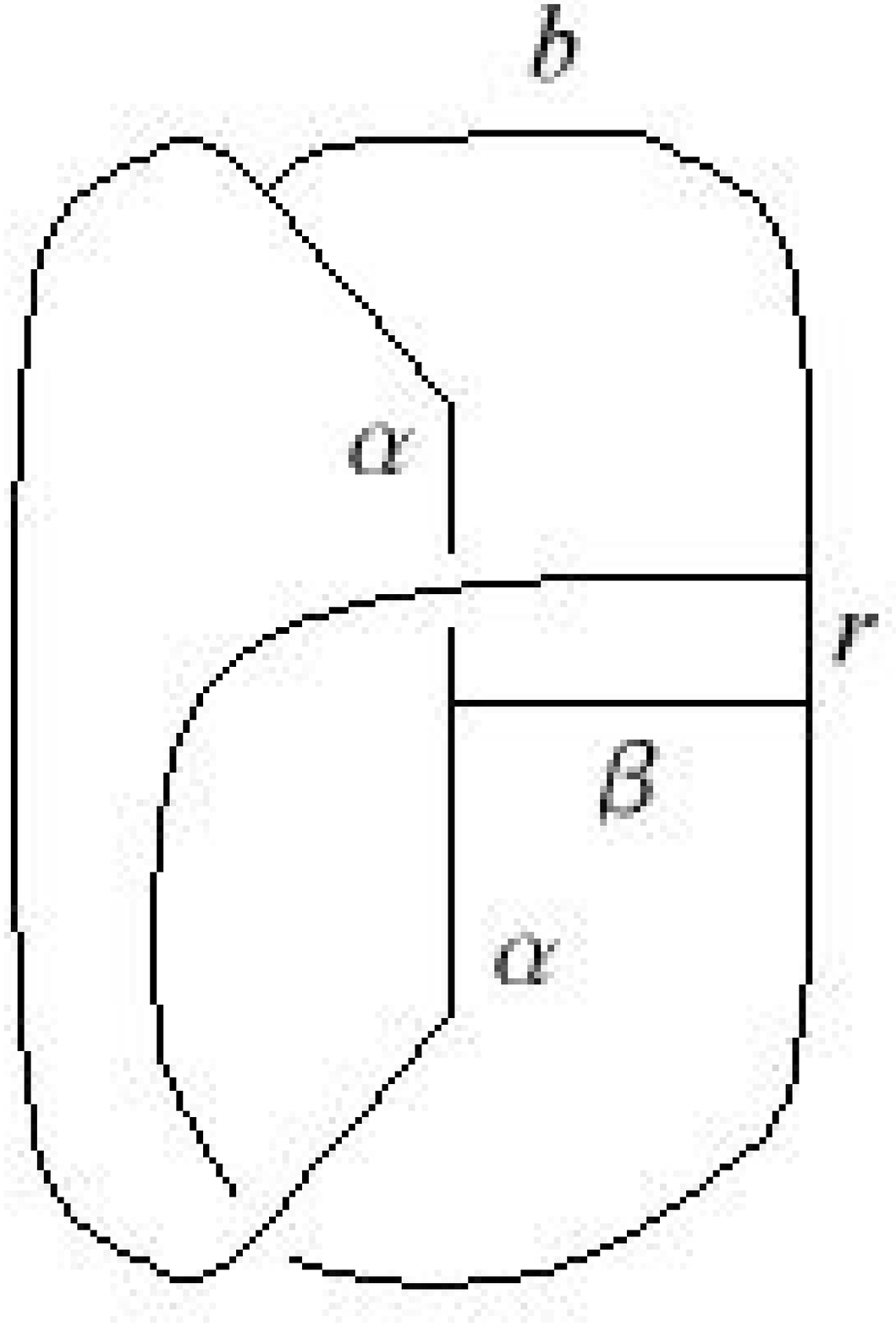}\end{minipage}
  = \sum_{i} \frac{
\Delta_i \;(\lambda_i^{\alpha \;1})^2 \tet r i \alpha b 1 \alpha \tet 1 i \alpha \beta r \alpha}
 {\lambda_r^{b \;1} \;\theta(\alpha,i,1) \;\theta(\alpha,r,i)}$$
and 
$$\begin{minipage}{1.175in}\includegraphics[width=1.175in]{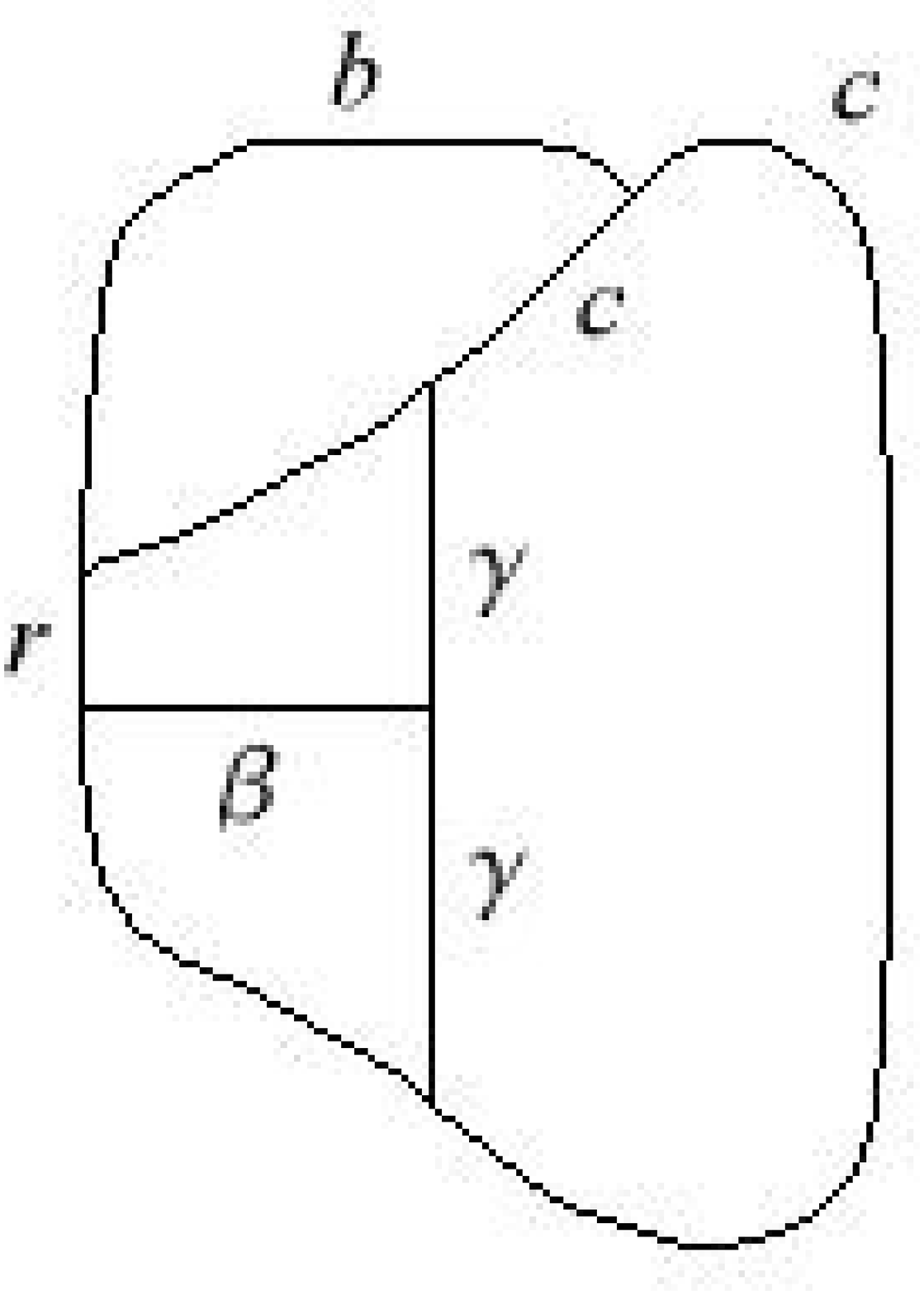}
\end{minipage}
 = \frac{\tet c r 1 c b \gamma \tet 1 r \gamma \gamma c \beta}{\theta(c,r,\gamma)}$$
\caption{
Computing the coefficients for the right side of Relation 1
}
\label{fig:rhsrel1}
\label{fig:rhsrel1_1}
\end{figure}

In Figures \ref{fig:lhsrel1} and \ref{fig:rhsrel1}, we demonstrate how the coefficients appearing in Relation 1
 are computed using recoupling theory and orthogonality.

In Figures \ref{fig:rel1}, \ref{fig:rel3}, \ref{fig:rel4}, and \ref{fig:rel6}, we present the relations obtained from slides 1,  3, 4, and 6.  In 
what follows, $<<a,b,c>>$ denotes $<(a,b,c),(a,b,c)>$ and terms in sums are taken to be zero if their denominators would be zero.
See \cite{Ha03} for a detailed derivation of the later relations from the slides.

\begin{figure}
\begin{flushleft}    $\sum_{a,b,c \;adm} \frac{\delta^a_\alpha \;\delta^b_\beta \;\theta(\alpha,\alpha,\beta) \tet c \gamma \gamma c \beta 1 }
 	{\Delta_\alpha \;\Delta_\beta \;\theta(c,\gamma,1)<<a,b,c>>} (a,b,c)$ \end{flushleft}

$= \sum_{a,b,c \;adm} \sum_{r} \frac{\delta_\alpha^a \Delta_r \tet c r 1 c b \gamma \tet 1 r \gamma \gamma c \beta  \sum_i Di }{\Delta_\alpha \;\theta(\beta,r,1) \;\theta(b,r,1) \;\theta(c,\gamma,1) {\theta(c,r,\gamma)} <<a,b,c>> } (a,b,c)$

\begin{flushleft} where $D_i= \frac{
\Delta_i \;(\lambda_i^{\alpha \;1})^2 \tet r i \alpha b 1 \alpha \tet 1 i \alpha \beta r \alpha}
 {\lambda_r^{b \;1} \;\theta(\alpha,i,1) \;\theta(\alpha,r,i)}$. \end{flushleft}

\caption{Relation 1}
\label{fig:rel1}
\end{figure}

\begin{figure}

\begin{flushleft} $\sum_{a,b,c \;adm} \frac{\tet 1 a a 1 b \alpha \tet c 1 1 c b \gamma}
 	{\theta(a,\alpha,1) \;\theta(b,1,1) \;\theta(c,\gamma,1)  <<a,b,c>>} (a,b,c)$ \end{flushleft}

$= \sum_{a,b,c \;adm} \frac{-A^3 \tet c 1 1 c b \gamma \sum_i D_i}{\lambda_a^{\alpha \;1} \theta(a,\alpha,1) \;\theta(b,1,1) \;\theta(c,\gamma,1) <<a,b,c>>} (a,b,c),$

\begin{flushleft} where $D_i =  \frac{\Delta_i \;\lambda_i^{a \;1} \tet i a a \alpha b 1 \tet \alpha 1 1 i b a}{\theta(a,i,1) \;\theta(b,\alpha,i)}$. \end{flushleft}

\caption{Relation 3}
\label{fig:rel3}
\end{figure}

\begin{figure}

\begin{flushleft} $\sum_{a,b,c \;adm} \frac{ \delta(a,\alpha,1) \;\delta^b_0 \;\delta^c_0}{<<a,b,c>>} \;(a,b,c)$ \end{flushleft}

$=\sum_{a,b,c \;adm} \sum_{r} \frac{-A^{-3}
\;\Delta_r \tet c r 1 c b 1 \tet a a r 1 \alpha b \sum_i D_i}{
\theta(a,\alpha,1) \;\theta(c,1,1) \;\theta(b,r,1) \;\theta(r,c,1) \;\theta(a,r,\alpha) <<a,b,c>>}  (a,b,c),$

\begin{flushleft} where $D_i = \frac{\Delta_i \;(\lambda_i^{a \;1})^2 \tet 1 i \alpha 1 c a \tet c i a r 1 \alpha}{\lambda_c^{r \;1} \;\theta(a,i,1) \;\theta(c,\alpha,i)}$. \end{flushleft}

\caption{Relation 4}
\label{fig:rel4}
\end{figure}

\begin{figure}
\bigskip

$$\sum_{a,b,c \;adm} \frac{\delta^a_\alpha \;\delta^b_0 \;\delta^c_\alpha \;\Delta_1}{<<a,b,c>>} \;(a,b,c) = \sum_{a,b,c \;adm}\sum_{p,q,q',r} \frac{C_{p,q,q',r}}{<<a,b,c>>} \;(a,b,c),$$

\begin{flushleft}for \end{flushleft}

$C_{p,q,q',r} = \frac{\Delta_p \;\Delta_q \;\Delta_{q'} \;\Delta_r \tet 1 c \alpha 1 {q'} r \tet {q'} c r 1 q \alpha \tet q b c r c 1}{\theta(\alpha,1,p) \;\theta(a,1,p) \;\theta(q',1,1) \;\theta(q,q',1) \;\theta(b,q,1) \;\theta(\alpha,r,1) \theta(c,r,1) \;\theta(c,\alpha,q') \;\theta(c,q,r)} \sum_i D_i,$

\begin{flushleft} where $D_i = \frac{\Delta_i \tet p 1 i 1 a \alpha }{\lambda_i^{\alpha \;1} \;\theta(\alpha,i,1) \;\theta(a,i,1) } \sum_j E_j,$ \end{flushleft}

\begin{flushleft}and $E_j = \frac{\Delta_j \;\lambda_j^{a \;1} \tet 1 {q'} i a j 1 \tet {q'} q \alpha i j 1 \tet q b p \alpha j 1 \tet j a a p b 1}{\theta(a,j,1) \;\theta(i,j,q') \;\theta(\alpha,j,q) \;\theta(p,b,j)}$. \end{flushleft}

\caption{Relation 6}
\label{fig:rel6}
\end{figure}

Let $r_i$ = (right-hand side of relation $i$) - (left-hand-side of relation $i$).

Note that, due to admissibility conditions, $(a,b,c)$ can only appear in $r_1$ if $a = \alpha, b \in \{\beta - 2, \beta, \beta + 2\}$, and $c = \gamma \pm 1$.  See Figure 
\ref{fig:rel1}.

$(a,b,c)$ can only appear in $r_2$ if $a = \alpha \pm 1, b \in \{\beta - 2, \beta, \beta + 2\}$, and $c = \gamma$. 

$(a,b,c)$ can only appear in $r_3$ if $a = \alpha \pm 1, b \in \{0,2\}$, and $c = \gamma \pm 1$.  See Figure \ref{fig:rel3}.
   
$(a,b,c)$ can only appear in $r_4$ if $a = \alpha \pm 1, b \in \{0,2,4\},$ and $c \in \{0,2\}$.  See Figure \ref{fig:rel4}.

$(a,b,c)$ can only appear in $r_5$ if $a \in \{0,2\}, b \in \{0,2,4\},$ and $c = \gamma \pm 1$.  

$(a,b,c)$ can only appear in $r_6$ if $a \in \{\alpha - 2, \alpha, \alpha + 2\}, b \in \{0,2,4\}$, and $c \in \{\alpha - 2, \alpha, \alpha + 2\}$.  See Figure 
\ref{fig:rel6}.

\begin{propo}
$S(M)$ is spanned by $(0,0,0)$, $(0,0,1)$, $(1,0,0)$, $(1,0,1)$, and $(0,0,2)$.
\end{propo}

\noindent {\em Proof.}  We proceed by induction on $(a,b,c)$.  Using the six sets of relations, we can rewrite $(x,y,z)$ as a linear combination of terms appearing earlier in the 
following ordering:

$(a',b',c') > (a,b,c)$ 
\begin{itemize}
	\item if $m' = \hbox{max}(a',c') > m = \hbox{max}(a',c')$, or 
	\item if $m' = m$ and $a' > a$, or 
	\item if $m' = m$, $a' = a$, and $c' > c$, or 
	\item if $m' = m$, $a' = a$, $c' = c$, and $b' > b$.
\end{itemize}

The proof splits into five cases:
\begin{enumerate}
	\item $x\geq1$, $y\geq2$, and $z\geq1$,
	\item $x\geq1$, $y=0$, $z\geq1$, and $z\neq x$,
	\item $x\geq2$ and $y=z=0$,
	\item $x=y=0$ and $z>2$, and 
	\item $x\geq2$, $y=0$, and $z=x$.
\end{enumerate}

After these have been rewritten, only $(0,0,0)$, $(0,0,1)$, $(1,0,0)$, $(1,0,1)$, and $(0,0,2)$ remain.

Case 1:  $x\geq1$, $y\geq2$, and $z\geq1$.

If we let $(\alpha, \beta, \gamma) = (x,y-2,z-1)$, then $(x,y,z)$ is the highest term appearing in relator $r_1$. Using Mathematica, we calculate that the coefficient of $(x,y,z)$ is 
$$ A^{-2-2x-y}(A^{2+2x}-A^y)(A^{2+2x}+A^y).$$ This is invertible, so $(x,y,z)$ is a linear combination of lesser terms.

Case 2:  $x\geq1$, $y=0$, $z\geq1$, and $z\neq x$.

Let $(\alpha, \beta, \gamma) = (x,0,z-1)$ in $r_1$, and let $(\alpha, \beta, \gamma) = (x-1,0,z)$ in $r_2$.  Then $(x,2,z)$ and $(x,0,z)$ are the two highest terms appearing in the relators.

Rearranging terms, we have
$$a_1 (x,2,z) + a_2 (x,0,z) = \hbox{lesser terms}$$ and 
$$b_1 (x,2,z) + b_2 (x,0,z) = \hbox{lesser terms},$$ 
where $a_i$ is the coefficient of the $i$th-highest term appearing in $r_1$ and $b_i$ is the coefficient of the $i$th-highest term appearing in $r_2$.

So, we can rewrite $(x,0,z)$ if $$\left|
\begin{array}[pos]{cc}
	a_1 & a_2 \\
	b_1 & b_2
\end{array}
\right|$$ is invertible.

Using Mathematica, we can see that this determinant is $-A^{-2-2x-2z} (-1+A^x)(1+A^x)(A^x-A^z)(-1+A^z)(1+A^z)(A^x+A^z),$ which is invertible when $x\geq1$, $z\geq1$, and $x \neq z$.

Case 3:  $x\geq2$, $y=z=0$.

Using $r_1$ with $(\alpha, \beta, \gamma) = (x,2,1)$, $r_2$ with $(\alpha, \beta, \gamma) = (x-1,0,2)$, $r_3$ with $(\alpha, 0, \gamma) = (x-1,0,1)$, and $r_4$ with $\alpha = x-1$, we obtain four relations with $(x,4,2), (x,2,2), (x,0,2)$, and $(x,0,0)$ appearing as the four highest terms.

Thus, we can rewrite $(x,0,0)$ if $$\left|
\begin{array}[pos]{cccc}
	a_1 & a_2 & a_3 & a_4 \\
	 0  & b_1 & b_2 &  0  \\
	 0  & c_1 & c_2 & c_3 \\
	d_1 & d_2 & d_3 & d_4
\end{array}
\right|$$ is invertible, where $a_i, b_i, c_i$, and $d_i$ are the coefficients of the $i$th-highest terms appearing in $r_1, r_2, r_3$, and $r_4$, respectively.

This determinant is $$-A^{-10-2x}(-1+A)(1+A)(1+A^2)(-A+A^x)^2(A+A^x)^2(A^2+A^{2x})^2,$$ which is invertible when $x \geq 2$.

Case 4:  $x=y=0$, $z>2$.

Using $r_1$ with $(\alpha, \beta, \gamma) = (2,2,z-1)$, $r_2$ with $(\alpha, \beta, \gamma) = (1,0,z)$, $r_3$ with $(\alpha, 0, \gamma) = (1,0,z-1)$, and $r_5$ with $\gamma = z-1$, we obtain four relations with $(2,4,z),(2,2,z),(2,0,z)$, and $(0,0,z)$ appearing as the four highest terms.
Note that this only holds for $z > 2$:  for $z = 2$, $(0,0,z)$ is no longer the fourth-highest term.
 
Thus, we can rewrite $(0,0,z)$ if $$\left|
\begin{array}[pos]{cccc}
	a_1 & a_2 & a_3 &  0  \\
	 0  & b_1 & b_2 & b_3 \\
	 0  & c_1 & c_2 & c_3 \\
	d_1 & d_2 & d_3 & d_4
\end{array}
\right|$$ is invertible, where $a_i, b_i, c_i$, and $d_i$ are the coefficients of the $i$th-highest terms appearing in $r_1, r_2, r_3$, and $r_5$, respectively.

This determinant is $-A^{-6-2z}(-1+A)(1+A)(1+A^2)(-1+A^z) (1+A^z)(-A+A^z)(A+A^z)(A^2+A^{2z}),$ which is invertible when $z > 2$.

Case 5:  $x\geq2$, $y=0$, and $z=x$.

Using $r_2$ with $(\alpha, \beta, \gamma) = (x-1,2,x)$, $r_3$ with $(\alpha, 0, \gamma) = (x-1,0,x-1)$, and $r_6$ with $\alpha = x-2$, we obtain three relations with $(x,4,x),(x,2,x)$, and $(x,0,x)$ appearing as the three highest terms.

Thus, we can rewrite $(x,0,x)$ if $$\left|
\begin{array}[pos]{ccc}
	a_1 & a_2 & a_3 \\
   0  & b_1 & b_2 \\
	c_1 & c_2 & c_3 
\end{array}
\right|$$ is invertible, where $a_i, b_i,$ and $c_i$ are the coefficients of the $i$th-highest terms appearing in $r_2, r_3$, and $r_6$, respectively.

This determinant is $A^{-4+2x}(-1+A^x)(1+A^x)$, which is invertible for $x\geq2$.

For all five cases, the determinants are indeed invertible, though their complexity compels us to use Mathematica for their evaluations.  See \cite{Ha03} for the code and the output. \qed

\noindent {\bf Remark} Note that all the determinants that we have computed are up to sign and powers of $A$
products of cyclotomic polynomials in $A.$  There is no obvious reason that these determinants should even be polynomials;  the entries of the matrices are generally quite complicated elements of $\cal{R}$.  It is this unexpected fact which allows us to make our computation over $\cal{R}.$

\section{Linear Independence}
\label{sec:LinearIndependence}

\begin{propo}
$(0,0,0)$, $(0,0,1)$, $(1,0,0)$, $(1,0,1)$, and $(0,0,2)$ are linearly independent in $S(M)$.
\end{propo}

\noindent {\em Proof.}  First, we recall that $S(M)$ is a direct sum of four submodules $S_1(M)$, $S_i(M)$, $S_j(M)$, and $S_k(M)$, and that the latter three are isomorphic.

Hence, our task is greatly simplified.  We only have to show that $(0,0,0)$ and $(0,0,2)$ are linearly independent, and that $(0,0,1)$ is nonzero.

\begin{defi}
The triple ${a,b,c}$ is said to be $r$-admissible if it is admissible and $a + b + c \leq 2r - 4$.
\end{defi}

The recoupling theory we have used in previous chapters works when $A$ is replaced by a primitive $2r$th root of unity, for odd $r>1$.  But we must replace ``admissible" with 
``$r$-admissible" and restrict our labels to the range {0, \ldots, r-2}:  our insistence on $r$-admissibility ensures that the fusion formula still makes sense. 

Let $A_r$ be the $2r$th root of unity $e^{\pi i / r}$ with $r>1$ odd, let $e_i$ denote the core of the solid torus, labelled $i$.  Let $\Omega_r = \sum_{i=0}^{\lfloor(r-3)/2\rfloor} \Delta_i e_i$, choose $\eta$ such that $\eta^2 <\Omega_r> = 1$, and let $\kappa$ be a root of unity such that $\kappa^6 = A_r^{-6-r(r+1)/2}$.  

\begin{defi}
For a framed link $K$ in a closed, connected, oriented 3-manifold $M$ described by surgery on a framed link $L \subset S^3$, we define the quantum invariant 
$$I_r(M,K) = \kappa^{3(b_-(L) - b_+(L))} <L(\eta \Omega_r) \cup K>,$$
where $b_+(L)$ and $b_-(L)$ are the numbers of the positive and negative eigenvalues of the linking matrix of $L$.
\end{defi}

We follow the notation of Masbaum and Roberts in \cite{MR97}.  See \cite{BHMV}, \cite{Li93}, \cite{KM}, \cite{RT91}, and \cite{Wi89} for the origins of this formula.

Note that for the quaternionic manifold $M$, with the surgery description $L$ presented in the introduction, $b_+(L)=b_-(L)$, and so $I_r(M,K) = \eta^2 <L(\Omega_r) \cup K>$.

\begin{propo} For odd $r>1$, 
$$(1-A_r^4) I_r(M) = \sum_{k=1}^{r-1} (-1)^k A_r^{2k^2+2k},$$ 
$$I_r(M,(0,0,1)) = (-1)^\frac{r-1}{2} \frac{A_r^{-2}}{A_r^2+1},$$ and  
$$A_r^4 I_r(M,(0,0,2)) = I_r(M) - 1.$$
\label{propo:invariant}
\end{propo}
\noindent {\em Proof.}  
$$\begin{minipage}{1.51in}\includegraphics[width=1.51in]{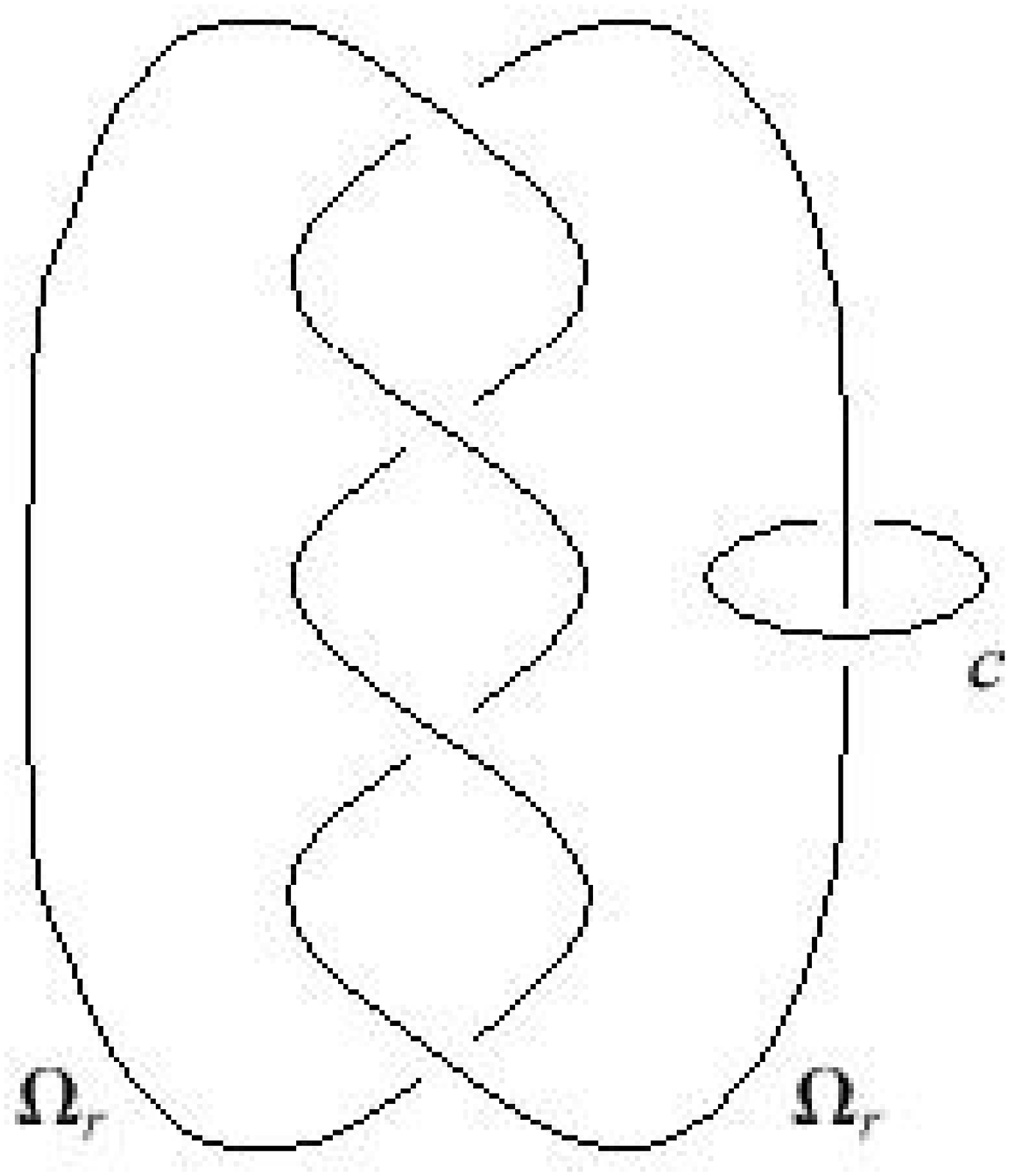}\end{minipage}
 = \sum_{i=0}^{(r-3)/2} (-A_r)^{i(i+2)} \Delta_i \begin{minipage}{1.6875in}\includegraphics[width=1.6875in]{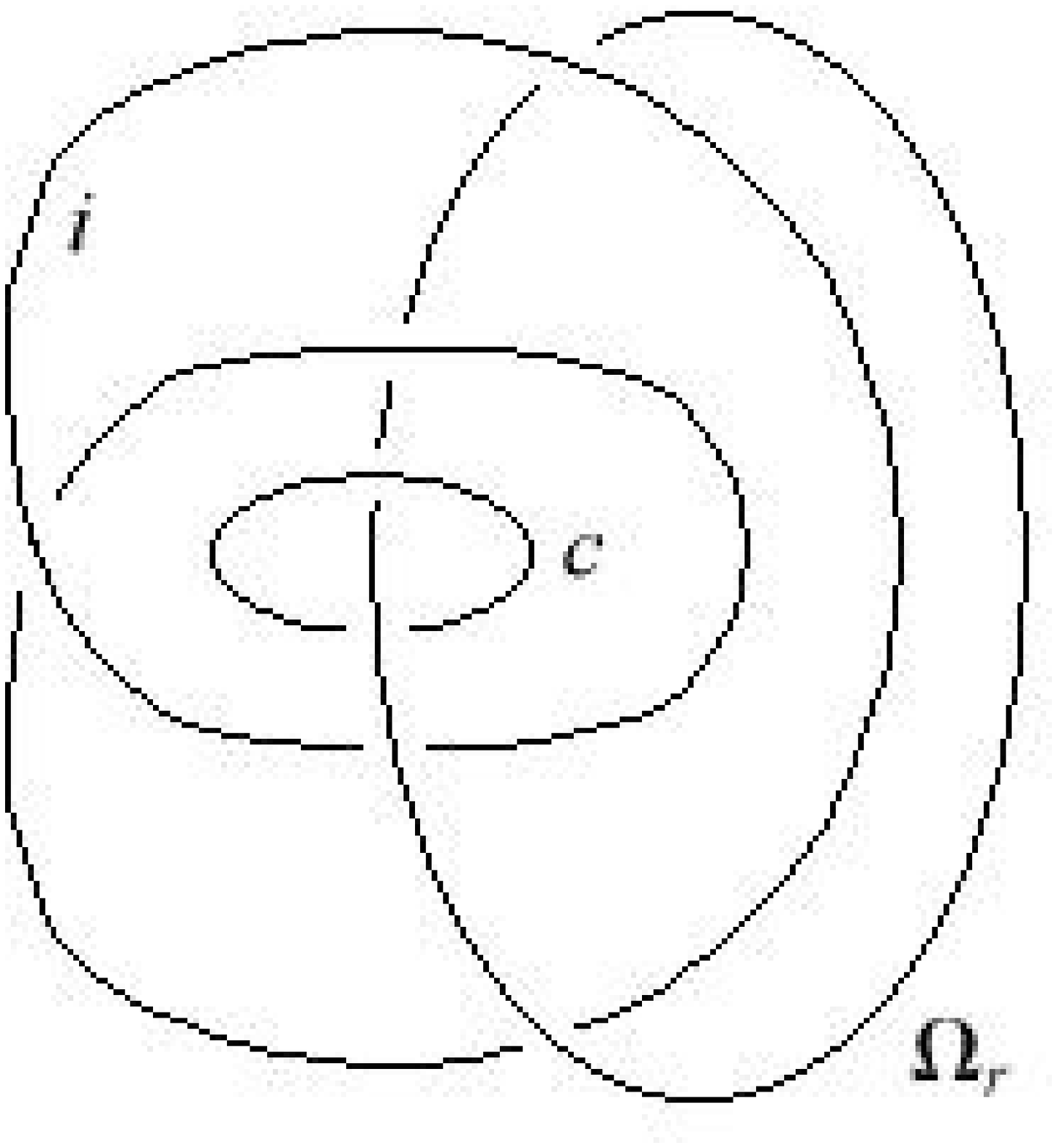}\end{minipage}$$
\begin{eqnarray*}
= \sum_{i=0}^{(r-3)/2} (-A_r)^{i(i+2)} \Delta_i \sum_{\alpha,\beta \;adm}
 \frac{\Delta_\alpha \Delta_\beta}{\theta(i,i,\alpha) \theta(\alpha,\beta,c)} \begin{minipage}{1.125in}\includegraphics[width=1.125in]{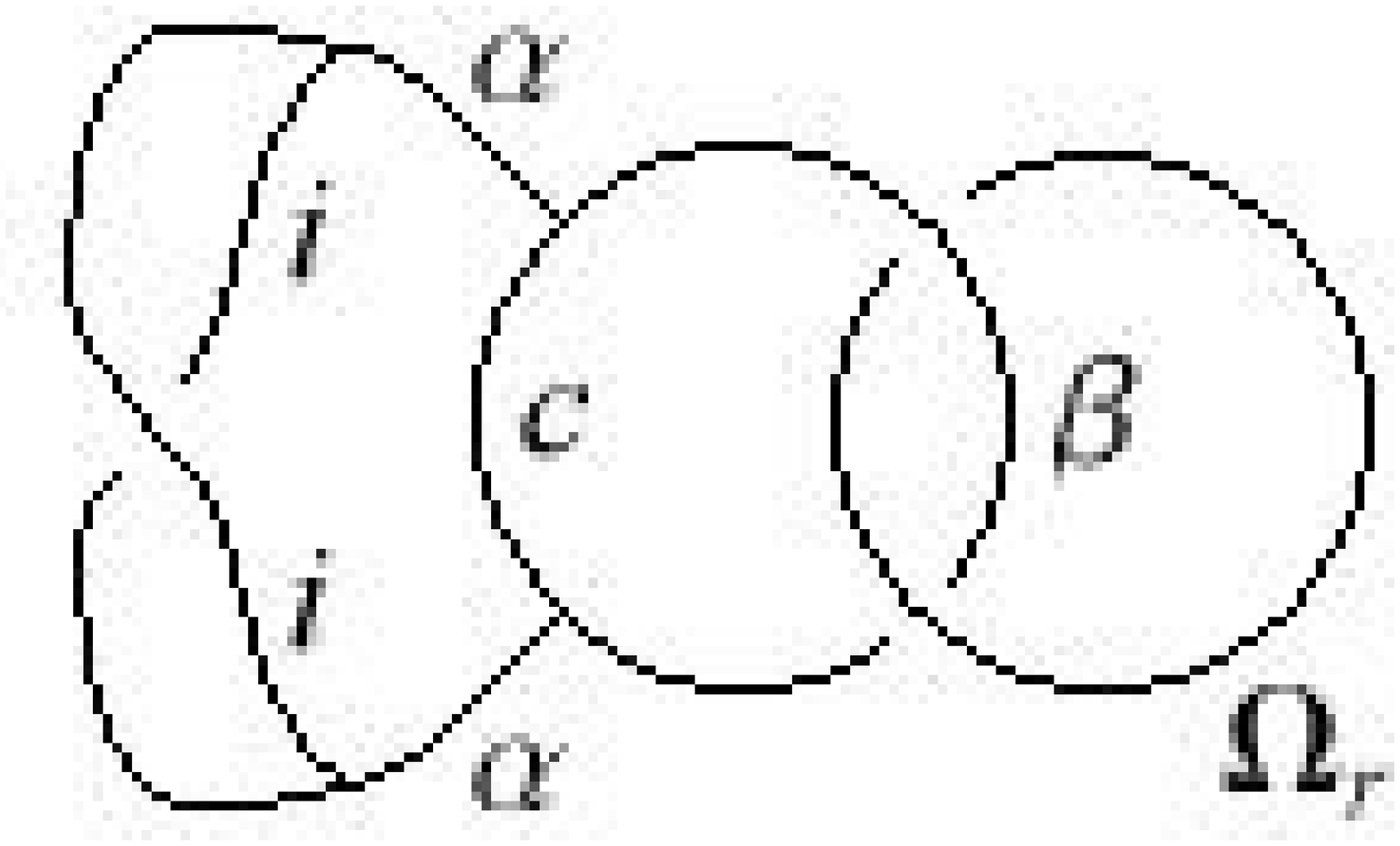}\end{minipage}\\
  = \sum_{i=0}^{(r-3)/2} (-A_r)^{i(i+2)} \Delta_i \sum_{\alpha,\beta \;as \;above}  
 \frac{\Delta_\beta \lambda_\alpha^{i \;i}}{\theta(\alpha,\beta,c)} \begin{minipage}{1.125in}\includegraphics[width=1.125in]{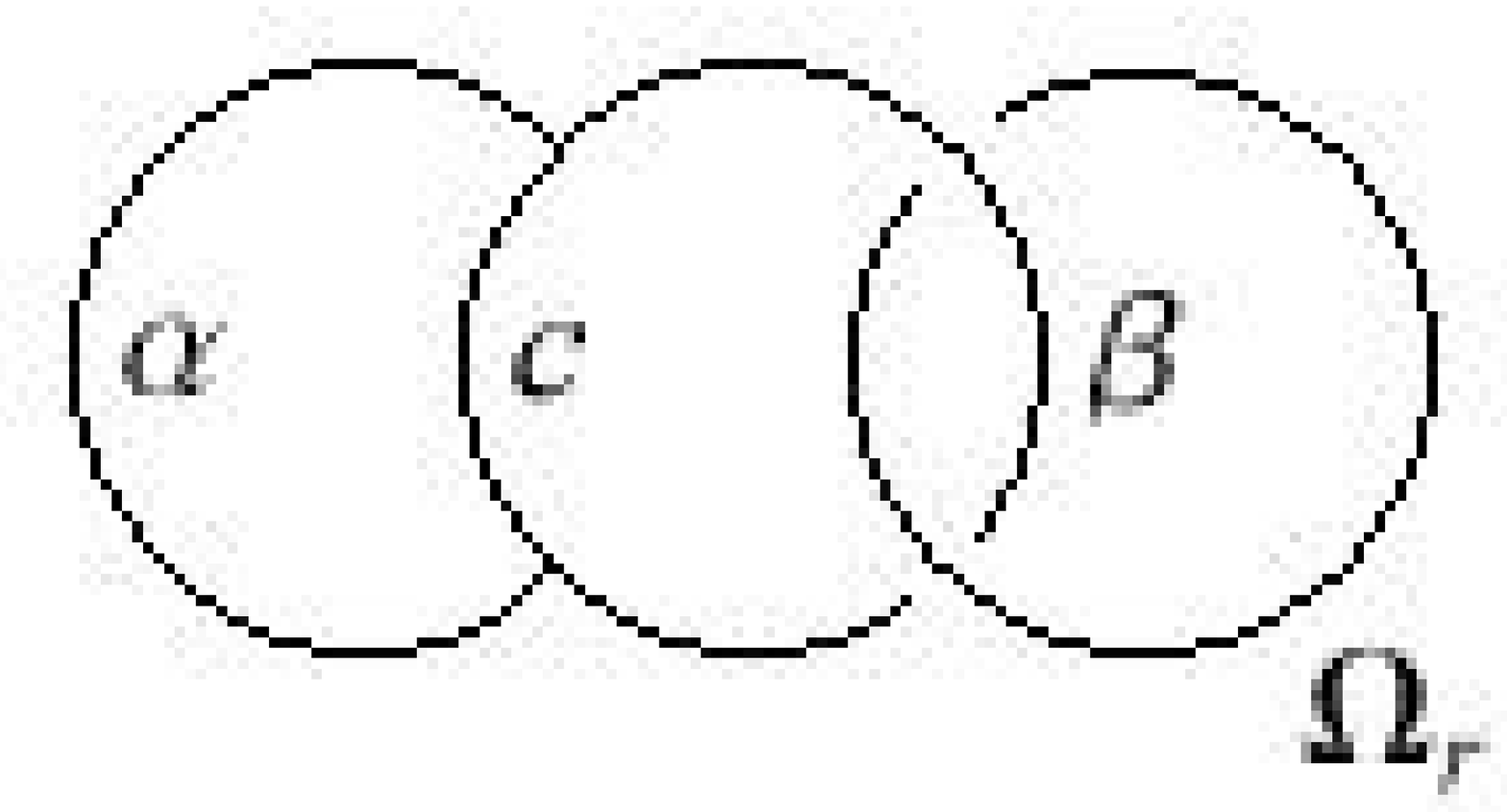}\end{minipage}\\
  = \sum_{i=0}^{(r-3)/2} (-A_r)^{i(i+2)} \Delta_i \sum_{\alpha,\beta \;as \;above} \lambda_\alpha^{i \;i}  \begin{minipage}{1.125in}\includegraphics[width=1.125in]{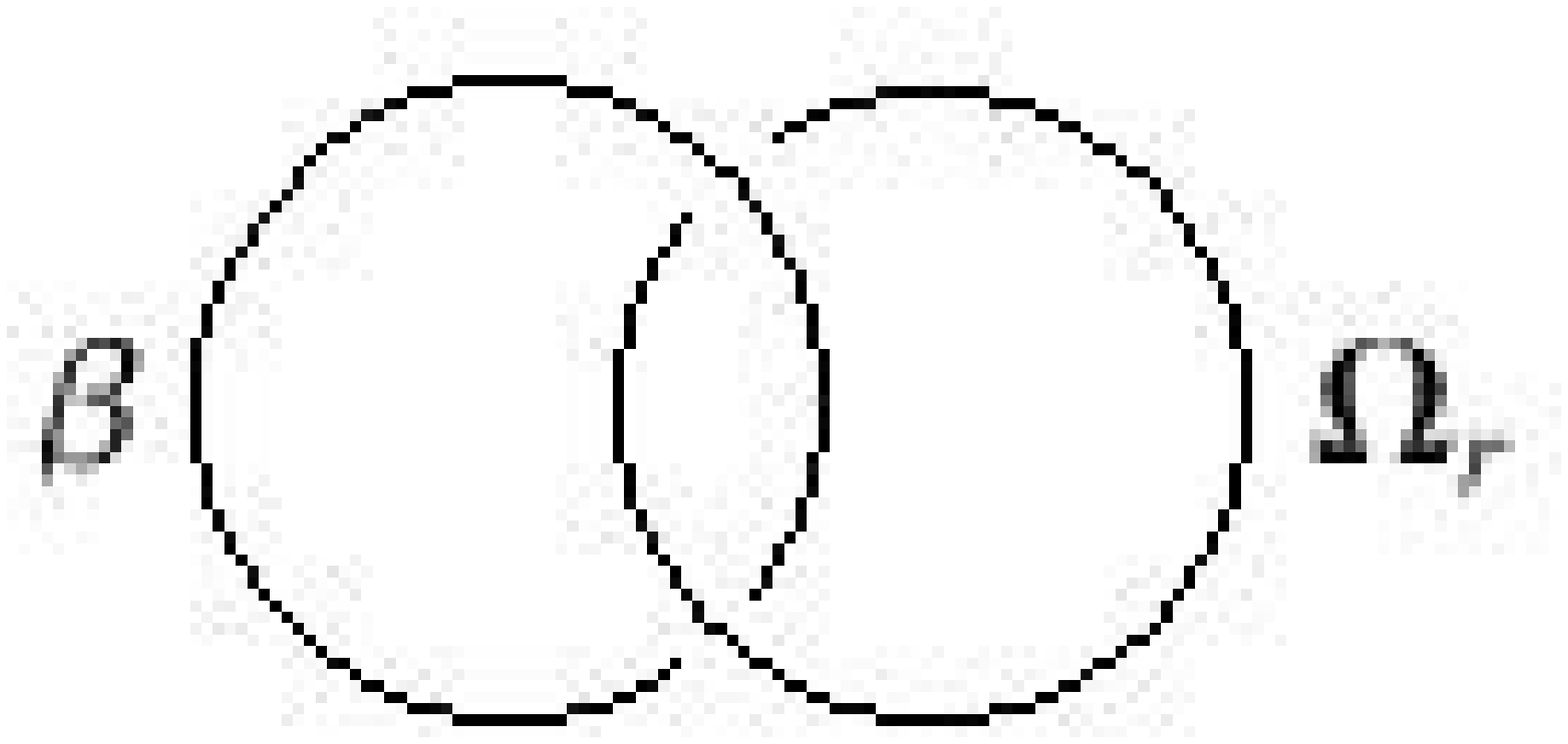}\end{minipage}
\end{eqnarray*}

By the encirclement lemma, 
\begin{minipage}{1.125in}\includegraphics[width=1.125in]{graphics/inv5}\end{minipage} $= \eta^{-2} \Delta_\beta$ if $\beta = 0$ or $\beta = r-2$, and $0$ otherwise.  See, for example, \cite{Li93}.  

Furthermore, the sums are restricted to $r$-admissible labellings, and so, $0 \leq \alpha \leq 2i$, $|\alpha - c| \leq \beta \leq \alpha + c$, and $\beta \equiv \alpha + c \mmod2$.  

Hence, a contribution to the sum can only occur when  $\alpha = 0, \beta = 0$ for $c = 0$, $\alpha = r-3, \beta = r-2, i = \frac{r-3}{2}$ for $c = 1$, and $\alpha = 2, \beta = 0$ for $c = 2$.

Hence, $$I_r(M) = \sum_{i=0}^{(r-3)/2} (-A_r)^{i(i+2)} \Delta_i \lambda_0^{i \;i},$$
$$I_r(M,(0,0,2)) = \sum_{i=1}^{(r-3)/2} (-A_r)^{i(i+2)} \Delta_i \lambda_2^{i \;i},$$ and 
$$I_r(M,(0,0,1)) = (-A_r)^{\frac{r-3}{2}(\frac{r-3}{2}+2)} \Delta_\frac{r-3}{2} \Delta_{r-2} \lambda_{r-3}^{\frac{r-3}{2} \;\frac{r-3}{2}} = (-1)^\frac{r-1}{2} \frac{A_r^{-2}}{A_r^2+1}.$$

Hence, $A_r^4 I_r(M,(0,0,2)) = I_r(M) - 1$, and letting $q=A_r^2$ and $k=i+1$, we obtain 
$$-q(q-q^{-1})I_r(M) = \sum_{k=1}^{(r-1)/2} (-1)^k (q^{k^2+k} - q^{k^2-k}).$$

Since $q^r=1$, $(-1)^{r-k} q^{(r-k)^2+(r-k)} = -(-1)^k q^{k^2-k}$. Hence, 

$$(1-q^2) I_r(M) = \sum_{k=1}^{r-1} (-1)^k q^{k^2+k}.$$  \qed

Before proceeding, we need more notation and a useful lemma:

\begin{nota}
Let $\zeta_N = e^{2 \pi i/N}$.  (Hence, $A_r = \zeta_{2r}$.)  Also, let $Q(A)$ denote the field of fractions  of $\Z[A].$ Elements of $Q(A)$ are called rational functions.
\end{nota}

\begin{lem}
If $F(A)$ is a rational function,
 then the imaginary part of $F(\zeta_{2r})$ cannot change sign infinitely often as $r \rightarrow \infty$.
\label{lem:rational}
\end{lem}

\noindent {\em Proof.}  
We have that $F(A^{-1}) \in Q(A)$ and so   $G(A) = F(A) - F(A^{-1}) \in Q(A)$. For $z$ on the unit circle, 
$\im F(z) = \frac{1}{2i}(F(z) - \overline {F(z)})= \frac{1}{2i}(F(z) -{F(\bar z)})=\frac{1}{2i} G(z)$.  
As $G(z)$ has at most finitely many zeros and poles,
$F(z)$  is zero or undefined for at most finitely many $z$ on the unit circle.
Thus, we can choose $N$ such that $F(z)$ has no poles  or zeros for $|z| = 1$ with $\mathrm{ Arg}z \in (0,\frac{2\pi}{N})$.  

Let $h(x) = \im F(e^{2 \pi i x})$, with domain $(0,\frac{1}{N})$.   Then $h$ is a  real-valued 
continuous function defined on an  interval with no zeros . Thus, $h$ cannot change sign. So  $ \im F(\zeta_{2r})$ cannot change sign infinitely often as $r \rightarrow \infty$.
\qed

\begin{propo} \label{2indep}
$(0,0,0)$ and $(0,0,2)$ are linearly independent.
\end{propo}

\noindent {\em Proof.}  
Suppose that $(0,0,0)$ and $(0,0,2)$ are linearly dependent.  Then we have 
$g_0$ and $g_2$ in $\cal{R}$ such that 
$$g_0(A) (0,0,0) + g_2(A) (0,0,2) = 0 \in S(M)$$

Hence, there are elementary  Kauffman relations  among framed links, $R_i$, 
 such that 
$$g_0(A) (0,0,0) + g_2(A)(0,0,2) = \sum_i h_i(A)R_i$$ in the free ${\cal R}$-module generated by all isotopy classes of framed links in $M$. Here, all the $h_i(A)$ are in ${\cal R}$.
Taking quantum invariants on both sides, we have, for $r$ sufficiently large,
$$g_0(A_r)I_r(M) + g_2(A_r) I_r(M,(0,0,2)) = 0.$$ 

Since $A_r^4 I_r(M,(0,0,2)) = I_r(M) - 1$ for odd $r>1$, 
$$g_0(A_r) I_r(M) +  A_r^{-4} g_2(A_r) (I_r(M) - 1) = 0,$$ for $r$ odd and sufficiently large.  

Hence, 
$$I_r(M) = \frac{g_2(A_r)}{A_r^4 g_0(A_r) + g_2(A_r)},$$ for $r$ odd and sufficiently large.

Thus, it suffices to show that there is no $f(A)$ in $Q(A)$ such that $f(\zeta_{2r})= I_r(M)$ for $r$ odd and sufficiently large.

By Proposition \ref{propo:invariant}, $(1-A_r^4) I_r(M) = \sum_{k=1}^{r-1} (-1)^k A_r^{2k^2+2k}$.

\begin{nota}
$g_N(m) = \frac{2}{\sqrt{N}} \sum_{k=0}^m \zeta_N^{k^2}$.
\end{nota}

The following lemma was provided by Paul van Wamelen.  See \cite[Appendix C]{Ha03}.

\begin{lem}{(van Wamelen)}
\label{lem:wamelen}
$$\sum_{k=1}^{r-1} \zeta_{2r}^{rk} \zeta_{2r}^{2k^2+2k} + 1 =  \zeta_{16r}^{-(r+2)^2} 
(\zeta_{16r}^{9r^2} - 2\sqrt{r}(2g_{16r}(r-1)-g_{4r}(\frac{r-1}{2}))).$$
\end{lem}

Hence, 
$(1-A_r^4) I_r(M) + 1 = \zeta_{16r}^{-(r+2)^2} 
(\zeta_{16r}^{9r^2} - 2\sqrt{r}(2g_{16r}(r-1)-g_{4r}(\frac{r-1}{2})))$.

Now, we can apply a result in \cite{Le76} by D.H. Lehmer to estimate each of the new sums.

\begin{thm}{(Lehmer)}
\label{thm:lehmer}
For $N \geq 100$ and $\sqrt{\frac{N}{2}} \leq m \leq \frac{N}{4}$, $g_N(m)$ lies within the circle with center $(h,k) = (C(\sqrt(2), S(\sqrt{2})-\frac{1}{\sqrt{2}\pi}) \approx (0.529,0.489)$ 
and radius $\frac{1}{\sqrt{2}\pi} + \frac{101}{40\sqrt{N}}$, where $C(u) = \int_0^u \cos(\frac{1}{2} \pi x^2) dx$ and $S(u) = \int_0^u \sin(\frac{1}{2} \pi x^2) dx$ are the Fresnel integrals.
\end{thm}

The proof of this result in Lehmer's paper is only sketched, and moreover, there are some minor errors.  However, Litherland has given a detailed proof \cite{Li02} along the lines indicated by Lehmer.  

\begin{figure}
\begin{minipage}{5in}\includegraphics[width=5in]{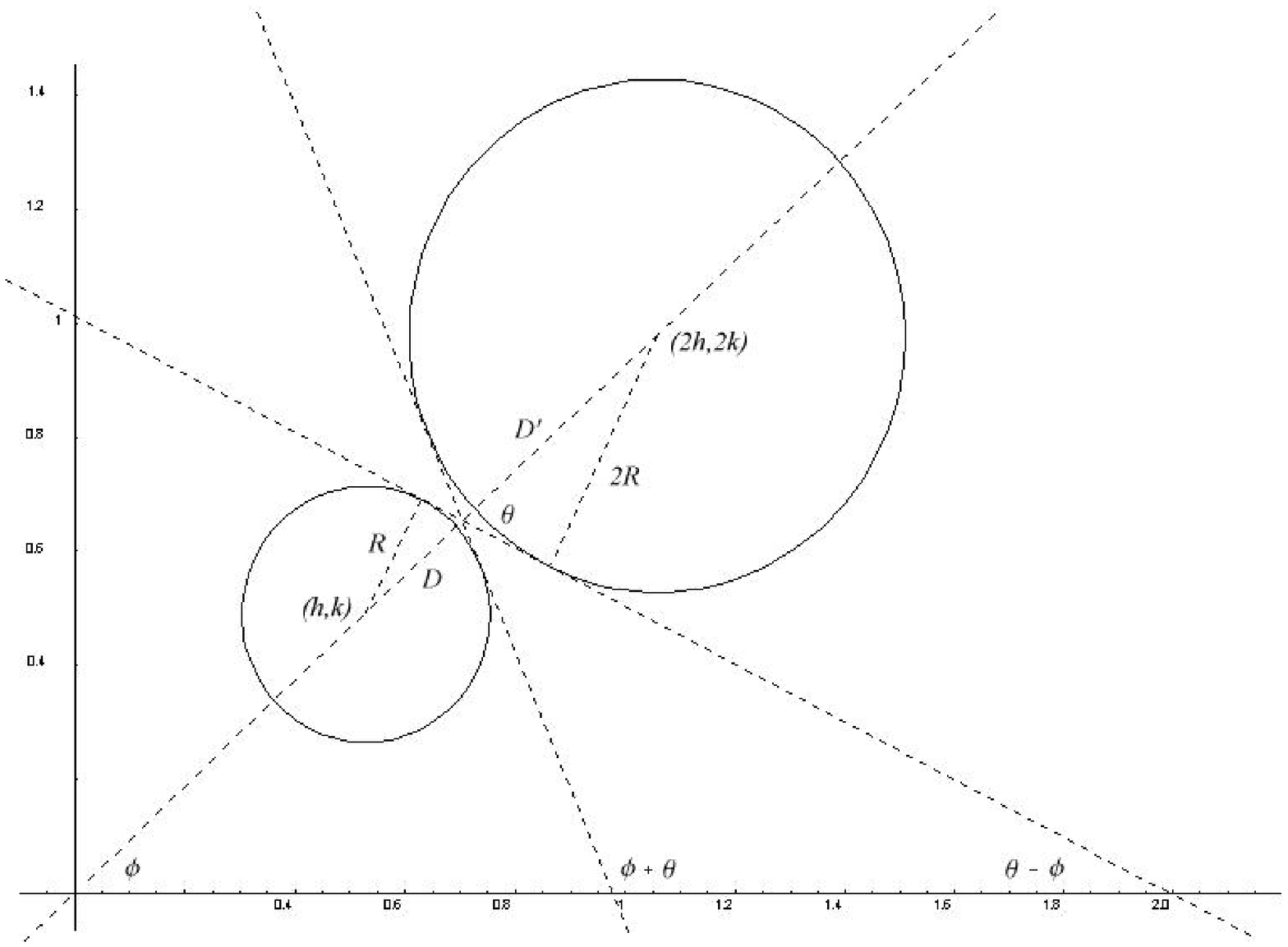}
\end{minipage}
\caption{Estimating $g_{4r}(\frac{r-1}{2})$ and $2g_{16r}(r-1)$}
\label{fig:circles}
\end{figure}

Let $R = \frac{1}{\sqrt{2}\pi} + 0.0001$.  Then, for sufficiently large $r$, 
$g_{4r}(\frac{r-1}{2})$ lies inside a circle with center $(h,k)$ and radius $R$, and $2g_{16r}(r-1)$ lies inside a circle with center $(2h, 2k)$ and radius $2R$.

Let $D$ be the distance from $(h,k)$ to the intersection of the tangent lines between the circles depicted in Figure \ref{fig:circles}, let $D'$ be the distance from $(2h,2k)$ to the point of intersection, let $L$ be the distance between the two centers, let $\theta$ be the angle between one of the tangent lines and the line joining the centers, and let $\phi$ be the angle between the line joining the centers and the $x$-axis.

$$\phi - \theta < \mathrm{ arg}(2g_{16r}(r-1)-g_{4r}(\frac{r-1}{2})) < \phi + \theta,$$ where $\mathrm{ arg}$ takes values in $[-90^\mathrm{o}, 270^\mathrm{o})$.

We will now show that multiplication by $\zeta_{16r}^{-(r+2)^2}$ rotates the difference into the upper half-plane for certain values of $r$ and into the lower half-plane for other values of $r$.
 
Since $\frac{R}{D} = \sin(\theta) = \frac{2R}{D'}$, $D' = 2D$.
$L = D + D'= 3D$, so 
$\sin(\theta) = \frac{3R}{3D} = \frac{3R}{L} = \frac{3R}{\sqrt{h^2+k^2}}$, and so $\theta = \sin^{-1} (\frac{3R}{\sqrt{h^2+k^2}}) \approx 69.7078^\mathrm{o}$.

Also, $\tan(\phi) = \frac{k}{h}$, so $\phi \approx 42.7495^\mathrm{o}$.

As $r \rightarrow \infty$, $\zeta_{16r}^{-(r+2)^2} \rightarrow \zeta_{16}^3$, for $r \equiv 9 \mmod16$ and 
$\zeta_{16r}^{-(r+2)^2} \rightarrow \zeta_{16}^{11}$, for $r \equiv 1 \mmod16$.

$$0^\mathrm{o} < 40.541^\mathrm{o} < \phi - \theta + 67.5^\mathrm{o} <  \mathrm{ arg}(2g_{16r}(r-1)-g_{4r}(\frac{r-1}{2})) + 67.5^\mathrm{o}$$
and
$$\mathrm{ arg}(2g_{16r}(r-1)-g_{4r}(\frac{r-1}{2})) + 67.5^\mathrm{o} < \phi + \theta + 67.5^\mathrm{o} < 179.96^\mathrm{o} < 180^\mathrm{o}.$$ 

Hence, for sufficiently large $r$, multiplication by $\zeta_{16r}^{-(r+2)^2}$ rotates $(2g_{16r}(r-1) - g_{4r}(\frac{r-1}{2}))$ into the upper half-plane for $r \equiv 9 \mmod16$ and into the lower half-plane for $r \equiv 1 \mmod16$.

Hence, for sufficiently large $r$, the imaginary part of 
$$\frac{(1-A_r^4) I_r(M) + 1}{2\sqrt{r}} = \zeta_{16r}^{-(r+2)^2} 
(\frac{\zeta_{16r}^{9r^2}}{2\sqrt{r}} - (2g_{16r}(r-1)-g_{4r}(\frac{r-1}{2})))$$ 
is positive for $r \equiv 1 \mmod16$ and is negative for $r \equiv 9 \mmod16$, and the same must hold for $(1-A_r^4) I_r(M) + 1$.

Hence, the imaginary part of a function $F$ such that $F(\zeta_{2r})=(1-\zeta_{2r}^4) I_r(M) + 1$ for all but finitely many odd $r$ changes sign infinitely often, and so, by Lemma \ref{lem:rational}, $F$ cannot be rational.  Hence, there can be no rational function $f$ such that $f(\zeta_{2r}) = I_r(M)$ for $r$ odd and sufficiently large, as required.  \qed

\begin{figure}
\includegraphics[width=3in]{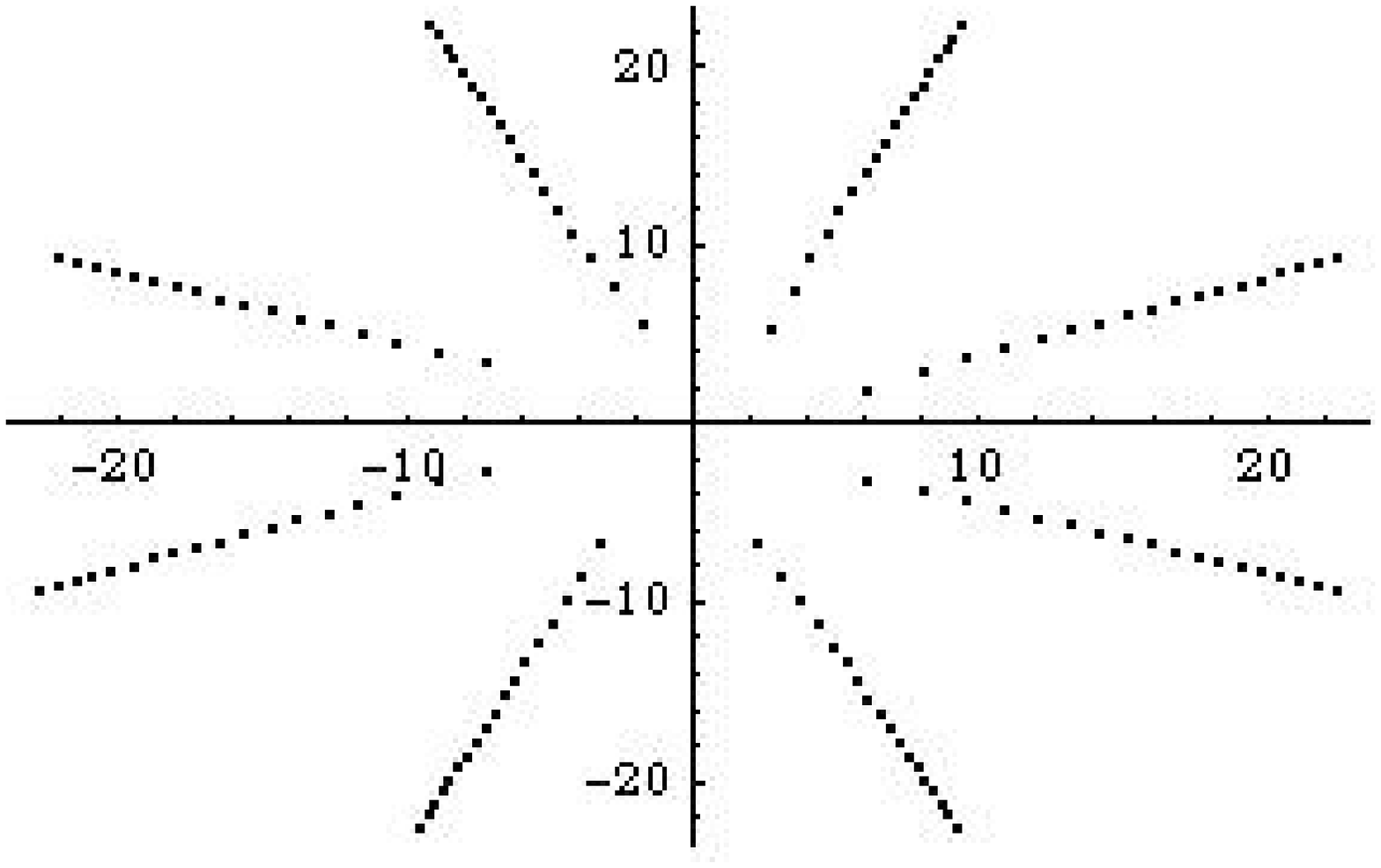}
\caption{Values of $(1-A_r^4) I_r(M)$ for $17 \leq r \leq 301$,  and $r$ odd.}
\label{fig:mod16}
\end{figure}

\begin{propo}
$(0,0,1) \neq 0$
\end{propo}

\noindent {\em Proof.}  If $(0,0,1)=0$, then there exists a linear combination of Kauffman bracket skein relators $R_i$ such that $(0,0,1) =
 \sum_i a_i({\cal A}) R_i$ 
in the free module over ${\cal R}$, and hence as in the previous proof, there is a  product of cyclotomic polynomials $k(A)$ such that $k(A_r) I_r(M,(0,0,1)) = 0$ for all $r$.

As shown in Proposition \ref{propo:invariant}, $I_r(M,(0,0,1)) = (-1)^\frac{r-1}{2} \frac{A^{-2}}{A^2+1} \neq 0$ for all odd $r>1$, and so, no such $k(A)$ can exist. 

Hence, $(0,0,1) \neq 0$. \qed

Thus, $(0,0,0)$, $(0,0,1)$, $(1,0,0)$, $(1,0,1)$, and $(0,0,2)$ are all linearly independent in
$S(M)$, and so, $\mathrm S(M) =  {\cal R}^5$.  \qed

Experiments run in Mathematica suggest that, restricted to each odd  congruence class of $r \mmod 16$, $I_r(M)$ takes on the values of some well-behaved function of $r$:  see Figure \ref{fig:mod16}.

\end{document}